\documentclass[a4paper,fleqn]{cas-sc}

\usepackage[authoryear,longnamesfirst]{natbib}

\def\tsc#1{\csdef{#1}{\textsc{\lowercase{#1}}\xspace}}
\tsc{WGM}
\tsc{QE}

\newtheorem{Proposition}{Proposition}
\usepackage{amssymb,amsmath}
\usepackage{cases}
\usepackage{float}
\usepackage{diagbox}
\usepackage{url}
\usepackage{graphicx}
\usepackage{caption}
\usepackage[left]{lineno}
\usepackage{xcolor}

\begin{document}
\let\WriteBookmarks\relax
\def\floatpagepagefraction{1}
\def\textpagefraction{.001}

\shorttitle{Curbside Congestion under Two-Tandem Bottlenecks}

\shortauthors{Y. Deng, Z. Li, S. Qian. and W. Ma.}  

\title [mode = title]{Modeling the Curbside Congestion Effects of Ride-hailing Services for Morning Commute using Bi-modal Two-Tandem Bottlenecks}  



%

\author[1,3]{Yao Deng}



\ead{hustdengyao@gmail.com}



\affiliation[1]{organization={School of Management, Huazhong University of Science and Technology},
            city={Wuhan},
            postcode={430074},
            country={China}}

\author[1]{Zhi-Chun Li}


\ead{smzcli@gmail.com}



\author[2]{Sean Qian}


\ead{seanqian@cmu.edu}



\affiliation[2]{organization={Department of Civil and Environmental Engineering, Carnegie Mellon University},
            city={Pittsburgh},
            postcode={15213}, 
            country={USA}}
     
\author[3]{Wei Ma}
\cormark[1]

\ead{wei.w.ma@polyu.edu.hk}



\affiliation[3]{organization={Department of Civil and Environmental Engineering, The Hong Kong Polytechnic University},
            city={Hong Kong SAR},
            postcode={999077}, 
            country={China}}
            
\cortext[1]{Corresponding author}



\begin{abstract}
With the proliferation of ride-hailing services, curb space in urban areas has become highly congested due to the massive passenger pick-ups and drop-offs. Particularly during peak hours, the massive ride-hailing vehicles waiting to drop off obstruct curb spaces and even disrupt the flow of mainline traffic. However, there is a lack of an analytical model that formulates and mitigates the congestion effects of ride-hailing drop-offs in curb spaces. To address this issue, this paper proposes a novel bi-modal two-tandem bottleneck model to depict the commuting behaviors of private vehicles (PVs) and ride-hailing vehicles (RVs) during the morning peak in a linear city. In the model, the upstream bottleneck models the congestion on highways, and the downstream curbside bottlenecks depict the congestion caused by RV drop-offs in curb spaces, PV queue on main roads, and the spillover effects between them in the urban area. The proposed model can be solved in a closed form under eight different scenarios. A time-varying optimal congestion pricing scheme, combined curbside pricing and parking pricing, is proposed to achieve the social optimum. It is found that potential waste of road capacity could occur when there is a mismatch between the highway and curbside bottlenecks, and hence the optimal pricing should be determined in a coordinated manner. A real-world case from Hong Kong shows that the limited curb space and main road in the urban area could be the major congestion bottleneck. Expanding the capacity of the curb space or the main road in the urban area, rather than the highway bottleneck, can effectively reduce social costs. This paper highlights the critical role of curbside management and provides policy implications for the coordinated management of highways and curb spaces.
\end{abstract}


\begin{highlights}
\item A novel bi-modal two-tandem bottleneck model is proposed and solved in closed-form
\item Novel insights into spatial-temporal congestion effects of ride-hailing drop-offs are provided
\item An optimal congestion pricing scheme addressing both curb space and highway bottlenecks is proposed
\item A real-world case study identifies the potential scenario in Hong Kong
\end{highlights}

\maketitle
\section{Introduction}\label{intro}
Curb space is becoming one of the busiest areas due to the rapid growth of various traffic demands such as car parking, passenger pick-up and drop-offs (PUDOs), truck loading, and bus stops \citep{mitman2018,castiglione2018,agarwal2023M&SOM}. Especially in view of the emergence of ride-hailing services in recent years, curb space is used in a way that has neither been envisaged nor designed for. Emerging and rapidly expanding ride-hailing services create substantial pick-up and drop-off demands on limited curb spaces \citep{lu2019,jaller2021}. This significantly adds to the growing congestion problem in curb spaces, demonstrating the critical roles of curb spaces in transportation road networks \citep{castiglione2018, erhardt2019Sci.Adv., agarwal2023M&SOM, tirachini2020Transportation, beojone2021Transp.Res.PartCEmerg.Technol., zhang2021AvailableSSRN3974957}. 

The congestion effects of PUDOs on curb spaces manifest in two aspects. First, the increasing demand for PUDOs could saturate limited curb spaces and cause travel delays for ride-hailing services. For example, in San Francisco, the average duration of PUDOs increases from around 60s in 2010 \citep{erhardt2019Sci.Adv., lu2019, jaller2021, rahaim2019} to 144.75s on major arterial and 79.49s on the minor arterial in 2016 \citep{erhardt2019Sci.Adv.,LiuXiaoHui2023TS}. Second, queuing congestion in curb spaces could spread to the main road, and disrupt mainstream traffic and even system-level network performance \citep{castiglione2018, butrina2020TransportPolicy, erhardt2019Sci.Adv., LiuJiaChao2023TRC}. A recent evaluation of the spillover effect of PUDOs by \cite{goodchild2019} revealed that one more unit of PUDO could lead to a decrease of approximately 1\% in average traffic speed for study lanes of Boren Avenue in Seattle. \citet{LiuXiaoHui2023TS} further demonstrated that, on average, 100 additional units of PUDOs would decrease traffic speed by 3.70 miles/h on weekdays and 4.54 miles/h on weekends in regions below West 110th Street in the Manhattan area. Furthermore, it is envisaged that the competition for limited curb space among the enormous rapidly growing demand of ride-hailing PUDOs will become more severe, especially with the further development of Mobility-as-a-Service (MaaS) \citep{smith2019}. Therefore, it is essential to understand and analyze the congestion effects of ride-hailing PUDOs on both curb space and system-level network performance.

Despite the increasingly severe congestion issues arising from ride-hailing PUDOs in curb spaces, research on curb space as a significant area for ride-hailing PUDOs is relatively scarce. Existing curb-related literature still primarily views the curb as a space for on-street parking \citep{yousif2004ICEProc.Transp., wijayaratna2015, cao2015Transp.Res.Rec., chen2016Transp.Res.Rec., lei2017Transp.Res.PartCEmerg.Technol., pu2017Transp.Res.PartCEmerg.Technol., dowling2020IEEETrans.Intell.Transp.Syst., abhishek2021Transp.Sci.}, bicycle lanes \citep{schimek2018Accid.Anal.Prev., ye2018Transp.Res.Rec.}, commercial vehicles loading/unloading \citep{jones2009ITEJ., chiara2020Transp.Policy}, and bus stop \citep{patkar2020Arab.J.Sci.Eng.}. Furthermore, these existing studies predominantly concentrate on the curb space itself, overlooking the potential congestion issues and the spillover effects on main roads and even the entire transportation system \citep{LiuJiaChao2023TRC}. Investigating the congestion caused by ride-hailing services in curb spaces and its possible spillover effects is crucial. Therefore, it is imperative to address the use of curb spaces by ride-hailing services. 

Although a limited number of studies have addressed curbside congestion and its impact on traffic flow, they mainly focus on on-street parking that involves a single traffic mode. For example, \cite{arnott2006J.UrbanEcon., arnott2010J.UrbanEcon.} and \cite{arnott2013Transp.Res.PartPolicyPract.} integrated curbside parking (also known as on-street parking) and traffic flow congestion into a downtown parking model. \cite{arnott2009Reg.Sci.UrbanEcon.} and \cite{arnott2015J.UrbanEcon.} further incorporated garage parking into the downtown parking model and examined the effectiveness of curbside parking policy. These studies depend on the steady-state spatial equilibrium of parking, overlooking how congestion influences travelers' travel mode choices and travel time schedules. Therefore, it is essential to examine the impact of curbside congestion on travelers' travel mode choices and travel time schedules in specific scenarios, especially during the morning peak period.

From the perspective of curbside management, various policy instruments have been extensively discussed in the literature, such as various charging schemes \citep{jones2009ITEJ., guo2013Transp.Policy,qian&rajagopal2013Procedia-SocialandBehavioralSciences, qian&rajagopal2014TransportationResearchPartC:EmergingTechnologies, millard-ball2014Transp.Res.PartPolicyPract., zheng2016Transp.Res.PartBMethodol., inan2019Transp.Res.PartBMethodol.}, parking time limits or permits \citep{arnott2013Transp.Res.PartPolicyPract., inan2019Transp.Res.PartBMethodol., stüger2022TransportationResearchRecordJournaloftheTransportationResearchBoard}, parking capacity regulation \citep{arnott2015J.UrbanEcon., inan2019Transp.Res.PartBMethodol., diehl2021Transp.Res.Rec.}, dedicated curb spaces for specialized use \citep{ye2018Transp.Res.Rec., shaheen2019, goodchild2019, mccormack2019Int.J.Transp.Dev.Integr., chiara2020Transp.Policy, diehl2021Transp.Res.Rec.}. In terms of practical applications, the U.S. Department of Transportation (DOT) released a curbside inventory report in 2021 \citep{abel2021}, and case studies have been conducted in Seattle and Washington D.C. \citep{goodchild2019, pochowski2022Transp.Res.Rec.}. The European Institute of Innovation and Technology (EIT) launched FlexCurb in 2022 to digitize and model supply and demand patterns in curb spaces \citep{fernandez2022EITUrbanMobility}. The Seattle Department of Transportation (SDOT) proposed and evaluated the effects of the RapidRide Roosevelt project on curbside usage in Seattle \citep{zimbabwe2018}. However, most existing research on curbside management has primarily concentrated on the curb space itself. Few studies have integrated curb space into general transportation networks that model a complete trip like the morning commute process. It is important to note that various curbside management strategies can profoundly influence travelers' behavior, such as mode choices and departure time preference, leading to substantially different traffic patterns within transportation networks.

In summary, existing curb-related literature has been relatively limited in studying the impact of emerging ride-hailing services. Most of these studies have focused on the ego performance of curb spaces, without considering curb space usage within a transportation network and a complete trip process. Consequently, these studies fail to capture the effects of increasing congestion on both the curb space itself and the entire transportation system, which hinders the system-level operation and management of curb spaces.

To fill the above gaps, this paper proposes a novel bi-modal two-tandem bottleneck model to describe the congestion effects of ride-hailing services on both curb space and the entire transportation system in morning commuting traffic. Intuitively, during the morning peak, there is a higher demand for ride-hailing drop-offs in the curb spaces of the central business district (CBD) than that for pick-ups. Therefore, this paper primarily focuses on the drop-off process of ride-hailing services\footnote{In contrast, the demand of ride-hailing pick-ups at curb spaces in the CBD will increase during the evening peak commuting hours. Therefore, the pick-up and drop-off processes of ride-hailing services can be examined together when addressing morning-evening commuting dynamics.}. To be specific, the model explicitly incorporates different traffic modes: ride-hailing vehicles (RVs) and private vehicles (PVs), and it characterizes potential bottleneck congestion in both curb space and the highway simultaneously. This innovative modeling approach allows us to evaluate the consequence of congestion on both RV and PV commuters from both temporal and spatial perspectives. 

The proposed model generalizes the two-tandem bottleneck model \citep{Kuwahara1990,Kim1999}. On the one hand, the model originates from real-world scenarios, specifically the mixed traffic flow during the morning peak commute involving private vehicles and ride-hailing vehicles. It incorporates shared road usage by multi-modal traffic and the congestion spillover phenomenon, enriching the application scenarios of the two-tandem bottleneck model.
On the other hand, in contrast to the previous studies that focused on merging Y-shaped bottlenecks, where multiple streams (usually with single traffic mode) merge at the downstream bottleneck \citep{Kuwahara1990,arnott1993TS,Kim1999,akamatsu2015TRP,he2022TS,wang2022TRE}, our model tackles a diverging Y-shaped bottleneck, where a multi-modal stream diverges at the downstream bottlenecks in curb space. Furthermore, our model captures the interactions among diverged streams by considering the congestion spillover effects, a factor that has been overlooked in previous studies. In addition, our model can be solved in closed-form under eight different scenarios, allowing us to develop curbside management strategies to effectively address traffic congestion and spillover effects. Therefore, an optimal congestion pricing mechanism that can eliminate the negative externalities is proposed. 
Furthermore, the performance of congestion pricing and capacity regulation is examined numerically in a real-world case from Hong Kong. 

The studies most closely related to this paper are \cite{LiuXiaoHui2023TS} and \cite{LiuJiaChao2023TRC}. The former study estimates the congestion effect of PUDOs on general networks from a data-driven perspective. However, their estimation is conducted at the regional level rather than the road level. As a result, it was not able to accurately characterize users' travel behavior such as travel time schedule. Furthermore, \citet{LiuJiaChao2023TRC} uses a static M/M/1 curbside queuing model to quantify the waiting time and queue length caused by curbside stops, and a bi-modal traffic assignment model is formulated to describe user travel behaviors and curbside usage. However, it also overlooks users' travel time schedule, which prevents it from characterizing the temporal distribution of congestion effects.

In summary, this paper makes the following contributions:
\begin{itemize}
\item To the best of our knowledge, this is the first study using tractable bottleneck models to describe the congestion effect of ride-hailing drop-offs in mode and departure time choices. This approach enables us to comprehensively discuss the congestion effects of ride-hailing services on both the curb space and the main road from temporal and spatial perspectives, as well as the further implementation of curbside management strategies.

\item This study introduces a novel bi-modal two-tandem bottleneck model, modeling the interaction between multi-modal traffic by considering the congestion spillover effects on a diverging Y-shaped bottleneck, which distinguishes it from existing two-tandem and merging Y-shaped bottleneck models.

\item An optimal congestion pricing scheme is proposed to eliminate congestion externalities within the system, achieving coordinated traffic management on highways and curb spaces. The effectiveness of congestion pricing and bottleneck expansion is examined using a real-world case study from Hong Kong.
\end{itemize}

The remainder of this paper is organized as follows. Section \ref{model} details the formulations of the bi-modal two-tandem bottleneck model and the user equilibrium conditions. Section \ref{section_different_scenarios} and \ref{section_bieffects} derive closed-form solutions for various equilibrium scenarios under unidirectional and bidirectional congestion spillover effects. Section \ref{section_congestion_pricing} introduces the optimal congestion pricing scheme. Section \ref{section_numerical_illustration} presents the numerical results to illustrate the properties of the proposed model and the findings from a real-world case study. Section \ref{section_LateArrival} examines the effects of relaxing the no-late-arrival assumption on user equilibrium results and the congestion pricing scheme. Finally, conclusions and future research directions are drawn in Section \ref{section_conclusion}.

\section{Formulations} \label{model}
In this section, we first present the basic model settings, which include configurations and attributes of the two-tandem bottlenecks corridor, such as bottleneck capacities and travel times. Next, we define the generalized travel costs for commuters using ride-hailing services and private vehicles. Lastly, we elaborate on the user equilibrium conditions in morning commuting.

\subsection{Basic Settings} \label{basicsetting}
\begin{figure}[!ht]
  \centering
  \includegraphics[width=0.8\linewidth]{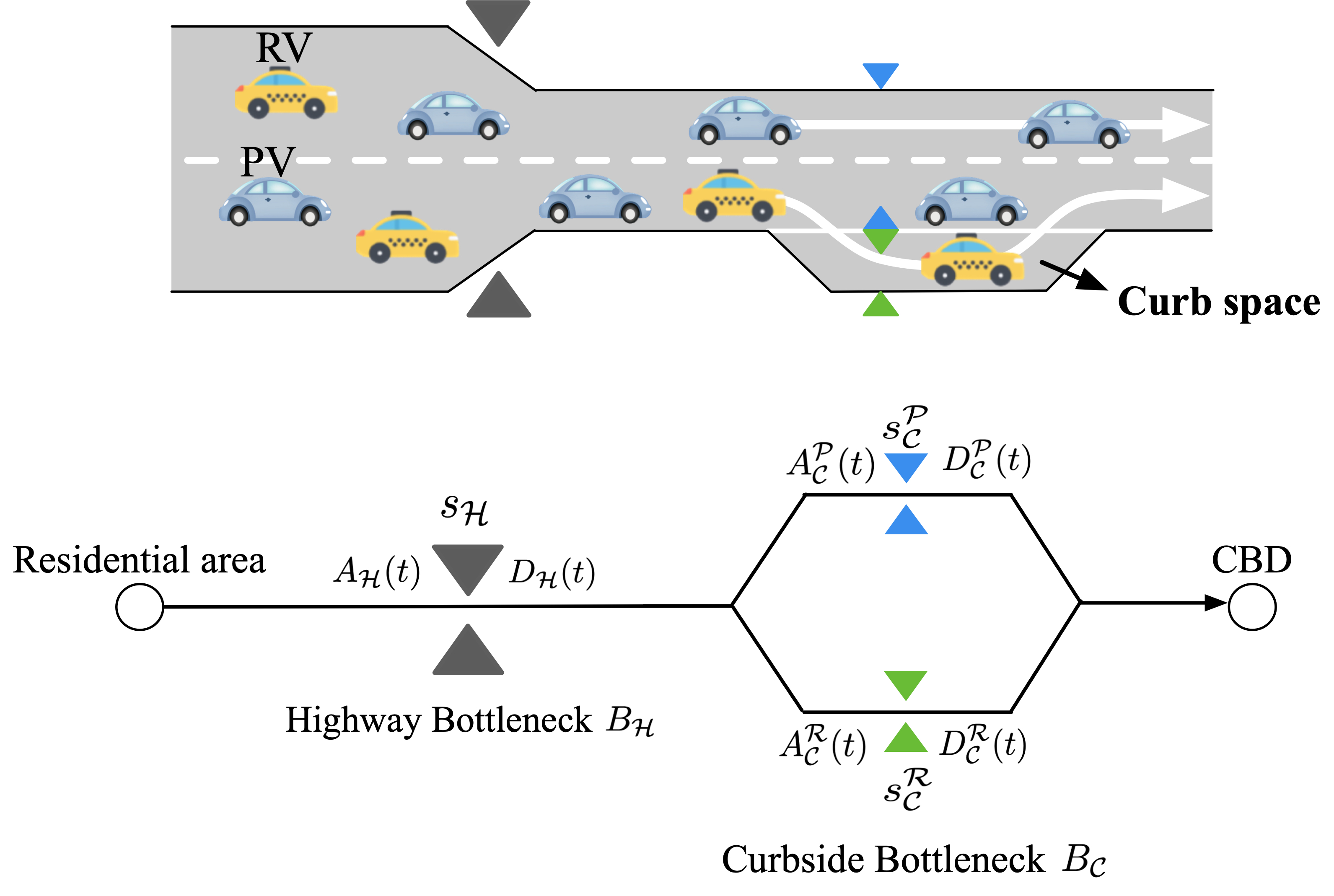}
  \caption{Transport corridor with two-tandem bottlenecks}
  \label{overview of the road network}
\end{figure}

We consider a transport corridor with bottleneck roads connecting the residential area (or home) and the CBD. Every morning, $N$ homogeneous commuters who have an identically preferred arrival time $t^*$ (e.g., 8:00 AM) at the CBD workplace depart from home and head to the worplace. Commuters can choose to drive alone or use ride-hailing services. Thus, there are two types of vehicles in the system: private vehicles (PVs) and ride-hailing vehicles (RVs). An overview of the transport corridor is presented in Fig. \ref{overview of the road network}. It consists of two tandem bottlenecks that model congestion delays on the highway and on the curbside, respectively. The highway bottleneck $B_\mathcal{H}$ with a service rate (or capacity) of $s_\mathcal{H}$, which is near the residential area, models the queuing delay on the highway caused by all types of vehicles (RVs and PVs); while the curbside bottleneck $B_\mathcal{C}$, which is near the CBD, models the queuing delay of RVs and PVs, as well as the delay caused by the spillover effects. Specifically, the curbside bottleneck can be divided into two bottlenecks for RVs and PVs, and their corresponding service rates are $s_\mathcal{C}^\mathcal{R}$ and $s_\mathcal{C}^\mathcal{P}$, respectively. 
Intuitively, the bottleneck for PVs is due to the limited capacity of main roads, while the bottleneck for RVs is due to the restricted curb spaces for drop-offs\footnote{Designating specific curb space for ride-hailing pick-up and drop-off has been increasingly implemented in various cities. For example, the city of Seattle has established Passenger Loading Zone (PLZ) spaces as designated pick-up/drop-off locations. Transportation Network Companies (TNCs) like Uber and Lyft have implemented geofences to direct drivers and passengers to these specific locations on a block \citep{goodchild2019}. The city of Brisbane has also designated spaces for passenger pick-up and drop-off with a two-minute time limit. As of April 2, 2024, the Hong Kong SAR has set up 88 taxi stands and 37 taxi pick-up/drop-off points for passenger use.}.
Since the curb space and local road capacities are typically lower than the highway capacity, we assume $s_\mathcal{C}^\mathcal{R}<s_\mathcal{H}$ and $s_\mathcal{C}^\mathcal{P}<s_\mathcal{H}$ in this paper.

Without loss of generality, we regard that the free flow travel time from the residential area to the CBD is zero, thus the travel time of commuters is characterized by the queuing delays experienced at the highway and curbside bottlenecks.
We define the cumulative arrival curves of RVs and PVs at the highway bottleneck (i.e., the cumulative departure curves from home) as $A_\mathcal{H}^\mathcal{R}(t)$ and $A_\mathcal{H}^\mathcal{P}(t)$, and the cumulative departure curves from the highway bottleneck as $D_\mathcal{H}^\mathcal{R}(t)$ and $D_\mathcal{H}^\mathcal{P}(t)$, respectively.
Thus, $A_\mathcal{H}(t)=A_\mathcal{H}^\mathcal{R}(t)+A_\mathcal{H}^\mathcal{P}(t)$ and $D_\mathcal{H}(t)=D_\mathcal{H}^\mathcal{R}(t)+D_\mathcal{H}^\mathcal{P}(t)$ represent the cumulative curves of all types of vehicles that arrive and depart from the highway bottleneck, respectively.
Based on the first-in-first-out (FIFO) principle, the queuing delays for RVs and PVs arriving at the highway bottleneck at time $t$ are the same, as expressed in Eq. \eqref{eq:w2}.
\begin{equation}\label{eq:w2}
    w_\mathcal{H}(t)=w_\mathcal{H}^\mathcal{R}(t)=w_\mathcal{H}^\mathcal{P}(t)=\frac{A_\mathcal{H}(t)-D_\mathcal{H}(t)}{s_\mathcal{H}},
\end{equation}
where the slopes of cumulative curves, $\dot{D}_\mathcal{H}(t)$ and $\dot{A}_\mathcal{H}(t)$, meet:
\begin{equation}\label{eq:rate2}
    \dot{D}_\mathcal{H}(t)=
        \begin{cases} 
            \dot{A}_\mathcal{H}(t),\ &\text{if}\  w_\mathcal{H}(t)=0 \ \text{and}\ \dot{A}_\mathcal{H}(t)\leq s_\mathcal{H}, \\
            s_\mathcal{H},\ &\text{else}.
        \end{cases}
\end{equation}

Different from the highway bottleneck, the curbside bottlenecks model the queuing delay of RVs and PVs, as well as the delay caused by the congestion spillover effects. We define the cumulative curves of RVs that arrive at and leave the curbside bottleneck as $A_\mathcal{C}^\mathcal{R}(t)$ and $D_\mathcal{C}^\mathcal{R}(t)$, respectively. Let the cumulative curves of PVs that arrive at and leave the curbside bottleneck be $A_\mathcal{C}^\mathcal{P}(t)$ and $D_\mathcal{C}^\mathcal{P}(t)$, respectively. Due to the assumption of zero free flow travel time, the cumulative arrival curve at the curbside bottlenecks is the cumulative departure curve from the highway bottleneck, as presented in Eq. \eqref{eq:A1D2}.
\begin{equation}\label{eq:A1D2}
        \begin{cases}
            A_\mathcal{C}^\mathcal{R}(t)=D_\mathcal{H}^\mathcal{R}(t),\\
            A_\mathcal{C}^\mathcal{P}(t)=D_\mathcal{H}^\mathcal{P}(t).
        \end{cases}
\end{equation}

\subsection{Spillover effects of curbside drop-offs}
After passing through the highway bottleneck, RVs should park and drop off passengers in the curb space. This process requires RVs to slow down and change lanes to the right-most one. Additionally, the limited curb space will cause RVs to queue. In this paper, we define the time required for RVs to complete this process and pass through the curbside bottleneck as the bottleneck queuing delay. Based on the FIFO principle, the queuing delay for RVs arriving at the curbside bottleneck at time $t$ is presented by Eq. \eqref{eq:wcr}.
\begin{equation}\label{eq:wcr}
    w_\mathcal{C}^\mathcal{R}(t)=\frac{A_\mathcal{C}^\mathcal{R}(t)-D_\mathcal{C}^\mathcal{R}(t)}{s_\mathcal{C}^\mathcal{R}},
\end{equation}
where the slopes of cumulative curves, $\dot{D}_\mathcal{C}^\mathcal{R}(t)$ and $\dot{A}_\mathcal{C}^\mathcal{R}(t)$, meet: 
\begin{equation}\label{eq:rate1r}
    \dot{D}_\mathcal{C}^\mathcal{R}(t)=
    \begin{cases}
        \dot{A}_\mathcal{C}^\mathcal{R}(t),\ &\text{if}\  w_\mathcal{C}^\mathcal{R}(t)=0 \ \text{and}\ \dot{A}_\mathcal{C}^\mathcal{R}(t)\leq s_\mathcal{C}^\mathcal{R},\\
        s_\mathcal{C}^\mathcal{R},\ &\text{else}.
    \end{cases}
\end{equation}

However, the complex behaviors of RVs, such as slowing down, changing lanes, and queuing may spill over to the main road and affect PV commuters driving through the main road. This phenomenon is quite similar to the queue spillover at busy highway off-ramps. Some previous studies on off-ramp queue spillover have already highlighted its significant impact in reducing the flow rate on the main road \citep{muñoz&daganzo2002TransportationResearchPartA:PolicyandPractice,cassidy&anani2002TransportationResearchPartA:PolicyandPractice,hall&2018JournalofPublicEconomics,chen&chen2021TransportationResearchPartC:EmergingTechnologies}. Additionally, some recent studies on curbside management \citep{goodchild2019,LiuJiaChao2023TRC, LiuXiaoHui2023TS} have found that the average traffic speed on main roads is negatively related to the number of RV drop-offs. Based on these studies, we model the congestion spillover effect by assuming that the queue spillover of RVs reduces the service rate of the curbside bottleneck serving PVs (i.e., the main road's capacity)\footnote{We acknowledge that the queue spillover impact may be mutual, as congestion caused by PVs on the main road may also limit RVs’ ability to reach the curbside. Therefore, we will discuss this issue later in Section \ref{section_bieffects}.}.

For illustration purposes, we define $\widetilde{s}_\mathcal{C}^\mathcal{P}(t)$ as the actual service rate of the curbside bottleneck serving PVs, which is discounted due to the congestion spillover effects. Under the point-queue assumption \citep{vickrey1969,arnott1990TransportationResearchPartB:Methodological,arnott1990JournalofUrbanEconomicsa} and the FIFO rule, we assume that PVs are only affected by the congestion spillover effects from those RVs departing at the same time\footnote{{Under the point-queue assumption, RV commuters who are already in the queue will not impact PV commuters who arrive at the curbside bottleneck entrance afterward. Under the FIFO rule, RV and PV commuters departing at the same time will arrive at the entrance of the curbside bottlenecks simultaneously. It is intuitive to assume that the impact of RV congestion spillover on PV commuters depends on their relative arrival rates at the entrance of the curbside bottlenecks.}}. Therefore, the service rate $\widetilde{s}_\mathcal{C}^\mathcal{P}(t)$ can be defined as:
\begin{equation}\label{eq:scp_actual}
        \widetilde{s}_\mathcal{C}^\mathcal{P}(t)=\begin{cases}
        \frac{\dot{A}_\mathcal{C}^\mathcal{P}(t)}{\dot{A}_\mathcal{C}^\mathcal{P}(t)+\delta^\mathcal{R}\dot{A}_\mathcal{C}^{\mathcal{R}}(t)}s_\mathcal{C}^\mathcal{P},\ &\text{if}\ \dot{A}_\mathcal{C}^\mathcal{P}(t)>0, \dot{A}_\mathcal{C}^\mathcal{R}(t)>0, \\
        s_\mathcal{C}^\mathcal{P},\ &\text{else.}\ \end{cases}\\
\end{equation} where $\delta^\mathcal{R}$ denotes the intensity of congestion spillover effects from the RV queue on PVs, and $\delta^\mathcal{R}\in(0,1)$ holds in this paper.

Therefore, the delay for PVs arriving at curbside bottleneck at time $t$ is presented by Eq. \eqref{eq:wcp}.
\begin{equation}\label{eq:wcp}
    w_\mathcal{C}^\mathcal{P}(t)=\frac{A_\mathcal{C}^\mathcal{P}(t)-D_\mathcal{C}^\mathcal{P}(t)}{\widetilde{s}_\mathcal{C}^\mathcal{P}(t)},
\end{equation} where the slopes of the cumulative curves, $\dot{D}_\mathcal{C}^\mathcal{P}(t)$ and $\dot{A}_\mathcal{C}^\mathcal{P}(t)$, meet: 
\begin{equation}\label{eq:u1c}
    \dot{D}_\mathcal{C}^\mathcal{P}(t)=
        \begin{cases}
            \dot{A}_\mathcal{C}^\mathcal{P}(t),\ &\text{if}\  w_\mathcal{C}^\mathcal{P}(t)=0 \ \text{and}\ \dot{A}_\mathcal{C}^\mathcal{P}(t)\leq\widetilde{s}_\mathcal{C}^\mathcal{P}(t),  \\
            
            \widetilde{s}_\mathcal{C}^\mathcal{P}(t),\ &\text{else}.
        \end{cases}\\
\end{equation}

Here we have discussed the attributes of the transport corridor, especially the travel times (i.e., the queuing delays at bottlenecks) represented by Eqs. \eqref{eq:w2}, \eqref{eq:wcr}, and \eqref{eq:wcp}. Next, we will define the generalized travel costs for RV and PV commuters during the morning peak. This allows us to establish user equilibrium conditions for the morning peak commuting process.

\subsection{Generalized travel cost} \label{travelcost}
The generalized travel cost for morning commuting consists of the travel time cost, schedule delay cost (associated with deviating from the preferred arrival time), and other expenses such as the ride-hailing service charges and the fixed cost for solo driving. To simplify the formulation, this paper focuses solely on early arrival scenarios, assuming that the penalty for late arrival to the workplace is infinitely high\footnote{This assumption has been widely adopted in previous studies \citep{arnott&depalma2011TransportationResearchPartB:Methodologicala, depalma&arnott2012TransportationResearchPartB:Methodologicala,xiao2013TRB,xiao&zhang2014TransportationScience,wu&huang2015TransportationResearchPartB:Methodological}, and we follow this approach here to simplify the discussion. We will relax this assumption in Section \ref{section_LateArrival} to examine its impact on the equilibrium patterns and congestion pricing.}. Therefore, the schedule delay cost is incurred only by arriving earlier than the preferred arrival time $t^*$. In addition, we assume that the walking cost for using ride-hailing services is zero due to the door-to-door service provided. Since PV users may need to walk from the parking lot to the workplace, the walking cost for PV users is non-zero and factored into their fixed cost.

For RV commuters leaving the residential area at time $t$, the generalized travel cost $C^\mathcal{R}(t)$ can be presented by Eq. \eqref{eq:costRV}.
\begin{equation}\label{eq:costRV}
    \begin{aligned}
        C^\mathcal{R}(t)&=\underbrace{\alpha \big(w_\mathcal{H}(t)+w_\mathcal{C}^\mathcal{R}\big(t+w_\mathcal{H}(t)\big)\big)}_\text{travel time cost}+\underbrace{\beta \big(t^*-t-w_\mathcal{H}(t)-w_\mathcal{C}^\mathcal{R}\big(t+w_\mathcal{H}(t)\big)\big)}_\text{schedule delay cost}\\
        &+\underbrace{\big(\lambda L+\pi \big(w_\mathcal{H}(t)+w_\mathcal{C}^\mathcal{R}\big(t+w_\mathcal{H}(t)\big)\big)+u^\mathcal{R}\big)}_\text{RV charges}\\
        &=(\alpha+\pi-\beta) \big(w_\mathcal{H}(t)+w_\mathcal{C}^\mathcal{R}\big(t+w_\mathcal{H}(t)\big)\big)+ \beta (t^*-t)+c^\mathcal{R}.
    \end{aligned}
\end{equation}

The first term on RHS (right-hand side) of Eq. \eqref{eq:costRV} represents the travel time cost, parameter $\alpha$ denotes the value of travel time. The second term is the schedule delay cost, and the parameter $\beta$ denotes the value of early arrival time. To ensure the existence of user equilibrium, it is assumed that $\beta < \alpha$ in this paper. This assumption has been validated by the empirical study \citep{small1982Am.Econ.Rev.} and widely adopted in previous studies \citep{zhang&huang2008TransportationResearchPartB:Methodological, xiao2013TRB,li2020TransportationResearchPartB:Methodological,xiao2021TransportationResearchPartB:Methodological}. The third term on RHS is charges using ride-hailing services, where $L$ is the distance from the residential area to the CBD drop-off area, $\lambda$ and $\pi$ denote the unit price for trip distance and time, respectively. $u^\mathcal{R}$ represents other related expenses, such as the flag-down fee. For simplicity, we use parameter $c^\mathcal{R}=\lambda L+u^\mathcal{R}$ to represent the fixed cost of RV, which is independent of departure time $t$.

For PV commuters, the generalized travel cost $C^\mathcal{P}(t)$ is presented by Eq. \eqref{eq:costPV}.
\begin{equation}\label{eq:costPV}
\begin{aligned}
    C^\mathcal{P}(t)&=\underbrace{\alpha \big(w_\mathcal{H}(t)+w_\mathcal{C}^\mathcal{P}\big(t+w_\mathcal{H}(t)\big)\big)}_\text{travel time cost} + \underbrace{\beta \big(t^*-t-w_\mathcal{H}(t)-w_\mathcal{C}^\mathcal{P} \big( t+w_\mathcal{H}(t) \big) \big)}_\text{schedule delay cost}+u^\mathcal{P}\\
    &=(\alpha-\beta)\big(w_\mathcal{H}(t)+w_\mathcal{C}^\mathcal{P}\big(t+w_\mathcal{H}(t)\big)\big)+\beta(t^*-t)+u^\mathcal{P},
\end{aligned}
\end{equation}
where the first two terms on RHS are costs of travel time and schedule delay. $u^\mathcal{P}$ represents other related expenses such as purchasing costs, maintenance expenses, vehicle depreciation, parking fee, and etc., collectively referred to as the fixed cost of PV.

\subsection{User equilibrium} \label{UEconditions}
In the context of bi-modal user equilibrium (UE), commuters choose either RV or PV based on the option with a lower generalized travel cost and determine the departure times from home based on a trade-off between the travel time cost and the schedule delay cost. Therefore, the intra-equilibrium of each mode concerning departure times and the inter-equilibrium between two modes must be considered simultaneously. 

The intra-equilibrium conditions ensure that the generalized travel costs for RVs and PVs remain constant with respect to the departure time $t$, as presented by Eq. \eqref{eq:innerUE}.
\begin{equation}\label{eq:innerUE}
\begin{aligned}
&\textbf{Intra-equilibrium conditions:}\\
    &\begin{cases}
        \frac{\partial{C^\mathcal{R}(t)}}{\partial t}=(\alpha+\pi-\beta)\big(\dot{w}_\mathcal{H}(t)+\dot{w}_\mathcal{C}^\mathcal{R}\big(t+w_\mathcal{H}(t)\big)\big(1+\dot{w}_\mathcal{H}(t)\big)\big)-\beta = 0, &\forall t\in(t_0^\mathcal{R},t_1^\mathcal{R}),\\[1.5ex]
        
        \frac{\partial{C^\mathcal{P}(t)}}{\partial t}=(\alpha-\beta)\big(\dot{w}_\mathcal{H}(t)+\dot{w}_\mathcal{C}^\mathcal{P}\big(t+w_\mathcal{H}(t)\big)\big(1+\dot{w}_\mathcal{H}(t)\big)\big)-\beta = 0,&\forall t\in(t_0^\mathcal{P},t_1^\mathcal{P}),
    \end{cases}
\end{aligned}
\end{equation}
where $t_0^\mathcal{R}$ and $t_1^\mathcal{R}$ denote the departure times from home of the first and last RV commuters, and $t_0^\mathcal{P}$ and $t_1^\mathcal{P}$ are the departure times from home of the first and last PV commuters, respectively.

The inter-equilibrium requires the equalization of generalized travel costs for RVs and PVs if both modes are utilized. Thus, the inter-equilibrium condition can be represented by Eq. \eqref{eq:costDiff}.
\begin{equation}
\label{eq:costDiff}
    \begin{aligned}
        C^\mathcal{R}(t_0^\mathcal{R})-C^\mathcal{P}(t_0^\mathcal{P})&=(\alpha+\pi-\beta)\big(w_\mathcal{H}(t_0^\mathcal{R})+w_\mathcal{C}^\mathcal{R}\big(t_0^\mathcal{R}+w_\mathcal{H}(t_0^\mathcal{R})\big)-(\alpha-\beta)\big(w_\mathcal{H}(t_0^\mathcal{P})+w_\mathcal{C}^\mathcal{P}\big(t_0^\mathcal{P}+w_\mathcal{H}(t_0^\mathcal{P})\big)\big)\\
        &+\beta(t_0^\mathcal{P}-t_0^\mathcal{R})-(u^\mathcal{P}-c^\mathcal{R})\\
        &=0.
    \end{aligned}
\end{equation}

At the beginning of the morning peak, there is no queue in the system, i.e., $w_\mathcal{H}(\cdot)=w_\mathcal{C}^\mathcal{R}(\cdot)=w_\mathcal{C}^\mathcal{P}(\cdot)=0$. By substituting this into Eq. \eqref{eq:costDiff}, it can be proven that in the initial phase of the morning peak, only one mode will be utilized if $ u^\mathcal{P}\neq c^\mathcal{R}$ holds. Due to the higher costs associated with private car ownership, including vehicle purchase, parking fees, and additional walking cost, we assume $u^\mathcal{P}>c^\mathcal{R}$ in this paper. Therefore, only RV will be utilized at the beginning of the morning peak and the condition of $t_0^\mathcal{R}<t_0^\mathcal{P}$ holds in our model. By combining the intra-equilibrium conditions of Eq. \eqref{eq:innerUE} with Eqs. \eqref{eq:w2}-\eqref{eq:rate1r}, we can derive the equilibrium departure rate from home $\dot{A}_\mathcal{H}(t)$ (i.e., arrival rate at the highway bottleneck) and the arrival rate at the curbside bottleneck $\dot{A}_\mathcal{C}^\mathcal{R}(t)$ for the beginning phase $(t_0^\mathcal{R},t_0^\mathcal{P})$ as,
\begin{equation}\label{eq:firstphase}
    \begin{cases}
        \dot{A}_\mathcal{H}(t)=\dot{A}_\mathcal{C}^\mathcal{R}(t)=\frac{(\alpha+\pi) s_\mathcal{C}^\mathcal{R}}{\alpha+\pi-\beta}, &\text{if}\  \frac{s_\mathcal{C}^\mathcal{R}}{s_\mathcal{H} }\leq \frac{\alpha+\pi-\beta}{\alpha+\pi},\\[1.5ex]
        \dot{A}_\mathcal{H}(t)=\frac{(\alpha+\pi) s_\mathcal{C}^\mathcal{R}}{\alpha+\pi-\beta},\dot{A}_\mathcal{C}^\mathcal{R}(t)=s_\mathcal{H}, &\text{if}\  \frac{s_\mathcal{C}^\mathcal{R}}{s_\mathcal{H} } > \frac{\alpha+\pi-\beta}{\alpha+\pi}.\\
    \end{cases}
\end{equation}

Eq. \eqref{eq:firstphase} indicates the arrival rate of commuters at bottlenecks and the resulting queuing scenarios during the period of  $(t_0^\mathcal{R},t_0^\mathcal{P})$. To be specific, as $\frac{s_\mathcal{C}^\mathcal{R}}{s_\mathcal{H} }\leq \frac{\alpha+\pi-\beta}{\alpha+\pi}$, $\dot{A}_\mathcal{H}(t)\leq s_\mathcal{H}$ and $\dot{A}_\mathcal{C}^\mathcal{R}(t)>s_\mathcal{C}^\mathcal{R}$ hold, RV queue forms exclusively at the curbside bottleneck during the period of $(t_0^\mathcal{R},t_0^\mathcal{P})$. Conversely, as $\frac{s_\mathcal{C}^\mathcal{R}}{s_\mathcal{H} }> \frac{\alpha+\pi-\beta}{\alpha+\pi}$, $\dot{A}_\mathcal{H}(t)>s_\mathcal{H}$ and $\dot{A}_\mathcal{C}^\mathcal{R}(t)>s_\mathcal{C}^\mathcal{R}$ hold, RV queue forms at both the highway and the curbside bottleneck. The following Proposition \ref{Proposition 1} summarizes the properties.

\begin{Proposition}[Initial queuing scenarios]\label{Proposition 1} Given the arrival rate of commuters at bottlenecks in Eq. \eqref{eq:firstphase},
(i) if $\frac{s_\mathcal{C}^\mathcal{R}}{s_\mathcal{H} }\leq \frac{\alpha+\pi-\beta}{\alpha+\pi}$, queues form solely at the curbside bottleneck at the beginning of the morning peak; Conversely, (ii) if $\frac{s_\mathcal{C}^\mathcal{R}}{s_\mathcal{H} }> \frac{\alpha+\pi-\beta}{\alpha+\pi}$, both the curbside bottleneck and the highway bottleneck are congested in the beginning phase of the morning peak.
\end{Proposition}

At user equilibrium, the first RV commuter, as the first arrival during the morning peak, faces no queue, thus we have $w_\mathcal{H}(t_0^\mathcal{R})=w_\mathcal{C}^\mathcal{R}({ t_0^\mathcal{R}+w_\mathcal{H}(t_0^\mathcal{R})})=0$. By contrast, the first PV commuter will experience a queue time of $w_\mathcal{H}(t_0^\mathcal{P})=\frac{A_\mathcal{H}(t_0^\mathcal{P})-s_\mathcal{H}(t_0^\mathcal{P}-t_0^\mathcal{R})}{s_\mathcal{H}}$ at the highway bottleneck if $\frac{s_\mathcal{C}^\mathcal{R}}{s_\mathcal{H}} > \frac{\alpha+\pi-\beta}{\alpha+\pi}$ and $t_0^\mathcal{P}\leq t_1^\mathcal{R}+w_\mathcal{H}(t_1^\mathcal{R})$; otherwise, the first PV commuter faces no queue at the highway bottleneck, i.e., $w_\mathcal{H}(t_0^\mathcal{P})=0$. This is because, according to Proposition \ref{Proposition 1}, the highway bottleneck is congested in the initial stage of the morning peak if $\frac{s_\mathcal{C}^\mathcal{R}}{s_\mathcal{H}} > \frac{\alpha+\pi-\beta}{\alpha+\pi}$ holds. In addition, as the earliest one through the curbside bottleneck serving PVs, the first PV commuter faces no queue at the curbside bottleneck, thus we have $w_\mathcal{C}^\mathcal{P}({ t_0^\mathcal{P}+w_\mathcal{H}(t_0^\mathcal{P})})=0$. Combining the queue times of the first RV and PV commuters at bottlenecks with Eq. \eqref{eq:costDiff}, the inter-equilibrium conditions can be further cast into Eq. \eqref{eq:jointUE}.
\begin{equation}\label{eq:jointUE}
\textbf{Inter-equilibrium conditions:}
    \begin{cases}
        t_0^\mathcal{P}-t_0^\mathcal{R}=\frac{u^\mathcal{P}-c^\mathcal{R}}{\beta}, &\text{if}\  \frac{s_\mathcal{C}^\mathcal{R}}{s_\mathcal{H} }\leq \frac{\alpha+\pi-\beta}{\alpha+\pi}, \\[1.5ex]
        
        t_0^\mathcal{P}-t_0^\mathcal{R}=
        \frac{u^\mathcal{P}-c^\mathcal{R}+(\alpha-\beta)w_\mathcal{H}(t_0^\mathcal{P})}{\beta}, &\text{if}\ \frac{s_\mathcal{C}^\mathcal{R}}{s_\mathcal{H} }> \frac{\alpha+\pi-\beta}{\alpha+\pi}.\\
    \end{cases}
\end{equation}

\section{Different scenarios}\label{section_different_scenarios}
In the preceding section, we examined the user equilibrium conditions and different queuing scenarios at the beginning of the morning peak. In this section, we will further explore possible departure curves and queuing scenarios that can arise throughout the morning peak hours.

Based on whether both travel modes are utilized, whether their departure intervals overlap, and the queuing situations at bottlenecks (as determined by Proposition \ref{Proposition 1}), we can classify all possible user equilibrium situations into 8 scenarios, as presented in Table \ref{table: 8 scenarios}. Except Scenarios 1 and 6, both RV and PV are utilized in the other scenarios. Specifically, the departure intervals of the two modes overlap in Scenarios 3, 5, and 8, thus there is a period when both modes depart from home simultaneously. All scenarios can be solved in closed-form, and the detailed derivation process is available in Appendix \ref{Appendix A}. 

Based on the equilibrium results of scenarios, Table \ref{table: conditions_for_8_scenarios} presents the conditions under which each scenario occurs. On top of this, Proposition \ref{Proposition 2} describes the critical conditions for the concurrent usage of both traffic modes during the morning peak, and the relationship between the departure intervals of RV and PV commuters from home when both modes are utilized. The proof of Proposition \ref{Proposition 2} is given in Appendix \ref{Appendix B}.

\begin{Proposition}[Utilization conditions of RVs and PVs]\label{Proposition 2}
Given the user equilibrium conditions in Eqs. \eqref{eq:innerUE} and \eqref{eq:jointUE} for the bi-modal two-tandem bottleneck model as Fig. \ref{overview of the road network}, (i) if the following condition holds: 
\begin{equation}\label{eq:both_modes_used}
    (u^\mathcal{P}-c^\mathcal{R})<\frac{\beta N}{s_\mathcal{C}^\mathcal{R}}.
\end{equation}
then both RV and PV will be utilized during the morning peak, otherwise, only RV will be utilized; (ii) when both modes are in use, there is no overlap in their departure intervals if the following conditions hold:
\begin{equation}\label{eq:departure_interval_separation}
    \begin{cases}
        (u^\mathcal{P}-c^\mathcal{R})\geq\frac{N\beta(\alpha+\pi-\beta)}{\beta s_\mathcal{C}^\mathcal{P}+(\alpha+\pi)s_\mathcal{C}^\mathcal{R}}, & \text{if}\  \frac{s_\mathcal{C}^\mathcal{R}}{s_\mathcal{H} }\leq \frac{\alpha+\pi-\beta}{\alpha+\pi}, \\[1.5ex]

        (u^\mathcal{P}-c^\mathcal{R})\geq\frac{N\big( \alpha(\alpha+\pi-\beta)s_\mathcal{H}-(\alpha+\pi)(\alpha-\beta)s_\mathcal{C}^\mathcal{R}\big)}{(\alpha+\pi)\big(s_\mathcal{H}(s_\mathcal{C}^\mathcal{R}+s_\mathcal{C}^\mathcal{P})-s_\mathcal{C}^\mathcal{R}s_\mathcal{C}^\mathcal{P}\big)}, & \text{if}\ \frac{s_\mathcal{C}^\mathcal{R}}{s_\mathcal{H} }> \frac{\alpha+\pi-\beta}{\alpha+\pi}.
    \end{cases}
\end{equation}
\end{Proposition}

\begin{table}[h]
    \caption{Overview of 8 scenarios}\label{table: 8 scenarios}
    \begin{center}
        \begin{tabular}{l l l l}\hline
        \diagbox{\textbf{Queuing situation}}{\textbf{Utilization}}
                            & \multirow{2}{*}{\textbf{only RV}}  &  \multicolumn{2}{c}{\textbf{RV+PV}}\\
                            &          &	\textbf{without overlap}	&	\textbf{with overlap} \\\hline
            Only $B_\mathcal{C}$ queuing	            & Scenario 1      &	Scenario 2	& Scenario 3\\
            $B_\mathcal{C}$ queuing, expanding to  $B_\mathcal{C}+B_\mathcal{H}$ queuing & $-$             &Scenario 4   &	Scenario 5	\\
            $B_\mathcal{C}+B_\mathcal{H}$	queuing	        & Scenario 6      &	Scenario 7	&	Scenario 8\\\hline
            \multicolumn{4}{l}{\textbf{Note:} Only $B_\mathcal{H}$ queuing will not happen as $s_\mathcal{H}>s_\mathcal{C}^\mathcal{P}$, $s_\mathcal{H}>s_\mathcal{C}^\mathcal{R}$ hold in this paper}
      \end{tabular}
    \end{center}
\end{table}

\begin{table}[h]
    \caption{Conditions for 8 scenarios}\label{table: conditions_for_8_scenarios}
        \begin{center}
        \resizebox{\linewidth}{!}{
        \begin{tabular}{l l l c}\hline
            \textbf{Occurrence Conditions}&   &     & \textbf{Scenario} \\\hline
            \multirow{7}{*}{$\frac{s_\mathcal{C}^\mathcal{R}}{s_\mathcal{H} }\leq \frac{\alpha+\pi-\beta}{\alpha+\pi}$}  
            & $(u^\mathcal{P}-c^\mathcal{R})\geq\frac{\beta N}{s_\mathcal{C}^\mathcal{R}}$ 
            & $\frac{s_\mathcal{C}^\mathcal{P}}{s_\mathcal{H}}\in (0,1)$ 
            &	 1\\[1.5ex]
            \cline{2-4}
            & \multirow{3}{*}{$(u^\mathcal{P}-c^\mathcal{R})\in \big[\frac{N\beta(\alpha+\pi-\beta)}{\beta s_\mathcal{C}^\mathcal{P}+(\alpha+\pi)s_\mathcal{C}^\mathcal{R}},\frac{\beta N}{s_\mathcal{C}^\mathcal{R}}\big)$} 
            & $\frac{s_\mathcal{C}^\mathcal{P}}{s_\mathcal{H}}\in \big(0,\frac{\alpha-\beta}{\alpha}\big]$
            &	 2\\[1.5ex]
            & 
            & $\frac{s_\mathcal{C}^\mathcal{P}}{s_\mathcal{H}}\in \big(\frac{\alpha-\beta}{\alpha},1\big)$ 
            &	 4\\[1.5ex]
            \cline{2-4}
            & \multirow{3}{*}{$(u^\mathcal{P}-c^\mathcal{R})\in \big(0,\frac{N\beta(\alpha+\pi-\beta)}{\beta s_\mathcal{C}^\mathcal{P}+(\alpha+\pi)s_\mathcal{C}^\mathcal{R}}\big)$}
            & $\frac{s_\mathcal{C}^\mathcal{P}}{s_\mathcal{H}}\in \big(0,\frac{\alpha-\beta}{\alpha}-\frac{(1-\delta^\mathcal{R})(\alpha+\pi)(\alpha-\beta)s_\mathcal{C}^\mathcal{R}}{\alpha(\alpha+\pi-\beta)s_\mathcal{H}}\big]$ 
            &	 3\\[1.5ex]
            &  
            & $\frac{s_\mathcal{C}^\mathcal{P}}{s_\mathcal{H}}\in \big(\frac{\alpha-\beta}{\alpha}-\frac{(1-\delta^\mathcal{R})(\alpha+\pi)(\alpha-\beta)s_\mathcal{C}^\mathcal{R}}{\alpha(\alpha+\pi-\beta)s_\mathcal{H}},1\big)$ 
            &	 5\\[1.5ex]
            \hline
            \multirow{5}*{$\frac{s_\mathcal{C}^\mathcal{R}}{s_\mathcal{H} }> \frac{\alpha+\pi-\beta}{\alpha+\pi}$}
            & \multicolumn{2}{l}{$(u^\mathcal{P}-c^\mathcal{R})\geq\frac{\beta N}{s_\mathcal{C}^\mathcal{R}}$} 
            &	 6\\
            &\multicolumn{2}{l}{$(u^\mathcal{P}-c^\mathcal{R})\in\bigg[\frac{N\big( \alpha(\alpha+\pi-\beta)s_\mathcal{H}-(\alpha+\pi)(\alpha-\beta)s_\mathcal{C}^\mathcal{R}\big)}{(\alpha+\pi)\big(s_\mathcal{H}(s_\mathcal{C}^\mathcal{R}+s_\mathcal{C}^\mathcal{P})-s_\mathcal{C}^\mathcal{R}s_\mathcal{C}^\mathcal{P}\big)},\frac{\beta N}{s_\mathcal{C}^\mathcal{R}}\bigg)$} 
            &	 7\\
            &\multicolumn{2}{l}{$(u^\mathcal{P}-c^\mathcal{R})\in\bigg(0,\frac{N\big( \alpha(\alpha+\pi-\beta)s_\mathcal{H}-(\alpha+\pi)(\alpha-\beta)s_\mathcal{C}^\mathcal{R}\big)}{(\alpha+\pi)\big(s_\mathcal{H}(s_\mathcal{C}^\mathcal{R}+s_\mathcal{C}^\mathcal{P})-s_\mathcal{C}^\mathcal{R}s_\mathcal{C}^\mathcal{P}\big)}\bigg)$} 
            &    8\\
            \hline
      \end{tabular}}
      \end{center}
  \end{table}

\begin{figure}[!ht]
    \centering
    \includegraphics[width=\linewidth]{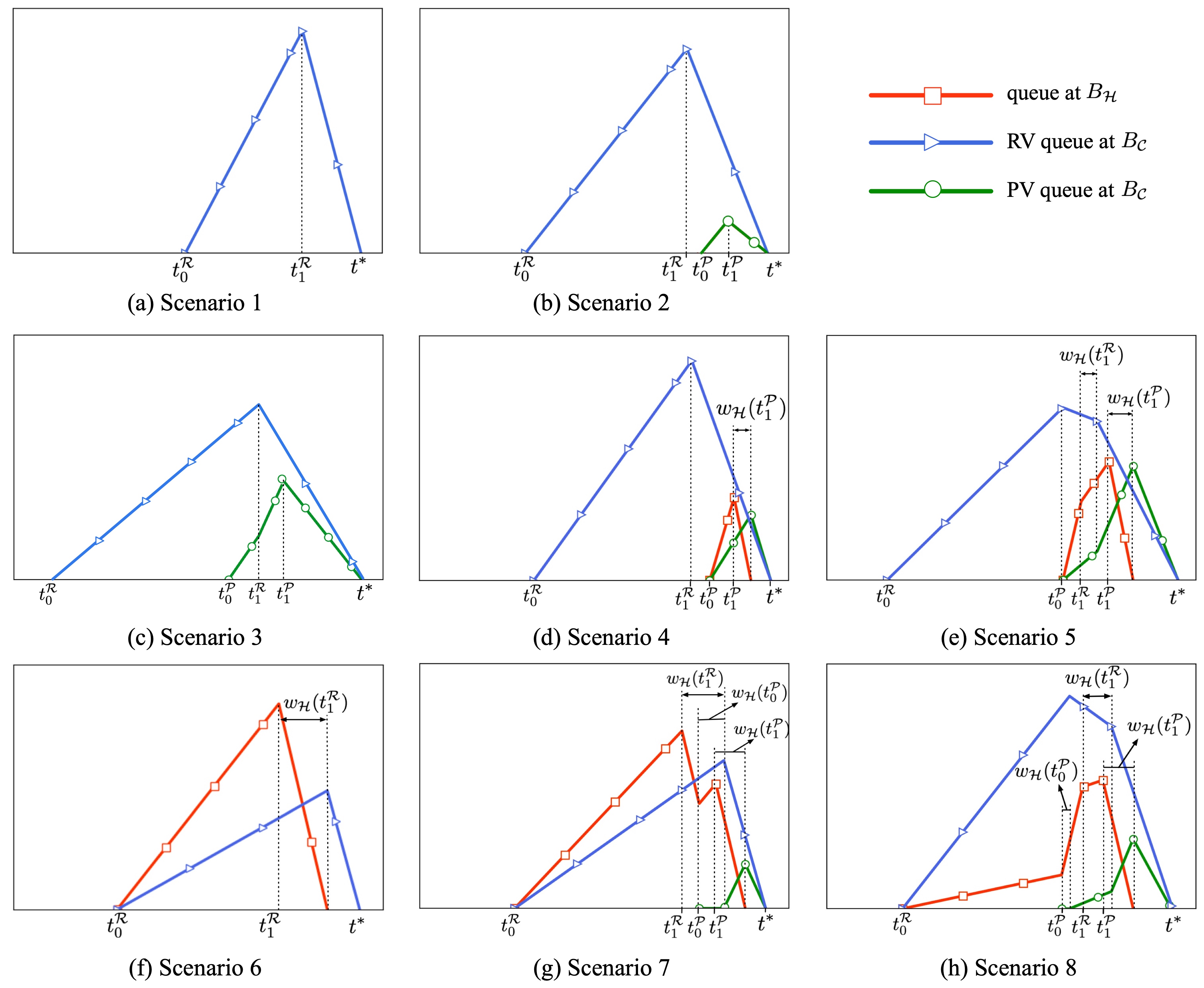}
    \caption{Queue length for different scenarios}
    \label{fig:eight_images}
\end{figure}

For ease of understanding, Fig. \ref{fig:eight_images} plots the changes in queue length at bottlenecks over time for these 8 scenarios. Here, we briefly discuss each scenario.

For Scenarios 1-5, shown in Fig. \ref{fig:eight_images}(a)-(e), only the RV mode is utilized and a queue forms exclusively at the curbside bottleneck $B_\mathcal{C}$ at the beginning of the morning peak. This is attributed to the lower fixed cost of RV compared to that of PV (i.e., $c^\mathcal{R} < u^\mathcal{P}$), and the significantly higher service rate of the highway bottleneck compared to the curbside bottleneck serving RVs (i.e., $s_\mathcal{H} \geq \frac{\alpha+\pi}{\alpha+\pi-\beta}s_\mathcal{C}^\mathcal{R}$). However, as the RV queue builds up at the curbside bottleneck, commuters will be willing to switch to the PV mode (except Scenario 1 where only RV is used due to the huge difference in the fixed costs of two modes). Specifically,
\begin{itemize}
    \item as $(u^\mathcal{P}-c^\mathcal{R})\geq\frac{\beta N}{s_\mathcal{C}^\mathcal{R}}$, Scenario 1 occurs and only the RV mode will be utilized during the morning peak. The queuing scenario shown in Fig. \ref{fig:eight_images}(a) is the same as the traditional bottleneck model \citep{vickrey1969}.
    
    \item as $(u^\mathcal{P}-c^\mathcal{R})<\frac{\beta N}{s_\mathcal{C}^\mathcal{R}}$, commuters will be willing to use the PV mode when the queuing delay for RV increases, as shown in Fig. \ref{fig:eight_images}(b)-(e). If $(u^\mathcal{P}-c^\mathcal{R})\in \big[\frac{N\beta(\alpha+\pi-\beta)}{\beta s_\mathcal{C}^\mathcal{P}+(\alpha+\pi)s_\mathcal{C}^\mathcal{R}},\frac{\beta N}{s_\mathcal{C}^\mathcal{R}}\big)$, then RV and PV will be used with non-overlapping departure intervals as shown in Fig. \ref{fig:eight_images}(b) and Fig. \ref{fig:eight_images}(d). More specifically, when the service rate of the curbside bottleneck serving PVs, $s_\mathcal{C}^\mathcal{P}$, is lower than $\frac{\alpha-\beta}{\alpha}$, Scenario 2 occurs; otherwise, Scenario 4 occurs. Conversely, if $(u^\mathcal{P}-c^\mathcal{R})\in \big(0,\frac{N\beta(\alpha+\pi-\beta)}{\beta s_\mathcal{C}^\mathcal{P}+(\alpha+\pi)s_\mathcal{C}^\mathcal{R}}\big)$, then RV and PV will be used with a co-departure interval of $(t_0^\mathcal{P},t_1^\mathcal{R})$, as shown in Fig. \ref{fig:eight_images}(c) and \ref{fig:eight_images}(e). More specifically, Scenario 3 occurs when $s_\mathcal{C}^\mathcal{P}$ is relatively low; otherwise, Scenario 5 occurs.
\end{itemize}

For Scenarios 6-8, queues will form at both the curbside bottleneck and the highway bottleneck at the beginning of the morning peak. This is because the capacity of the curbside bottleneck serving RVs is closer to the highway bottleneck in these scenarios. Specifically,
\begin{itemize}
    \item as $(u^\mathcal{P}-c^\mathcal{R})\geq\frac{\beta N}{s_\mathcal{C}^\mathcal{R}}$,
    Scenario 6 occurs and only the RV mode will be used during the morning peak as shown in Fig. \ref{fig:eight_images}(f).
    
    \item as $(u^\mathcal{P}-c^\mathcal{R})<\frac{\beta N}{s_\mathcal{C}^\mathcal{R}}$,
    commuters will also use the PV mode when the queuing delay for RV increases. More specifically, if $(u^\mathcal{P}-c^\mathcal{R})\geq\frac{N\big( \alpha(\alpha+\pi-\beta)s_\mathcal{H}-(\alpha+\pi)(\alpha-\beta)s_\mathcal{C}^\mathcal{R}\big)}{(\alpha+\pi)\big(s_\mathcal{H}(s_\mathcal{C}^\mathcal{R}+s_\mathcal{C}^\mathcal{P})-s_\mathcal{C}^\mathcal{R}s_\mathcal{C}^\mathcal{P}\big)}$, then
    Scenario 7 occurs and the departure intervals of RV and PV commuters do not overlap (i.e., $t_1^\mathcal{R}\leq t_0^\mathcal{P}$), as shown in Fig. \ref{fig:eight_images}(g).
    Conversely, if $(u^\mathcal{P}-c^\mathcal{R})<\frac{N\big( \alpha(\alpha+\pi-\beta)s_\mathcal{H}-(\alpha+\pi)(\alpha-\beta)s_\mathcal{C}^\mathcal{R}\big)}{(\alpha+\pi)\big(s_\mathcal{H}(s_\mathcal{C}^\mathcal{R}+s_\mathcal{C}^\mathcal{P})-s_\mathcal{C}^\mathcal{R}s_\mathcal{C}^\mathcal{P}\big)}$, then
    Scenario 8 occurs and the departure intervals of RV and PV commuters overlap (i.e., $t_1^\mathcal{R}>t_0^\mathcal{P}$), as shown in Fig. \ref{fig:eight_images}(h).
\end{itemize}

\section{Bidirectional congestion spillover effects}\label{section_bieffects}
In the previous section, we considered the unidirectional congestion spillover from the RV queue on PVs. Inevitably, congestion caused by PVs on the main road may also limit RVs’ ability to reach the curbside, implying that the congestion spillover effects may be mutual\footnote{We would like to thank the two anonymous reviewers for this insightful suggestion.}. Therefore, we will now discuss the effects of bidirectional congestion spillover on the morning peak commuting equilibrium.

\subsection{Spillover effects of PV queue}
Similarly, we assumed that the actual service rate of the curbside bottleneck serving RVs is discounted due to the congestion caused by the PV queue on the main road, and the actual service rate of the curbside bottleneck serving RVs, $\widetilde{s}_\mathcal{C}^\mathcal{R}(t)$, can be defined as:
\begin{equation}\label{eq:scr_actual}
        \widetilde{s}_\mathcal{C}^\mathcal{R}(t)=\begin{cases}
        \frac{\dot{A}_\mathcal{C}^\mathcal{R}(t)}{\dot{A}_\mathcal{C}^\mathcal{R}(t)+\delta^\mathcal{P}\dot{A}_\mathcal{C}^{\mathcal{P}}(t)}s_\mathcal{C}^\mathcal{R},\ &\text{if}\ \dot{A}_\mathcal{C}^\mathcal{P}(t)>0, \dot{A}_\mathcal{C}^\mathcal{R}(t)>0, \\
        s_\mathcal{C}^\mathcal{R},\ &\text{else.} \end{cases}
\end{equation} where $\delta^\mathcal{P}$ represent the intensity of the spillover congestion effects from the PV queue on RVs. Specifically, $\delta^\mathcal{P} = 0$ means that the PV queue has no impact on RVs and $\widetilde{s}_\mathcal{C}^\mathcal{P}(t)={s}_\mathcal{C}^\mathcal{P}$ always holds, which corresponds to the case with unidirectional congestion spillover discussed in the previous sections.

Therefore, the queuing delay for RVs arriving at the curbside bottleneck at time $t$ is presented by Eq. \eqref{eq:wcr_bidirectional}.
\begin{equation}\label{eq:wcr_bidirectional}
    w_\mathcal{C}^\mathcal{R}(t)=\frac{A_\mathcal{C}^\mathcal{R}(t)-D_\mathcal{C}^\mathcal{R}(t)}{\widetilde{s}_\mathcal{C}^\mathcal{R}(t)},
\end{equation}
where the slopes of cumulative curves, $\dot{D}_\mathcal{C}^\mathcal{R}(t)$ and $\dot{A}_\mathcal{C}^\mathcal{R}(t)$, meet: 
\begin{equation}\label{eq:rate1r_bidirectional}
    \dot{D}_\mathcal{C}^\mathcal{R}(t)=
    \begin{cases}
        \dot{A}_\mathcal{C}^\mathcal{R}(t),\ &\text{if}\  w_\mathcal{C}^\mathcal{R}(t)=0 \ \text{and}\ \dot{A}_\mathcal{C}^\mathcal{R}(t)\leq \widetilde{s}_\mathcal{C}^\mathcal{R}(t),\\
        \widetilde{s}_\mathcal{C}^\mathcal{R}(t),\ &\text{else}.
    \end{cases}
\end{equation} Together with Eqs. \eqref{eq:w2} and \eqref{eq:wcp}, Eq. \eqref{eq:wcr_bidirectional} defines the travel times (i.e., queuing delays at bottlenecks) for RV and PV commuters under bidirectional congestion spillover effects. 

\subsection{User equilibrium}
Using Eq. \eqref{eq:costRV}, Eq. \eqref{eq:costPV}, and the condition of $\frac{\partial C^\mathcal{R}(t)}{\partial t}=\frac{\partial C^\mathcal{P}(t)}{\partial t}=0$, the intra-equilibrium conditions for RV and PV commuters under bidirectional spillover effects can be derived. The intra-equilibrium conditions are expressed as differential equations regarding the queuing delays at bottlenecks $w_\mathcal{H}(t)$, $w_\mathcal{C}^\mathcal{R}(t+w_\mathcal{H}(t))$, and $w_\mathcal{C}^\mathcal{P}(t+w_\mathcal{H}(t))$, which have the same form as those under unidirectional congestion spillover, as shown in Eq. \eqref{eq:innerUE}.
        

As $u^\mathcal{P} > c^\mathcal{R}$, only RV commuters depart from home in the initial stage of the morning peak. By combining Eqs. \eqref{eq:w2}, \eqref{eq:wcr_bidirectional}, and the intra-equilibrium conditions of Eq. \eqref{eq:innerUE}, one can derive and find that the equilibrium departure rate for the initial stage remains the same as that in the case with unidirectional congestion spillover, as given in Eq. \eqref{eq:firstphase}. Therefore, Proposition \ref{Proposition 1}, which identifies the queuing situations in the initial stage of morning peak, also holds when considering the bidirectional congestion spillover effects.

The inter-equilibrium requires the equalization of the generalized travel costs for RVs and PVs if both modes are utilized, i.e., $C^\mathcal{R}(t_0^\mathcal{R})=C^\mathcal{P}(t_0^\mathcal{P})=0$ as defined in Eq. \eqref{eq:costDiff}. At time $t_0^\mathcal{R}$ and $t_0^\mathcal{P}$, there is no PV queue at the curbside bottleneck, and thus the bidirectional congestion spillover have not yet occurred. Therefore, the inter-equilibrium conditions remain the same as in the case with unidirectional congestion spillover, as defined in Eq. \eqref{eq:jointUE}.
        

The above analysis indicates that, when considering bidirectional congestion spillover effects, the queuing situation and equilibrium departure rate in the initial stage of the morning peak, as well as the inter-equilibrium conditions for RV and PV commuters, are the same as in the case with unidirectional congestion spillover. In addition, the bidirectional congestion spillover does not occur when only one travel mode is used, or when both modes are used with non-overlapping departure intervals. Therefore, Proposition \ref{Proposition 2} remains applicable for the case with bidirectional congestion spillover, which identifies the critical conditions for scenarios where only one travel mode is used and where both modes are used with non-overlapping departure intervals.

\subsection{Typical scenarios}\label{subsec_bieffect_Typical_scenarios}
For the case with bidirectional congestion spillover, Proposition \ref{Proposition 1} and Proposition \ref{Proposition 2} still hold, thus the morning peak commuting equilibrium can also be categorized into 8 scenarios, similar to the case with unidirectional congestion spillover. The classification rule is the same, as shown in Table \ref{table: 8 scenarios}. Here, we select three scenarios where both RV and PV are utilized with overlapping departure intervals\footnote{Here we focus on the cases with $t_1^\mathcal{P} > t_1^\mathcal{R}$, meaning that only PV commuters depart after $t_1^\mathcal{R}$ during the morning peak. The analysis for the cases with $t_1^\mathcal{P} < t_1^\mathcal{R}$, where only RV commuters depart after $t_1^\mathcal{P}$, follows a similar process and will not be elaborated here.}, i.e., Scenarios 3, 5, and 8, to discuss the impact of bidirectional congestion spillover effects on the equilibrium patterns.

In Scenario 3, RV and PV are utilized with overlapping departure intervals, and queuing congestion only occurs at the curbside bottlenecks during the morning peak. This implies that $w_\mathcal{H}(t)=0$ holds for any departure time $t$. The entire morning peak period can be divided into three stages: in the initial stage, only RVs depart; in the second stage, both RVs and PVs depart, and thus the bidirectional congestion spillover effects occur; in the final stage, only PVs depart. In the initial and final stages, the conditions of $\widetilde{s}_\mathcal{C}^\mathcal{R}(t)=s_\mathcal{C}^\mathcal{R}$ and $\widetilde{s}_\mathcal{C}^\mathcal{P}(t)=s_\mathcal{C}^\mathcal{P}$ hold. The equilibrium departure rates for RVs and PVs can be derived based on their respective intra-equilibrium condition of Eq. \eqref{eq:innerUE}, and thus we have\footnote{Because there is no queuing delay at the highway bottleneck, the equilibrium departure rates for RV and PV commuters from home are also their arrival rates at the curbside bottlenecks.},
\begin{equation}\label{eq:bieffect_DepartRate_S3}
    \begin{cases}
        \dot{A}_\mathcal{C}^\mathcal{R}(t) = \frac{\alpha+\pi}{\alpha+\pi-\beta}s_\mathcal{C}^\mathcal{R},\ &t\in(t_0^\mathcal{R},t_0^\mathcal{P}),\\
        \dot{A}_\mathcal{C}^\mathcal{P}(t) = \frac{\alpha}{\alpha-\beta}s_\mathcal{C}^\mathcal{P},\ &t\in(t_1^\mathcal{R},t_1^\mathcal{P}).
    \end{cases}
\end{equation}

In the co-departure stage of $t\in(t_0^\mathcal{P},t_1^\mathcal{R})$, both RV and PV commuters depart from home. Due to the bidirectional congestion spillover effects, the actual service rates of the curbside bottlenecks serving PVs and RVs are discounted, and can be determined by Eqs. \eqref{eq:scp_actual} and \eqref{eq:scr_actual}. By substituting Eqs. \eqref{eq:scp_actual}-\eqref{eq:u1c} and Eqs. \eqref{eq:scr_actual}-\eqref{eq:rate1r_bidirectional} into the user equilibrium conditions given by Eq. \eqref{eq:innerUE}, the equilibrium departure rates for RVs and PVs in the co-departure stage can be derived as follows:
\begin{equation}\label{eq:bieffect_CoDepartRate_S3}
    \begin{cases}
        \dot{A}_\mathcal{C}^\mathcal{R}(t) = \frac{(\alpha+\pi)s_\mathcal{C}^\mathcal{R}}{(\alpha+\pi-\beta)(1-\delta^\mathcal{R}\delta^\mathcal{P})} - \frac{\alpha\delta^\mathcal{P}s_\mathcal{C}^\mathcal{P}}{(\alpha-\beta)(1-\delta^\mathcal{R}\delta^\mathcal{P})},\ &t\in(t_0^\mathcal{P},t_1^\mathcal{R}),\\
        
        \dot{A}_\mathcal{C}^\mathcal{P}(t) = \frac{\alpha s_\mathcal{C}^\mathcal{P}}{(\alpha-\beta)(1-\delta^\mathcal{R}\delta^\mathcal{P})}-\frac{(\alpha+\pi)\delta^\mathcal{R}s_\mathcal{C}^\mathcal{R}}{(\alpha+\pi-\beta)(1-\delta^\mathcal{R}\delta^\mathcal{P})},\ &t\in(t_0^\mathcal{P},t_1^\mathcal{R}).
    \end{cases}
\end{equation} where $\delta^\mathcal{R}$ and $\delta^\mathcal{P}$ represent the impact intensity of the RV queue on PVs and the PV queue on RVs, respectively. Specifically, $\delta^\mathcal{P} = 0$ means that the PV queue has no effect on RVs. In this case, the equilibrium departure rates during the co-departure stage are consistent with that in the case with unidirectional congestion spillover (see Eq. \eqref{eq:UEresults_S3_2} in Appendix \ref{Appendix A}). Substituting Eq. \eqref{eq:bieffect_CoDepartRate_S3} into Eqs. \eqref{eq:scp_actual} and \eqref{eq:scr_actual}, one can further obtain the actual service rates of the curbside bottlenecks.

It can be proven that, compared to the case with unidirectional congestion spillover (i.e., Eq. \eqref{eq:UEresults_S3_2}), the bidirectional congestion spillover results in a higher departure rate of PV commuters and a lower departure rate of RV commuters from home (i.e., arrival rates at the curbside bottlenecks, $\dot{A}_\mathcal{C}^\mathcal{P}(t)$ and $\dot{A}_\mathcal{C}^\mathcal{R}(t)$, respectively) during the co-departure stage. It also results in a lower total departure rate of RV and PV commuters during the co-departure stage (i.e., $\dot{A}_\mathcal{C}^\mathcal{R}(t)+\dot{A}_\mathcal{C}^\mathcal{P}(t)$). As there is no queue at the highway bottleneck, the total departure rate of RV and PV commuters should be lower than the highway bottleneck capacity, i.e., $\dot{A}_\mathcal{C}^\mathcal{R}(t)+\dot{A}_\mathcal{C}^\mathcal{P}(t)<s_\mathcal{H}$. Therefore, one can obtain the condition under which no queue forms at the highway bottleneck during the co-departure stage:
\begin{equation}\label{eq:bieffect_s3&s5_condition}
    s_\mathcal{H}>\frac{(\alpha+\pi)(1-\delta^\mathcal{R})s_\mathcal{C}^\mathcal{R}}{(\alpha+\pi-\beta)(1-\delta^\mathcal{R}\delta^\mathcal{P})}+\frac{\alpha(1-\delta^\mathcal{P})s_\mathcal{C}^\mathcal{P}}{(\alpha-\beta)(1-\delta^\mathcal{R}\delta^\mathcal{P})}.
\end{equation} 

Eq. \eqref{eq:bieffect_s3&s5_condition} also distinguishes Scenario 3 from Scenario 5 when considering bidirectional congestion spillover effects. Setting  $\delta^\mathcal{P} = 0$ yields the boundary condition between Scenarios 3 and 5 under unidirectional congestion spillover, as shown in Table \ref{table: conditions_for_8_scenarios}.

In Scenario 5, both RV and PV are utilized with overlapping departure intervals, and queuing congestion forms at the highway bottleneck after the first PV commuter departs from home. By contrast, in Scenario 8, queuing congestion occurs at the highway bottleneck from the beginning of the morning peak. Following a similar analytical process to that of Scenario 3, we can derive the equilibrium departure rates for Scenarios 5 and 8 as follows:
\begin{itemize}
    \item \text{Scenario 5:}
    \begin{equation}\label{eq:bieffect_DepartRate_S5}
    \begin{aligned}
        &\dot{A}_\mathcal{H}(t) = 
            \begin{cases}
                \frac{(\alpha+\pi)s_\mathcal{C}^\mathcal{R}}{\alpha+\pi-\beta}\leq s_\mathcal{H},\ &t\in(t_0^\mathcal{R},t_0^\mathcal{P}),\\[1.5ex]

                \frac{\alpha(1-\delta^\mathcal{P})s_\mathcal{C}^\mathcal{P}}{(\alpha-\beta)(1-\delta^\mathcal{R}\delta^\mathcal{P})}+\frac{(\alpha+\pi)(1-\delta^\mathcal{R})s_\mathcal{C}^\mathcal{R}}{(\alpha+\pi-\beta)(1-\delta^\mathcal{R}\delta^\mathcal{P})},\ &t\in(t_0^\mathcal{P},t_1^\mathcal{R}),\\[1.5ex]
                
                \frac{\alpha s_\mathcal{C}^\mathcal{P}}{\alpha-\beta},\ &t\in(t_1^\mathcal{R},t_1^\mathcal{P}).
            \end{cases}\\
        &\dot{A}_\mathcal{C}^\mathcal{R}(t) =
            \begin{cases}
                \frac{(\alpha+\pi)s_\mathcal{C}^\mathcal{R}}{\alpha+\pi-\beta},\ &t\in(t_0^\mathcal{R},t_0^\mathcal{P}),\\[1.5ex]
                
                \frac{(\alpha+\pi)(\alpha-\beta)s_\mathcal{C}^\mathcal{R}-\alpha(\alpha+\pi-\beta)\delta^\mathcal{P}s_\mathcal{C}^\mathcal{P}}{(\alpha+\pi)(\alpha-\beta)(1-\delta^\mathcal{R})s_\mathcal{C}^\mathcal{R}+\alpha(\alpha+\pi-\beta)(1-\delta^\mathcal{P})s_\mathcal{C}^\mathcal{P}}s_\mathcal{H},\ &t\in\big(t_0^\mathcal{P},t_1^\mathcal{R}+w_\mathcal{H}(t_1^\mathcal{R})\big).
            \end{cases}\\
        &\dot{A}_\mathcal{C}^\mathcal{P}(t) =
            \begin{cases}
                s_\mathcal{H} - \dot{A}_\mathcal{C}^\mathcal{R}(t),\ &t\in\big(t_0^\mathcal{P},t_1^\mathcal{R}+w_\mathcal{H}(t_1^\mathcal{R})\big),\\
                
                s_\mathcal{H},\ &t\in\big(t_1^\mathcal{R}+w_\mathcal{H}(t_1^\mathcal{R}),t_1^\mathcal{P}+w_\mathcal{H}(t_1^\mathcal{P})\big).
            \end{cases}
    \end{aligned}
    \end{equation}
    \item \text{Scenario 8:}
    \begin{equation}\label{eq:bieffect_DepartRate_S8}
        \begin{aligned}
        &\dot{A}_\mathcal{H}(t) = 
            \begin{cases}
                \frac{(\alpha+\pi)s_\mathcal{C}^\mathcal{R}}{\alpha+\pi-\beta}>s_\mathcal{H},\ &t\in(t_0^\mathcal{R},t_0^\mathcal{P}),\\[1.5ex]
                
                \frac{\alpha(1-\delta^\mathcal{P})s_\mathcal{C}^\mathcal{P}}{(\alpha-\beta)(1-\delta^\mathcal{R}\delta^\mathcal{P})}+\frac{(\alpha+\pi)(1-\delta^\mathcal{R})s_\mathcal{C}^\mathcal{R}}{(\alpha+\pi-\beta)(1-\delta^\mathcal{R}\delta^\mathcal{P})},\ &t\in(t_0^\mathcal{P},t_1^\mathcal{R}),\\[1.5ex]
                
                \frac{\alpha s_\mathcal{C}^\mathcal{P}}{\alpha-\beta},\ &t\in(t_1^\mathcal{R},t_1^\mathcal{P}).
            \end{cases}\\
        &\dot{A}_\mathcal{C}^\mathcal{R}(t) =
            \begin{cases}
                s_\mathcal{H},\ &t\in(t_0^\mathcal{R},t_0^\mathcal{P}+w_\mathcal{H}(t_0^\mathcal{P})),\\
                
                \frac{(\alpha+\pi)(\alpha-\beta)s_\mathcal{C}^\mathcal{R}-\alpha(\alpha+\pi-\beta)\delta^\mathcal{P}s_\mathcal{C}^\mathcal{P}}{(\alpha+\pi)(\alpha-\beta)(1-\delta^\mathcal{R})s_\mathcal{C}^\mathcal{R}+\alpha(\alpha+\pi-\beta)(1-\delta^\mathcal{P})s_\mathcal{C}^\mathcal{P}}s_\mathcal{H},\ &t\in\big(t_0^\mathcal{P}+w_\mathcal{H}(t_0^\mathcal{P}),t_1^\mathcal{R}+w_\mathcal{H}(t_1^\mathcal{R})\big).
            \end{cases}\\
        &\dot{A}_\mathcal{C}^\mathcal{P}(t) =
            \begin{cases}
                s_\mathcal{H} - \dot{A}_\mathcal{C}^\mathcal{R}(t),\ &t\in\big(t_0^\mathcal{P}+w_\mathcal{H}(t_0^\mathcal{P}),t_1^\mathcal{R}+w_\mathcal{H}(t_1^\mathcal{R})\big),\\
                
                s_\mathcal{H},\ &t\in\big(t_1^\mathcal{R}+w_\mathcal{H}(t_1^\mathcal{R}),t_1^\mathcal{P}+w_\mathcal{H}(t_1^\mathcal{P})\big).
            \end{cases}
        \end{aligned}
    \end{equation}
\end{itemize}

For Scenario 5, we compare Eqs. \eqref{eq:bieffect_DepartRate_S5} and \eqref{eq:arrivalrate_S5_t1c>t1r}; for Scenario 8, we compare Eqs. \eqref{eq:bieffect_DepartRate_S8} and \eqref{eq:arrivalrate_S8_t1c>t1r}. Based on the comparisons, one can prove that incorporating the effects of bidirectional congestion spillover results in a lower total departure rate from home of PV and RV commuters $\dot{A}_\mathcal{H}(t)$ during the co-departure stage. Additionally, it leads to a lower arrival rate of RVs at the curbside bottleneck $\dot{A}_\mathcal{C}^\mathcal{R}(t)$ and a higher arrival rate of PVs at the curbside bottleneck $\dot{A}_\mathcal{C}^\mathcal{P}(t)$, compared to the case with unidirectional congestion spillover. Combining this with the findings from Scenario 3, the following proposition summarizes these observations.

\begin{Proposition}[Effects of Bidirectional Congestion Spillover]\label{Proposition 3}
    When both travel modes are utilized with overlapping departure intervals, the introduction of bidirectional congestion spillover results in a reduced total departure rate of RV and PV commuters during the co-departure stage. The arrival rate of those RV commuters at the curbside bottleneck decreases, while the arrival rate of those PV commuters at the curbside bottleneck increases.
\end{Proposition}

It is intuitive that the congestion spillover from the PV queue reduces the actual service rate of the curbside bottleneck serving RVs, leading to a decrease in the equilibrium arrival rate of RVs at the curbside bottleneck. Furthermore, the reduced arrival rate of RVs at the curbside bottleneck lessens their spillover congestion effect on PVs during the co-departure stage, allowing the arrival rate of PVs at the curbside bottleneck to increase. However, the bidirectional congestion spillover effects reduce the overall service rate of the curbside bottlenecks compared to the case with unidirectional congestion spillover, thus it could result in a lower total departure rate of RV and PV commuters during the co-departure stage.

\begin{figure}[!h]
    \centering
    \includegraphics[width=\linewidth]{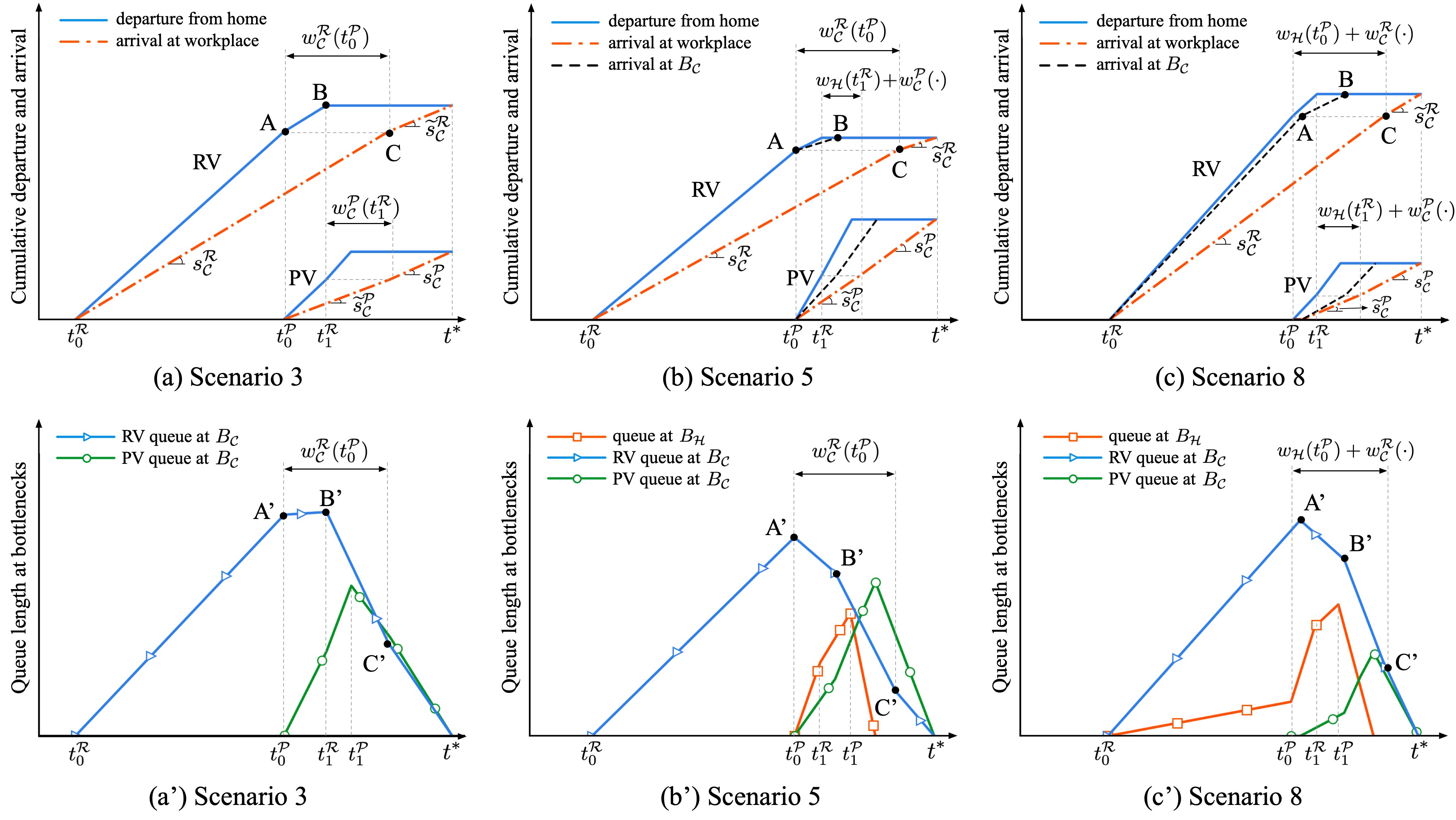}
    \caption{Equilibrium patterns with bidirectional congestion spillover effects: (a)-(c) illustrate the cumulative departure and arrival of commuters; (a’)-(c’) depict the queue length at bottlenecks over time}
    \label{fig:bieffects_S358}
\end{figure}

For ease of understanding, Fig. \ref{fig:bieffects_S358} illustrates equilibrium patterns under bidirectional congestion spillover effects, showing the cumulative departure and arrival of commuters over time and changes in queue length at bottlenecks during the morning peak for Scenarios 3, 5, and 8. The derivation process for the equilibrium results, including the critical time points (i.e., $t_0^\mathcal{R},t_0^\mathcal{P},t_1^\mathcal{R},t_1^\mathcal{P}$) and the mode split result $N^\mathcal{R}$, can be found in Appendix \ref{Appendix C}. For presentation purpose, we define $\widetilde{s}_\mathcal{C}^\mathcal{R} = \frac{\dot{A}_\mathcal{C}^\mathcal{R}(t)}{\dot{A}_\mathcal{C}^\mathcal{R}(t)+\delta^\mathcal{P}\dot{A}_\mathcal{C}^{\mathcal{P}}(t)}s_\mathcal{C}^\mathcal{R}$ and $\widetilde{s}_\mathcal{C}^\mathcal{P} = \frac{\dot{A}_\mathcal{C}^\mathcal{P}(t)}{\dot{A}_\mathcal{C}^\mathcal{P}(t)+\delta^\mathcal{R}\dot{A}_\mathcal{C}^{\mathcal{R}}(t)}s_\mathcal{C}^\mathcal{P}$ to represent the discounted service rates of the curbside bottlenecks, and thus $\widetilde{s}_\mathcal{C}^\mathcal{R}<{s}_\mathcal{C}^\mathcal{R}$ and $\widetilde{s}_\mathcal{C}^\mathcal{P}<{s}_\mathcal{C}^\mathcal{P}$ hold according to Eqs. \eqref{eq:scp_actual} and \eqref{eq:scr_actual}.

Fig. \ref{fig:bieffects_S358}(a) and (a') indicate, respectively, the cumulative departure and arrival of RV and PV commuters, and the changes in queue length at bottlenecks over time in Scenario 3. In Fig. \ref{fig:bieffects_S358}(a), the slopes of the cumulative departure (arrival) line represent the departure (arrival) rates of commuters. It can be observed that the RV departure rate (i.e., slopes of blue solid line) decreases after time $t_0^\mathcal{P}$ (see point A) and drops to zero after the last RV commuter departs at time $t_1^\mathcal{R}$ (see point B). In addition, the arrival rate of RV commuters at the workplace (i.e., slopes of red dash-dotted line) decreases from  $s_\mathcal{C}^\mathcal{R}$  to  $\widetilde{s}_\mathcal{C}^\mathcal{R}$ at time $t_0^\mathcal{P}+w_\mathcal{C}^\mathcal{R}(t_0^\mathcal{P})$ (see point C).
Consequently, as shown at points A’, B’, and C’ in Fig. \ref{fig:bieffects_S358}(a’), the RV queue length at the curbside bottleneck changes at the corresponding time points. Obviously, the changes in queue length at points A’ and C’ are due to the spillover effects from the PV queue. This variation does not occur in the case with unidirectional congestion spillover (see Fig. \ref{fig:eight_images}(c)), because the departure and arrival rates of RV commuters will not change during the morning peak (see Eq. \eqref{eq:UEresults_S3_2}). By contrast, the changes in cumulative departure, cumulative arrival, and queue length at bottlenecks for PV commuters are similar to those in the unidirectional congestion spillover case.

Fig. \ref{fig:bieffects_S358}(b) and (b') indicate, respectively, the cumulative departure and arrival of RV and PV commuters, and the changes in queue length at bottlenecks over time in Scenario 5. In Fig. \ref{fig:bieffects_S358}(b), the RV departure rate (i.e., slopes of blue solid line) decreases after time $t_0^\mathcal{P}$ and drops to zero after the last RV commuter departs at time $t_1^\mathcal{R}$. Because RV commuters experience queuing delays at the highway bottleneck before reaching the curbside bottleneck $B_\mathcal{C}$, their arrival rate at the curbside bottleneck $B_\mathcal{C}$ (i.e., slopes of black dashed lines) decreases at times $t_0^\mathcal{P}$ and $t_1^\mathcal{R} + w_\mathcal{H}(t_1^\mathcal{R})$ (see points A and B in Fig. \ref{fig:bieffects_S358}(b)). In addition, the RV arrival rate at the workplace (i.e., slopes of red dash-dotted line) decreases from $s_\mathcal{C}^\mathcal{R}$ to $\widetilde{s}_\mathcal{C}^\mathcal{R}$ at time $t_0^\mathcal{P} + w_\mathcal{C}^\mathcal{R}(t_0^\mathcal{P})$ (see point C in Fig. \ref{fig:bieffects_S358}(b)). Consequently, the RV queue length at the curbside bottleneck changes at the corresponding time points, as shown at points A’, B’, and C’ in Fig. \ref{fig:bieffects_S358}(b’). 
Notably, the observed decrease in the dissipation rate of RV queue length at point C’ in Fig. \ref{fig:bieffects_S358}(b') is due to the decrease in the RV arrival rate at point C in Fig. \ref{fig:bieffects_S358}(b). This phenomenon does not occur in the unidirectional congestion spillover case (see Fig. \ref{fig:eight_images}(e)).

Fig. \ref{fig:bieffects_S358}(c) and (c') indicate, respectively, the cumulative departure and arrival of RV and PV commuters, and the changes in queue length at bottlenecks over time in Scenario 8. Because the queue forms at the highway bottleneck from the beginning of the morning peak, the black dashed line shows that the RV arrival rate at the curbside bottleneck decreases at $t_0^\mathcal{P} +w_\mathcal{H}(t_0^\mathcal{P})$ and reaches zero at $t_1^\mathcal{R} + w_\mathcal{H}(t_1^\mathcal{R})$ (see points A and B in Fig. \ref{fig:bieffects_S358}(c)). In addition, the RV arrival rate at the workplace decreases from $s_\mathcal{C}^\mathcal{R}$ to $\widetilde{s}_\mathcal{C}^\mathcal{R}$ at time $t_0^\mathcal{P} + w_\mathcal{H}(t_0^\mathcal{P}) + w_\mathcal{C}^\mathcal{R}(t_0^\mathcal{P}+w_\mathcal{H}(t_0^\mathcal{P}))$, as shown at point C in Fig. \ref{fig:bieffects_S358}(c). Consequently, the RV queue length at the curbside bottleneck changes at the corresponding time points, as shown at points A’, B’, and C’ in Fig. \ref{fig:bieffects_S358}(c’). Obviously, the observed decrease in the dissipation rate of the RV queue length at point C’ in Fig. \ref{fig:bieffects_S358}(c') does not occur in the unidirectional congestion spillover case (see Fig. \ref{fig:eight_images}(h)).

\section{Congestion pricing}\label{section_congestion_pricing}
In the previous section, we discussed user equilibrium scenarios under unidirectional and bidirectional congestion spillover effects without congestion tolls. In the no-toll cases, two types of congestion externalities arise during the morning peak: the queue of RVs and PVs at the highway and curbside bottlenecks, as well as the congestion spillover effects among RVs and PVs. Due to the tendency of commuters to disregard the adverse externalities they generate, the entire transportation system experiences inefficiencies. Previous studies have demonstrated that the congestion externalities at bottlenecks can be eliminated by implementing dynamic time-varying tolls \citep{vickrey1969,arnott1990TransportationResearchPartB:Methodological,arnott1990JournalofUrbanEconomicsa}. Therefore, in this section, we investigate congestion pricing schemes that aim to eliminate congestion externalities at bottlenecks during the morning peak. It is important to note that this paper primarily focuses on optimal congestion pricing in bi-modal utilization scenarios. The optimal congestion pricing for single-mode utilization scenarios, including Scenario 1 and Scenario 6, is omitted as it is in line with previous research \citep{vickrey1969,arnott1990TransportationResearchPartB:Methodological,arnott1990JournalofUrbanEconomicsa}. For a comprehensive review of the bottleneck model and congestion pricing, interested readers can refer to \cite{li&huang2020TransportationResearchPartB:Methodological}.

In previous studies featuring a single bottleneck \citep{li&lam2017TransportationResearchPartB:Methodological,ren2020Transp.Policy, yu2022TransportationResearchPartC:EmergingTechnologies, wu2023TransportationResearchPartC:EmergingTechnologies}, congestion pricing is typically applied at the entrance of the highway. The congestion toll levied corresponds to the vehicle passage time, and the toll level is equivalent to the queuing delay cost in the no-toll case. However, in the bi-modal two-tandem bottleneck model, two traffic modes share the road both temporally and spatially, and queues form at bottlenecks on both the highway and the curb space. In addition, the congestion spillover effects further exacerbate the congestion delay experienced by RV and PV commuters. These disparities suggest that the conventional pricing scheme based solely on vehicle passage time at highway entrances may not be suitable for this intricate bi-modal two-tandem bottleneck scenario. Consequently, it is imperative to develop alternative pricing schemes that can effectively address the unique characteristics and challenges of the proposed model. 

One possible solution to this problem is to design different congestion pricing strategies for RVs and PVs, respectively. Taking into account the parking needs of PVs upon their arrival at the CBD, parking pricing can be introduced to alleviate PV congestion. By contrast, for RVs, their drop-offs can be viewed as transient parking in the curb space \citep{schaller2011,LiuXiaoHui2023TS}. The introduction of curbside pricing, similar to parking pricing, can address the external congestion caused by RVs\footnote{Charging for curb space usage by ride-hailing vehicles has attracted the attention of researchers and has been applied in real-world scenarios. For instance, \cite{LiuJiaChao2023TRC} optimized the pricing for each ride-hailing temporary stop (pick-ups and drop-offs) to minimize system costs. \cite{lim&masoud2024} considered multiple curb space users, including passenger pick-ups and drop-offs, and examined curb space allocation and pricing from the perspective of curb space operators. Additionally, the Council of the city of Vancouver enacted a by-law mandating that only ride-hailing providers holding a permit can pick up or drop off passengers within the Metro core between 7:00 am to 10:00 pm. A fee of \$0.25 to \$0.50 is charged for each pick-up or drop-off \citep{Vancouver2024}.}. On top of this, an optimal congestion pricing scheme is proposed, which integrates curbside pricing, charging RVs for short-term curbside usage, and parking pricing, charging PVs for parking upon arrival at the CBD. Different from congestion pricing at highway entrances, the combined pricing scheme investigated in this paper charges at both curbside areas and CBD parking lots. It can effectively eliminate the queuing delays at highway and curbside bottlenecks simultaneously, leading to the social optimal. In reality, this combined congestion pricing scheme can be implemented using on-street meters as well as mobile Apps (e.g., HKeMeter).

\subsection{Optimal congestion pricing} 
The optimal congestion pricing scheme can completely eliminate the queuing delays and minimize the total social cost of the system, leading to the social optimum. Here, we provide an intuitive derivation of the social optimum. First, commuters face the congestion fee instead of the queuing delays and choose their departure times from home (i.e., arrival times at CBD since there is no queue) to minimize their generalized travel cost. Second, the arrival intervals of commuters should be continuous and close to $t^*$; otherwise, the total costs of schedule delay are unnecessarily high. These two results together imply that bottlenecks should be fully utilized throughout rush hour \citep{arnott1990JournalofUrbanEconomicsa}. 

It is imperative to note that queues may inevitably occur at one bottleneck when the other bottleneck reaches its full capacity in a two-tandem bottleneck corridor. For example, if the capacity of the curbside bottlenecks exceeds that of the highway bottleneck (i.e., $s_\mathcal{C}^\mathcal{R}+s_\mathcal{C}^\mathcal{P}>s_\mathcal{H}$), full utilization of the curbside bottlenecks will inevitably lead to a queue at the highway bottleneck. Conversely, if $s_\mathcal{C}^\mathcal{R}+s_\mathcal{C}^\mathcal{P}\leq s_\mathcal{H}$, both RV and PV commuters can leave home and reach CBD at their respective maximum service rates at the curbside bottlenecks (i.e., $s_\mathcal{C}^\mathcal{R}$ and $s_\mathcal{C}^\mathcal{P}$ respectively) without causing any queues in the system, even during the co-departure phase. This underscores the imperative to comprehensively consider the attributes of the transport corridor, such as the service rates of bottlenecks, when designing congestion pricing schemes. Therefore, we will discuss the above two different cases, which are distinguished by the relative size of the bottlenecks. 
 
\subsubsection*{(1) the case of $s_\mathcal{C}^\mathcal{R}+s_\mathcal{C}^\mathcal{P}\leq s_\mathcal{H}$}

Under the optimal congestion pricing scheme, all commuters face no queue, i.e., $w_\mathcal{H}(t)=w_\mathcal{C}^\mathcal{R}(t)=w_\mathcal{C}^\mathcal{P}(t)=0$. Therefore, the generalized travel cost for commuters only consists of the schedule delay cost, the fixed cost, and the congestion fee. Let the time-varying congestion fee for RVs and PVs be $f^\mathcal{R}(t)$ and $f^\mathcal{P}(t)$, respectively. Let $\hat{C}^\mathcal{R}(t)$ and $\hat{C}^\mathcal{P}(t)$ be the generalized travel cost of RV and PV commuters under the optimal congestion pricing scheme. According to Eqs. \eqref{eq:costRV} and \eqref{eq:costPV}, $\hat{C}^\mathcal{R}(t)$ and $\hat{C}^\mathcal{P}(t)$ can be expressed as:
\begin{equation}\label{eq:cost_toll}
    \begin{cases}
        \hat{C}^\mathcal{R}(t)= \beta (t^*-t)+c^\mathcal{R} + f^\mathcal{R}(t), t\in[\hat{t}_0^\mathcal{R},t^*],\\
        \hat{C}^\mathcal{P}(t)=\beta (t^*-t)+u^\mathcal{P} + f^\mathcal{P}(t), t\in[\hat{t}_0^\mathcal{P},t^*],
    \end{cases}
\end{equation}
where $\hat{t}_0^\mathcal{R}$ and $\hat{t}_0^\mathcal{P}$ denote the arrival times of the first RV and the first PV after charging the congestion fee.

All commuters have an equal generalized travel cost under the optimal congestion pricing scheme. Set $\frac{\partial\hat{C}^{\mathcal{R}}(t)}{\partial t}=\frac{\partial\hat{C}^{\mathcal{P}}(t)}{\partial t}=0$, one can obtain,
\begin{equation}\label{eq:df/dt_uc1+ur1<u2}
    \frac{\partial f^{\mathcal{R}}(t)}{\partial t}=\frac{\partial f^{\mathcal{P}}(t)}{\partial t}=\beta.
\end{equation}

Eq. \eqref{eq:df/dt_uc1+ur1<u2} shows that the congestion fees $f^\mathcal{R}(t)$ and $f^\mathcal{P}(t)$ are linear functions with respect to the arrival time $t$. Based on this, by specifying the level of initial fees $f^\mathcal{R}(\hat{t}_0^\mathcal{R})$ and $f^\mathcal{P}(\hat{t}_0^\mathcal{P})$, the optimal congestion pricing functions can be determined. Therefore, the optimization of congestion pricing can be formulated as the following minimization problem for the social cost $SC$:
\begin{alignat}{2}
\mathop{\min}_{f^\mathcal{R}(\hat{t}_0^\mathcal{R}),f^\mathcal{P}(\hat{t}_0^\mathcal{P})} SC &=\int_{\hat{t}_0^\mathcal{R}}^{t^*} \big(\hat{C}^\mathcal{R}(t)-f^\mathcal{R}(t)\big) s_\mathcal{C}^\mathcal{R} dt+\int_{\hat{t}_0^\mathcal{P}}^{t^*} \big(\hat{C}^\mathcal{P}(t)-f^\mathcal{P}(t)\big) s_\mathcal{C}^\mathcal{P} dt,\label{eq:minimizationSC_ur1+uc1<u2}\\
\mbox{s.t.}\quad
&\quad\beta(t^*-\hat{t}_0^\mathcal{R})+f^\mathcal{R}(\hat{t}_0^\mathcal{R})+c^\mathcal{R} = \beta(t^*-\hat{t}_0^\mathcal{P})+f^\mathcal{P}(\hat{t}_0^\mathcal{P})+u^\mathcal{P},\label{eq:constrain1b} \\
&\quad(t^*-\hat{t}_0^\mathcal{R})s_\mathcal{C}^\mathcal{R} = N^\mathcal{R},\label{eq:constrain1c} \\
&\quad(t^*-\hat{t}_0^\mathcal{P})s_\mathcal{C}^\mathcal{P} = N^\mathcal{P},\label{eq:constrain1d} \\
&\quad N^\mathcal{R}+N^\mathcal{P} = N, \label{eq:constrain1e}
\end{alignat}
where Eq. \eqref{eq:constrain1b} indicates that the first RV commuter and the first PV commuter have an equal equilibrium generalized travel cost. Eqs. \eqref{eq:constrain1c}-\eqref{eq:constrain1e} are population constraints.

By solving Eqs. \eqref{eq:constrain1b}-\eqref{eq:constrain1e} and substituting the results of $\hat{t}_0^\mathcal{R}$ and $\hat{t}_0^\mathcal{P}$ into Eq. \eqref{eq:minimizationSC_ur1+uc1<u2}, then setting $\frac{\partial SC}{\partial f^\mathcal{R}(\hat{t}_0^\mathcal{R})} = \frac{\partial SC}{\partial f^\mathcal{P}(\hat{t}_0^\mathcal{P})} = 0$, one can obtain the first-order optimality condition for this minimization problem as\footnote{It can be proven that the Hessian matrix of Eq. \eqref{eq:minimizationSC_ur1+uc1<u2} is semi-positive definite, indicating that the objective function is convex with respective to ${f^\mathcal{R}(\hat{t}_0^\mathcal{R})}$ and ${f^\mathcal{P}(\hat{t}_0^\mathcal{P})}$. Therefore, the optimal solution to the minimization problem can be determined by solving the first-order optimality conditions.}:
\begin{equation}\label{eq:minimizationSC_ur1+uc1<u2_first_order_condition}
    f^\mathcal{R}(\hat{t}_0^\mathcal{R})=f^\mathcal{P}(\hat{t}_0^\mathcal{P}).
\end{equation}

Eq. \eqref{eq:minimizationSC_ur1+uc1<u2_first_order_condition} indicates an identical initial fee for RV and PV commuters. Let $f_0$ denote an arbitrary initial fee, i.e., $f^\mathcal{R}(\hat{t}_0^\mathcal{R})=f^\mathcal{P}(\hat{t}_0^\mathcal{P}) = f_0$. Combining it with Eq. \eqref{eq:df/dt_uc1+ur1<u2} yields the optimal congestion pricing for this case:
\begin{equation}\label{eq:toll_level_ur1+uc1<u2}
    \begin{aligned}
        &f^\mathcal{R}(t)=
            \begin{cases}
            f_0,\ & t<\hat{t}_0^\mathcal{R},\\[1.5ex]
            f_0+\beta(t-\hat{t}_0^\mathcal{R}),\ &\hat{t}_0^\mathcal{R}\leq t\leq t^*.\\
            \end{cases}\\
        &f^\mathcal{P}(t)=
            \begin{cases}
            f_0,\ &t<\hat{t}_0^\mathcal{P},\\[1.5ex]
            f_0+\beta(t-\hat{t}_0^\mathcal{P}),\ &\hat{t}_0^\mathcal{P}\leq t\leq t^*.\\
            \end{cases}\\
    \end{aligned}
\end{equation} 

It is important to clarify that $f_0$ represents an arbitrary initial fee. Since both modes of transportation (i.e., RV and PV) have the same initial fee, the value of $f_0$ does not affect the modal split and the performance of the optimal congestion pricing. 

Substitute Eqs. \eqref{eq:toll_level_ur1+uc1<u2} into Eqs. \eqref{eq:constrain1b}-\eqref{eq:constrain1e}, one can further obtain the equilibrium results with optimal congestion pricing (i.e., the social optimum) as:
\begin{equation}\label{eq:ueresults_toll_ur1+uc1<u2}
    \begin{cases}
        \hat{t}_0^\mathcal{R} = t^*-\frac{N\beta+(u^\mathcal{P}-c^\mathcal{R})s_\mathcal{C}^\mathcal{P}}{\beta(s_\mathcal{C}^\mathcal{R}+s_\mathcal{C}^\mathcal{P})},\\[1.5ex]
        \hat{t}_0^\mathcal{P} = t^*-\frac{N\beta-(u^\mathcal{P}-c^\mathcal{R})s_\mathcal{C}^\mathcal{R}}{\beta(s_\mathcal{C}^\mathcal{R}+s_\mathcal{C}^\mathcal{P})},\\[1.5ex]
        \hat{N}^\mathcal{R}=\frac{N\beta+(u^\mathcal{P}-c^\mathcal{R})s_\mathcal{C}^\mathcal{P}}{\beta(s_\mathcal{C}^\mathcal{P}+s_\mathcal{C}^\mathcal{R})}s_\mathcal{C}^\mathcal{R},\\[1.5ex]
        \hat{C}^\mathcal{R}(t)=\hat{C}^\mathcal{P}(t) = \frac{N\beta+u^\mathcal{P}s_{C}^\mathcal{P}+c^\mathcal{R}s_\mathcal{C}^\mathcal{R}}{s_\mathcal{C}^\mathcal{R}+s_\mathcal{C}^\mathcal{P}}+f_0.
    \end{cases}
\end{equation}


\subsubsection*{(2) the case of $s_\mathcal{C}^\mathcal{R}+s_\mathcal{C}^\mathcal{P} >s_\mathcal{H}$}
In the case of $s_\mathcal{C}^\mathcal{R}+s_\mathcal{C}^\mathcal{P} >s_\mathcal{H}$, it is crucial to ensure that both the individual departure rates of RV and PV commuters, as well as their combined departure rate during the co-departure stage, do not exceed the service rates of the bottlenecks. To minimize the overall social cost without creating queues, for the individual travel stage of RVs or PVs, the departure rate of commuters should be equal to $s_\mathcal{C}^\mathcal{R}$ or $s_\mathcal{C}^\mathcal{P}$, respectively. For the co-departure stage, the sum of their departure rates should be equal to $s_\mathcal{H}$. For the sake of illustration, we introduce a parameter $\theta$ to denote the ratio of the departure rate of PVs to the curbside bottleneck capacity during the co-departure stage (referred to as the share ratio), i.e., $\theta = \dot {A}_\mathcal{C}^\mathcal{P}(t)/s_\mathcal{C}^\mathcal{P}$. Consequently, during the co-departure stage, the departure rates of PV and RV commuters can be expressed as $\theta s_\mathcal{C}^\mathcal{P}$ and $s_\mathcal{H}-\theta s_\mathcal{C}^\mathcal{P}$, respectively. The parameter $\theta$ falls within the range of $\big[(s_\mathcal{H}-s_\mathcal{C}^\mathcal{R})/s_\mathcal{C}^\mathcal{P},1\big]$, ensuring full utilization of both bottlenecks without the occurrence of queues. 

Similarly, the generalized travel cost of RV and PV commuters under the optimal congestion pricing scheme in this case can be expressed as Eq. \eqref{eq:cost_toll}. Since all commuters have an equal generalized travel cost, setting $\frac{\partial\hat{C}^{\mathcal{R}}(t)}{\partial t}=\frac{\partial\hat{C}^{\mathcal{P}}(t)}{\partial t}=0$, one can obtain,
\begin{equation}\label{eq:df/dt_uc1+ur1>u2}
    \frac{\partial f^{\mathcal{R}}(t)}{\partial t}=\frac{\partial f^{\mathcal{P}}(t)}{\partial t}=\beta.
\end{equation}

Eq. \eqref{eq:df/dt_uc1+ur1>u2} shows that the congestion fees $f^\mathcal{R}(t)$ and $f^\mathcal{P}(t)$ are linear functions with respect to the arrival time. Based on this, by further specifying the level of initial fees $f^\mathcal{R}(\hat{t}_0^\mathcal{R})$ and $f^\mathcal{P}(\hat{t}_0^\mathcal{P})$, the optimal congestion pricing functions can be determined. Therefore, the optimization for congestion pricing in this case can be formulated as the following social cost minimization problem.
\begin{alignat}{2}
    \begin{split}
    \mathop{\min}_{f^\mathcal{R}(\hat{t}_0^\mathcal{R}),f^\mathcal{P}(\hat{t}_0^\mathcal{P})} SC &=\int_{\hat{t}_0^\mathcal{R}}^{\hat{t}_0^\mathcal{P}} \big(\hat{C}^\mathcal{R}(t)-f^\mathcal{R}(t)\big) s_\mathcal{C}^\mathcal{R} dt + \int_{\hat{t}_0^\mathcal{P}}^{t^*}\big(\hat{C}^\mathcal{R}(t)-f^\mathcal{R}(t)\big) (s_\mathcal{H}-\theta s_{C}^\mathcal{P}) dt\\&+\int_{\hat{t}_0^\mathcal{P}}^{t^*} \big(\hat{C}^\mathcal{P}(t)-f^{P}(t)\big) \theta s_\mathcal{C}^\mathcal{P} dt,
    \end{split}\label{eq:minimizationSC_ur1+uc1>u2}\\
\mbox{s.t.}\quad
    &\quad\beta(t^*-\hat{t}_0^\mathcal{R})+f^\mathcal{R}(\hat{t}_0^\mathcal{R})+c^\mathcal{R} = \beta(t^*-\hat{t}_0^\mathcal{P})+f^\mathcal{P}(\hat{t}_0^\mathcal{P})+u^\mathcal{P}\label{eq:constrain2b}, \\
    &\quad(\hat{t}_0^\mathcal{P}-\hat{t}_0^\mathcal{R})s_\mathcal{C}^\mathcal{R} + (t^*-\hat{t}_0^\mathcal{P})(s_\mathcal{H}-\theta s_\mathcal{C}^\mathcal{P})=N^\mathcal{R}\label{eq:constrain2c}, \\
    &\quad(t^*-\hat{t}_0^\mathcal{P})\theta s_\mathcal{C}^\mathcal{P} = N^\mathcal{P}\label{eq:constrain2d}, \\
    &\quad N^\mathcal{R}+N^\mathcal{P} = N \label{eq:constrain2e},
\end{alignat}
where Eq. \eqref{eq:constrain2b} indicates that the first RV commuter and the first PV commuter have an equal generalized travel cost. Eqs. \eqref{eq:constrain2c}-\eqref{eq:constrain2e} are population constraints. 

Solving Eqs. \eqref{eq:constrain2b}-\eqref{eq:constrain2e} and substituting the results of $\hat{t}_0^\mathcal{R}$ and $\hat{t}_0^\mathcal{P}$ into Eq. \eqref{eq:minimizationSC_ur1+uc1>u2}, one can obtain,
\begin{equation}\label{eq:SC_ur1+uc1>u2}
    \begin{aligned}
        SC &= \frac{s_\mathcal{C}^\mathcal{R}}{2\beta s_\mathcal{H}}\bigg( \big(f^\mathcal{R}(\hat{t}_0^\mathcal{R})-f^\mathcal{P}(\hat{t}_0^\mathcal{P})-u^\mathcal{P}+c^\mathcal{R}\big) \big(2\theta s_\mathcal{C}^\mathcal{P}(u^\mathcal{P}-c^\mathcal{R})+(s_\mathcal{H}-s_\mathcal{C}^\mathcal{R})(f^\mathcal{R}(\hat{t}_0^\mathcal{R})-f^\mathcal{P}(\hat{t}_0^\mathcal{P})-u^\mathcal{P}+c^\mathcal{R}) \big)   \bigg)\\
        &+\frac{N}{2s_\mathcal{H}}\big(N\beta+2s_{\mathcal{H}} C^{\mathcal{R}}+2\theta s_\mathcal{C}^\mathcal{P}(u^\mathcal{P}-c^\mathcal{R})\big).
    \end{aligned}
\end{equation}

setting $\frac{\partial{SC}}{\partial{f^\mathcal{R}(\hat{t}_0^\mathcal{R})}}=\frac{\partial{SC}}{\partial{f^\mathcal{P}(\hat{t}_0^\mathcal{P})}}=0$, one can obtain the first-order optimality condition for the objective function in Eq. \eqref{eq:minimizationSC_ur1+uc1>u2} as\footnote{The Hessian matrix of Eq. \eqref{eq:SC_ur1+uc1>u2} is semi-positive definite, indicating that the objective function is convex with respective to ${f^\mathcal{R}(\hat{t}_0^\mathcal{R})}$ and ${f^\mathcal{P}(\hat{t}_0^\mathcal{P})}$.}: 
\begin{equation}\label{eq:firstoptimalitycondition}
    f^\mathcal{P}(\hat{t}_0^\mathcal{P})-f^\mathcal{R}(\hat{t}_0^\mathcal{R}) =\frac{(u^\mathcal{P}-c^\mathcal{R})(s_\mathcal{C}^\mathcal{R}+\theta s_\mathcal{C}^\mathcal{P}-s_\mathcal{H})}{s_\mathcal{H}-s_\mathcal{C}^\mathcal{R}}.
\end{equation}

For presentation purpose, define $\Delta f = f^\mathcal{P}(\hat{t}_0^\mathcal{P})-f^\mathcal{R}(\hat{t}_0^\mathcal{R})$, which represents the difference in the initial fee between PV and RV commuters. Eq. \eqref{eq:firstoptimalitycondition} indicates that $\Delta f$ has effects on the share ratio $\theta$. As $\theta\in\big[{(s_\mathcal{H}-s_\mathcal{C}^\mathcal{R})}/{s_\mathcal{C}^\mathcal{P}},1\big]$, ensuring full utilization of both bottlenecks without the occurrence of queues, the value of $\Delta f$ should fall within the range of $\big[0,{(u^\mathcal{P}-c^\mathcal{R})(s_\mathcal{C}^\mathcal{R}+ s_\mathcal{C}^\mathcal{P}-s_\mathcal{H})}/{(s_\mathcal{H}-s_\mathcal{C}^\mathcal{R})}\big]$. Substituting Eq. \eqref{eq:firstoptimalitycondition} into Eq. \eqref{eq:SC_ur1+uc1>u2}, and taking the first-order derivative of $SC$ with respective to $\Delta f$, we have\footnote{In summary, $\frac{\partial SC}{\partial \Delta f} \geq 0$ for $\Delta f \leq \frac{\beta N}{s_\mathcal{C}^\mathcal{R}} - (u^\mathcal{P} - c^\mathcal{R})$, and $\frac{\partial SC}{\partial \Delta f} < 0$ for $\Delta f > \frac{\beta N}{s_\mathcal{C}^\mathcal{R}} - (u^\mathcal{P} - c^\mathcal{R})$. However, when $\Delta f > \frac{\beta N}{s_\mathcal{C}^\mathcal{R}} - (u^\mathcal{P} - c^\mathcal{R})$, the cost difference between two modes becomes so large that $(u^\mathcal{P} - c^\mathcal{R}) + \Delta f > \frac{\beta N}{s_\mathcal{C}^\mathcal{R}}$ holds, leading all commuters to choose the RV mode (see Proposition \ref{Proposition 2}). Therefore, the condition $\frac{\partial SC}{\partial \Delta f} \geq 0$ always holds for the bi-modal utilization scenarios.},
\begin{equation}\label{eq:dSC/df_Delta}
    \frac{\partial SC}{\partial \Delta f}\geq 0.
\end{equation}

Eq. \eqref{eq:dSC/df_Delta} indicates that the optimal $\Delta f^*$ to minimize the social cost takes the value of its left boundary, i.e., $\Delta f^* = 0$. It implies an identical initial fee of the optimal congestion pricing scheme for PV and RV commuters. Let $f_0$ represent an arbitrary initial fee, i.e., $f^\mathcal{R}(\hat{t}_0^\mathcal{R})=f^\mathcal{P}(\hat{t}_0^\mathcal{P})=f_0$. Combining it with Eq. \eqref{eq:df/dt_uc1+ur1>u2} yields the optimal congestion pricing for this case:
\begin{equation}\label{eq:toll_level_ur1+uc1>u2}
    \begin{aligned}
        &f^\mathcal{R}(t)=
            \begin{cases}
            f_0,\ &t<\hat{t}_0^\mathcal{R},\\[1.5ex]
            f_0+\beta(t-\hat{t}_0^\mathcal{R}),\ &\hat{t}_0^\mathcal{R}\leq t\leq t^*.\\
            \end{cases}\\
        &f^\mathcal{P}(t)=
            \begin{cases}
            f_0,\ &t<\hat{t}_0^\mathcal{P},\\[1.5ex]
            f_0+\beta(t-\hat{t}_0^\mathcal{P}),\ &\hat{t}_0^\mathcal{P}\leq t\leq t^*.\\
            \end{cases}\\
    \end{aligned}
\end{equation}Similarly, because both modes of transportation (i.e., RV and PV) have the same initial fee, the value of $f_0$ does not affect the modal split and the performance of the optimal congestion pricing.

Combining Eqs. \eqref{eq:toll_level_ur1+uc1>u2} with Eqs. \eqref{eq:constrain2b}-\eqref{eq:constrain2e} and \eqref{eq:firstoptimalitycondition}, one can further derive the equilibrium results with optimal congestion pricing (i.e., the social optimum) for this case as:
\begin{equation}\label{eq:ueresults_toll_ur1+uc1>u2}
    \begin{cases}
        \hat{t}_0^\mathcal{R} = t^*-\frac{N\beta+(u^\mathcal{P}-c^\mathcal{R})(s_\mathcal{H}-s_\mathcal{C}^\mathcal{R})}{\beta s_\mathcal{H}},\\[1ex]
        \hat{t}_0^\mathcal{P} = t^*-\frac{N\beta-(u^\mathcal{P}-c^\mathcal{R})s_\mathcal{C}^\mathcal{R}}{\beta s_\mathcal{H}},\\[1ex]
        \hat{N}^\mathcal{R}=\frac{\beta N + (u^\mathcal{P}-c^\mathcal{R})(s_\mathcal{H}-s_\mathcal{C}^\mathcal{R})}{\beta s_\mathcal{H}}s_\mathcal{C}^\mathcal{R},\\[1ex]
        \hat{C}^\mathcal{R}(t)=\hat{C}^\mathcal{P}(t) = \frac{N\beta+u^\mathcal{P}s_\mathcal{H}-(u^\mathcal{P}-c^\mathcal{R})s_\mathcal{C}^\mathcal{R}}{s_\mathcal{H}}+f_0.
    \end{cases}
\end{equation}

\subsection{Summary}
It is evident that while the equilibrium results of the social optimum differ between the cases of $s_\mathcal{C}^\mathcal{R}+s_\mathcal{C}^\mathcal{P}\leq s_\mathcal{H}$ and $s_\mathcal{C}^\mathcal{R}+s_\mathcal{C}^\mathcal{P}> s_\mathcal{H}$ (see Eqs. \eqref{eq:ueresults_toll_ur1+uc1<u2} and \eqref{eq:ueresults_toll_ur1+uc1>u2}), the optimal congestion pricing that achieves social optimum is structurally consistent, as defined by Eqs. \eqref{eq:toll_level_ur1+uc1<u2} and \eqref{eq:toll_level_ur1+uc1>u2}. Specifically, the optimal congestion pricing scheme applies a uniform initial fee $f_0$  (i.e., the fee paid by the first RV and the first PV commuters) and an identical marginal charge rate $\beta$ for both RVs and PVs. In addition, based on Eqs. \eqref{eq:toll_level_ur1+uc1<u2}-\eqref{eq:ueresults_toll_ur1+uc1<u2} and \eqref{eq:toll_level_ur1+uc1>u2}-\eqref{eq:ueresults_toll_ur1+uc1>u2}, for the co-departure stage of RV and PV commuters, we have,
\begin{equation}
    \label{eq:toll_level_diff_co_departure}
    f^\mathcal{R}(t)-f^\mathcal{P}(t) = \beta(\hat{t}_0^\mathcal{P}-\hat{t}_0^\mathcal{R})=u^\mathcal{P}-c^\mathcal{R},\ t\in[\hat{t}_0^\mathcal{P},t^*].
\end{equation}

\begin{figure}[!h]
    \centering
    \includegraphics[width=0.5\linewidth]{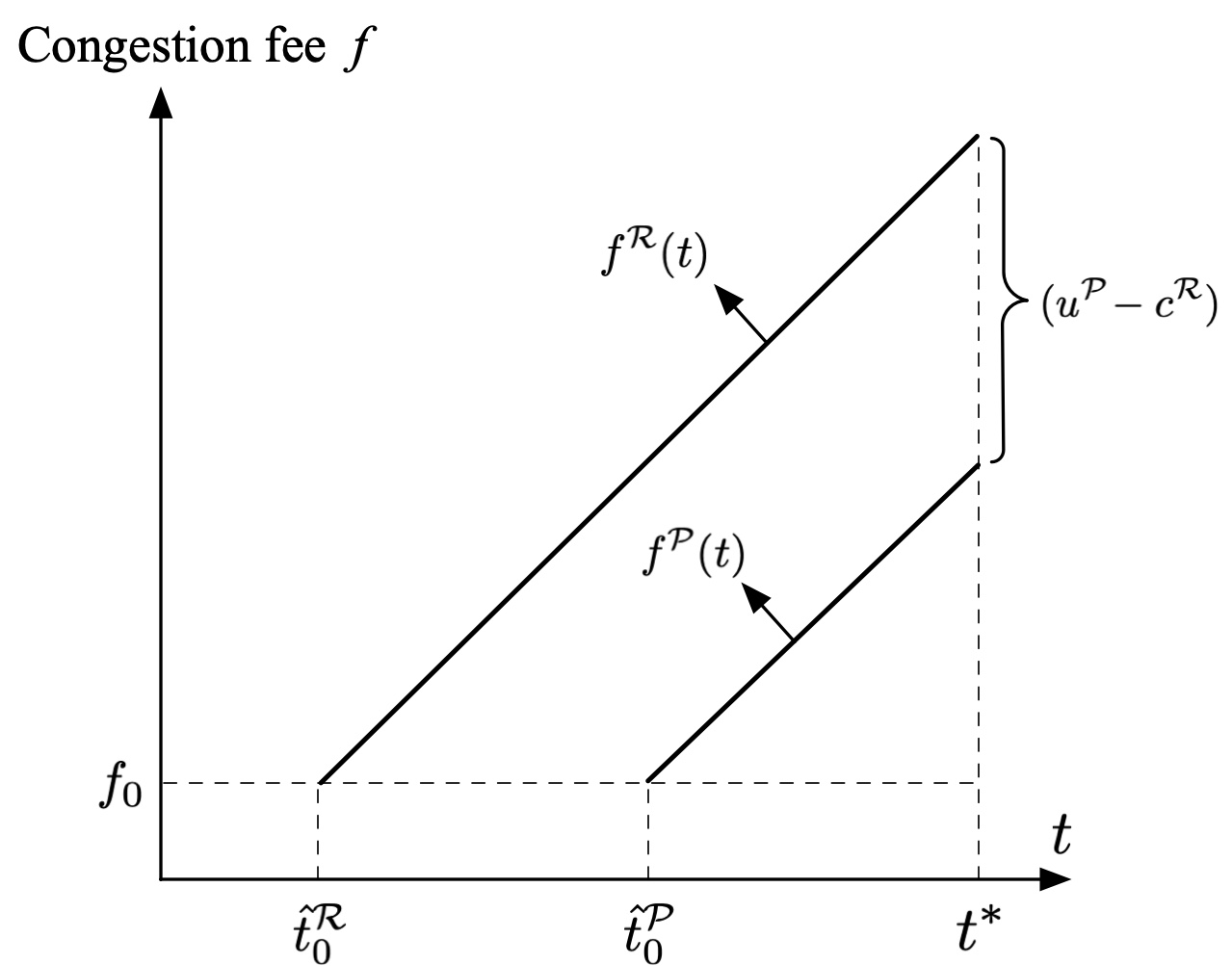}
    \caption{The optimal congestion pricing for RVs and PVs}
    \label{fig:congestion_fee}
\end{figure}

Eq. \eqref{eq:toll_level_diff_co_departure} indicates that the difference in congestion fees between RV and PV commuters during the co-departure phase is constant as $u^\mathcal{P}-c^\mathcal{R}$. This simplifies the practical implementation of the pricing scheme, allowing the RV to be charged with just a markup of $u^\mathcal{P}-c^\mathcal{R}$ over the PV. For illustration purposes, Fig. \ref{fig:congestion_fee} plots the optimal congestion fee with respect to time $t$. It is intriguing that despite the presence of multiple equilibrium scenarios shown in Fig. \ref{fig:eight_images}, the optimal congestion pricing scheme leading to social optimum remains remarkably straightforward. The following proposition summarizes the properties of the optimal congestion pricing scheme.
\begin{Proposition}[Optimal congestion pricing]\label{Proposition 4}
The optimal congestion pricing scheme has a uniform initial fee and an identical marginal charge rate for RVs and PVs. During the co-departure phase, the difference in the optimal congestion fee for RVs and PVs is constant as $u^\mathcal{P}-c^\mathcal{R}$.
\end{Proposition}

\section{Numerical illustration}\label{section_numerical_illustration}
In this section, we will illustrate the properties of the proposed model and its practical applications through two numerical examples. The first example is used to verify the eight equilibrium scenarios derived analytically in previous sections, as well as the performance of congestion pricing schemes. The second example is based on real data in Hong Kong, which analyzes equilibrium travel patterns during the morning peak commuting period in the study area and discusses feasible congestion mitigation strategies.

\subsection{A synthetic example}
We first use a synthetic example, as depicted in Fig. \ref{overview of the road network}, to illustrate the properties of the proposed model. The input data for the model parameters are set as follows. Following \cite{tang2021TransportationResearchPartE:LogisticsandTransportationReview} and \cite{xiao2021TransportationResearchPartB:Methodological}, the values of travel time $\alpha$ and schedule delay time $\beta$ are $\$ 6.4/h$ and $\$ 3.9/h$, respectively. The unit price for the trip time using RV is $\pi = \$ 8/h$. The total number of commuters is $N=3,000$. The service rate of the highway bottleneck $B_\mathcal{H}$ is 2,500 vehicles per hour, i.e., $s_\mathcal{H} = 2,500$veh/h. The parameters $\delta^\mathcal{R}$ and $\delta^\mathcal{P}$, which measure the intensity of spillover effects, are set to 0.1. The initial fee for the congestion pricing scheme is set to zero, i.e., $f_0=0$. In the following analysis, unless specifically stated otherwise, the input data are the same as the values detailed above.

\subsubsection{Equilibrium scenarios division}\label{numerical_equilibrium_scenarios}

According to Proposition \ref{Proposition 2} and Table \ref{table: conditions_for_8_scenarios}, the equilibrium scenarios will be determined by the value of $u^\mathcal{P}-c^\mathcal{R}$ and the capacities of the highway and curbside bottlenecks. In this section, we discuss how these parameters influence the occurrence of equilibrium scenarios. 

\begin{figure}[!ht]
    \centering
    \includegraphics[width=\linewidth]{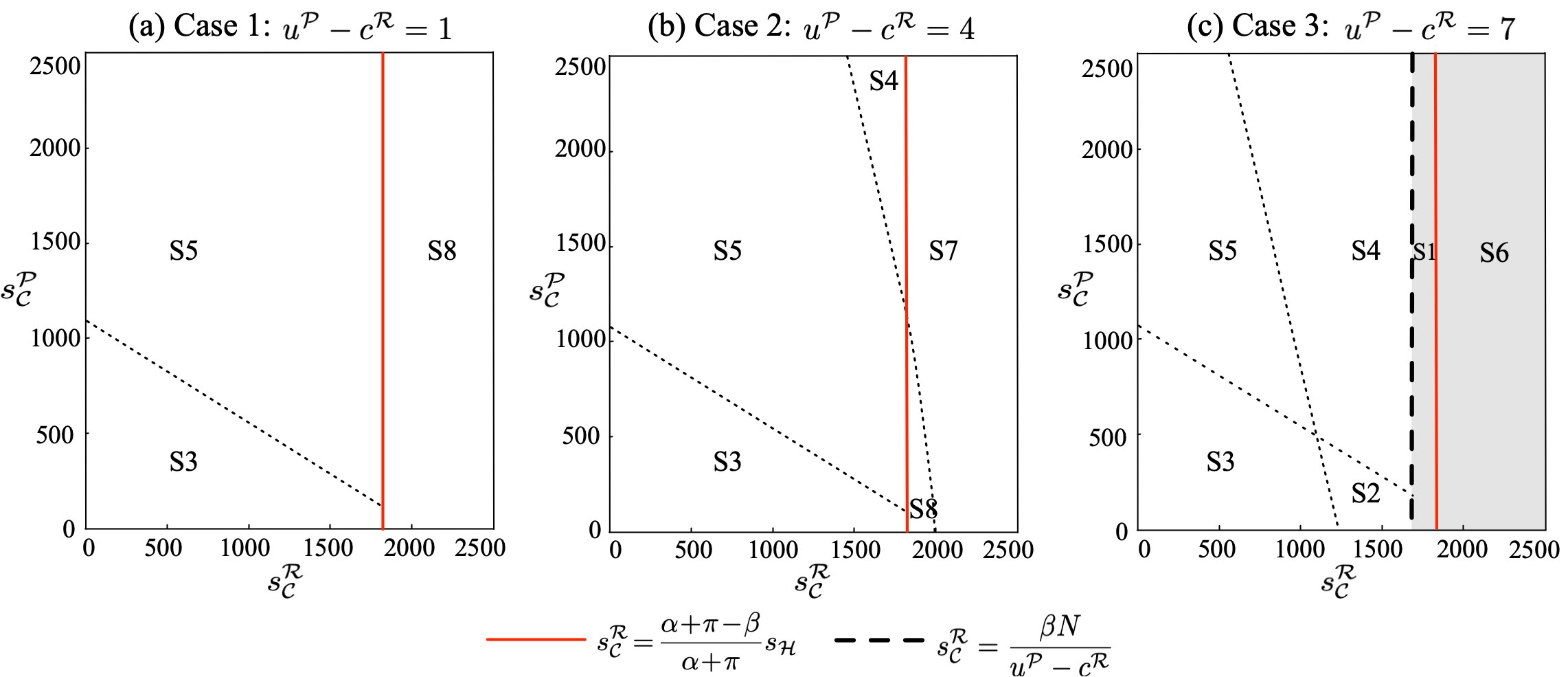}
    \caption{Equilibrium scenarios division with respect to $s_\mathcal{C}^\mathcal{R}$, $s_\mathcal{C}^\mathcal{P}$, and $u^\mathcal{P}-c^\mathcal{R}$}
    \label{fig:equilibrium_scenarios}
\end{figure}

Fig. \ref{fig:equilibrium_scenarios} illustrates the partitioning of different equilibrium scenarios on the two-dimensional space $(s_\mathcal{C}^\mathcal{R}, s_\mathcal{C}^\mathcal{P})$. Areas S1 to S8 correspond to Scenarios 1 to 8, respectively. Fig. \ref{fig:equilibrium_scenarios}(a)-(c) correspond to the following cases: (a) Case 1: the value of $u^\mathcal{P}-c^\mathcal{R}$ is low, thus both modes of transportation (i.e., RV and PV) are used with overlapping departure intervals, leading to Scenarios 3, 5, and 8; (b) Case 2: the value of $u^\mathcal{P}-c^\mathcal{R}$ is moderate, thus the scenarios where RV and PV are used with non-overlapping departure intervals appear (i.e., Scenarios 4 and 7); (c) Case 3: the value of $u^\mathcal{P}-c^\mathcal{R}$ is large enough that the scenarios where only RV is used appear (i.e., Scenarios 1 and 6). In the following, each case will be briefly introduced.

In Case 1, Scenarios 3, 5, and 8 emerge, corresponding to areas S3, S5 and S8 in Fig. \ref{fig:equilibrium_scenarios}(a). Both RV and PV will be used with overlapping departure intervals (referring to the features of Scenarios 3, 5, and 8 shown in Table \ref{table: 8 scenarios} and Fig. \ref{fig:eight_images}). If $s_\mathcal{C}^\mathcal{R}$ and $s_\mathcal{C}^\mathcal{P}$ are small enough (as indicated by the area S3 in Fig. \ref{fig:equilibrium_scenarios}(a)), then Scenario 3 occurs and queues form solely at the curbside bottlenecks. However, as $s_\mathcal{C}^\mathcal{R}$ and $s_\mathcal{C}^\mathcal{P}$ increase, Scenario 5 or Scenario 8 will occur and queues extend to the highway bottleneck. More specifically, as $s_\mathcal{C}^\mathcal{R}>\frac{\alpha+\pi-\beta}{\alpha+\pi}s_\mathcal{H}$, Scenario 8 occurs and queues form at both the highway and the curbside bottlenecks from the start of the morning peak.
    
In Case 2, compared to Case 1, Scenarios 4 and 7 emerge, corresponding to areas S4 and S7 in Fig. \ref{fig:equilibrium_scenarios}(b). Referring to Table \ref{table: 8 scenarios} and Fig. \ref{fig:eight_images}, both RV and PV will be utilized in Scenarios 4 and 7. Unlike Scenarios 3, 5, and 8, the departure intervals of RV and PV commuters do not overlap in Scenarios 4 and 7. For the queuing situation, except Scenario 3 where queues only form at the curbside bottleneck, queues also form at the highway in other scenarios. More specifically, in Scenarios 4 and 5, queues initially form solely at the curbside bottleneck, and later extend to the highway bottleneck after the first PV commuter departs from home. Conversely, in Scenarios 7 and 8, queues emerge at both the highway and the curbside bottlenecks from the start of the morning peak.

In Case 3, compared to Case 2, Scenarios 7 and 8 do not occur, while Scenarios 1, 2, and 6 emerge, corresponding to areas S1, S2, and S6 in Fig. \ref{fig:equilibrium_scenarios}(c). For the utilization of travel modes, RV and PV are used with overlapping departure intervals in Scenarios 3 and 5, while in Scenarios 2 and 4, RV and PV are used with non-overlapping departure intervals. When $s_\mathcal{C}^\mathcal{R}>\frac{\beta N}{(u^\mathcal{P}-c^\mathcal{R})}$, Scenario 1 and Scenario 6 will occur and only the RV mode will be used. For the queuing situation, in Scenarios 1, 2, and 3, queues form solely at the curbside bottleneck. In Scenarios 4 and 5, queues initially form solely at the curbside bottleneck, and later extend to the highway bottleneck after the first PV commuter departs from home. In Scenarios 6, queues emerge at both bottlenecks from the initial phase.

In general, Fig. \ref{fig:equilibrium_scenarios} indicates that as the service rates $s_\mathcal{C}^\mathcal{R}$ and $s_\mathcal{C}^\mathcal{P}$ increase, all three cases exhibit a tendency of queuing congestion extending from the curbside bottlenecks to the highway bottleneck. This suggests that the expansion of the curb space or the main road may result in a spillover of traffic congestion to the highway. However, such a shift in congestion from the curb space or main road to the highway does not necessarily result in an increase in the total social cost or traffic congestion across the entire system. On the contrary, as discussed later in Section \ref{subsection_num_Capacity}, the total social cost is expected to decrease as the capacity of the curb space or main road increases. Furthermore, the expansion of the curb space or the main road causes the equilibrium scenarios to shift from joint utilization of RV and PV with overlapping departure intervals (i.e., Scenarios 3, 5, and 8) to joint utilization with non-overlapping departure intervals (i.e., Scenarios 2, 4, and 7), and finally to exclusive utilization of RVs (i.e., Scenarios 1 and 6). 

\subsubsection{System performance}\label{numerical_system_performance}
\begin{figure}[!ht]
    \centering
    \includegraphics[width=0.95\linewidth]{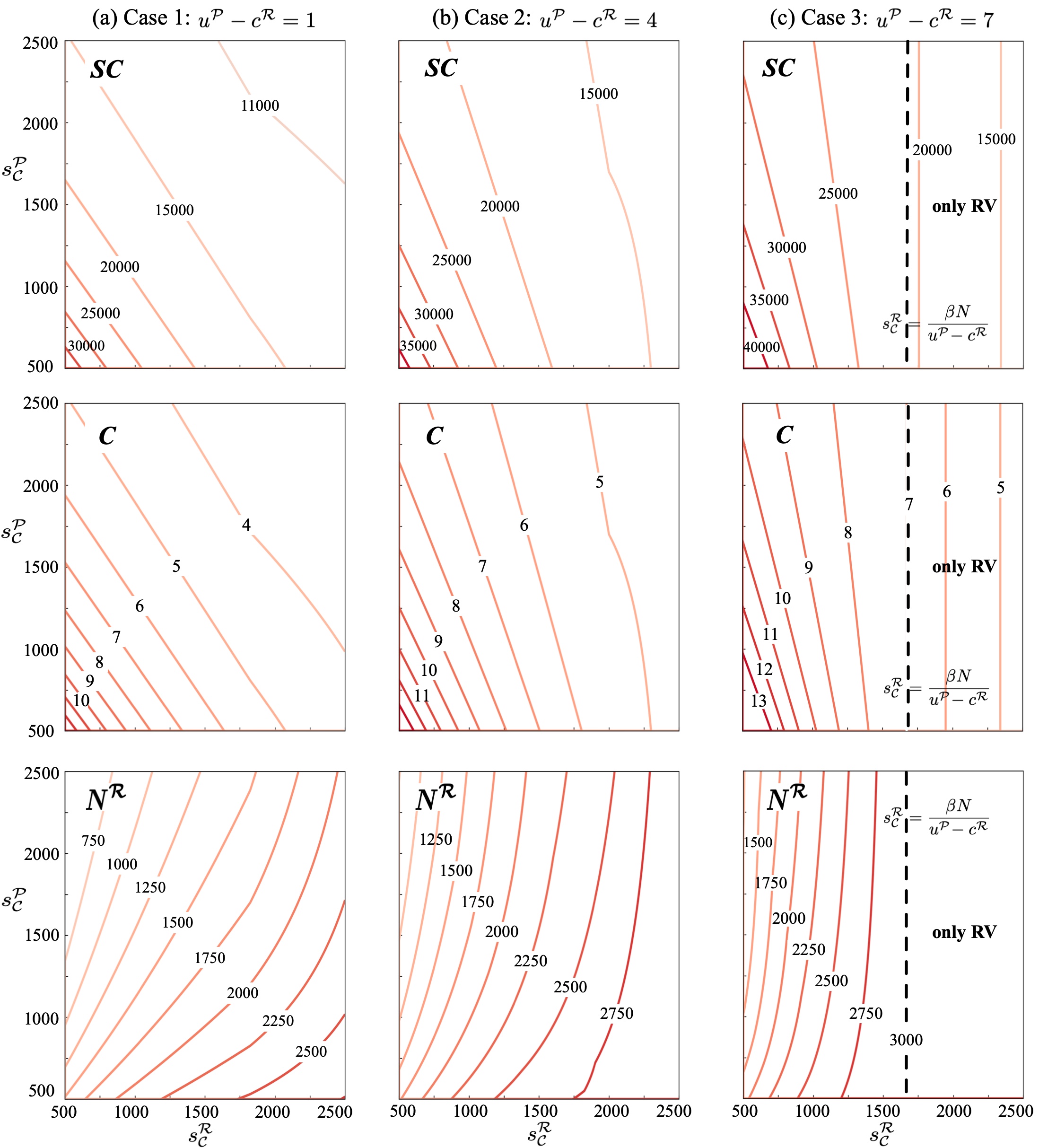}
    \caption{Contour maps for $SC$, $C$, and $N^\mathcal{R}$ with respect to $s_\mathcal{C}^\mathcal{R}$ and $s_\mathcal{C}^\mathcal{P}$}
    \label{fig:SC_C_Nr}
\end{figure}

Fig. \ref{fig:SC_C_Nr} indicates changes in social cost $SC$, individual travel cost $C$, and modal split $N^\mathcal{R}$ (i.e., the number of RV users) with respect to the service rates of the curbside bottlenecks $s_\mathcal{C}^\mathcal{R}$ and $s_\mathcal{C}^\mathcal{P}$. In general, as $s_\mathcal{C}^\mathcal{R}$ and $s_\mathcal{C}^\mathcal{P}$ increase, the social cost $SC$ and the individual travel cost $C$ decrease. The modal split $N^\mathcal{R}$ increases with an increase in $s_\mathcal{C}^\mathcal{R}$ or a decrease in $s_\mathcal{C}^\mathcal{P}$. Moreover, if the value of $u^\mathcal{P}-c^\mathcal{R}$ is larger than $\frac{\beta N}{s_\mathcal{C}^\mathcal{R}}$, only the RV mode will be used during the morning peak. In this case, both $SC$ and $C$ decrease as $s_\mathcal{C}^\mathcal{R}$ increases, whereas the increase in  $s_\mathcal{C}^\mathcal{P}$ has no impact on $SC$ or $C$. These findings suggest that the expansion of the curb space or the main road is capable of effectively reducing the overall social cost and the individual travel cost of commuters. 

Fig. \ref{fig:SC_C_Nr} also shows that as $s_\mathcal{C}^\mathcal{R}$ and $s_\mathcal{C}^\mathcal{P}$ increase, the contour lines of $SC$ and $C$ become more sparse, indicating that the marginal benefit of expanding the curb space or the main road is reduced. These findings are in line with the reality that larger curb spaces or wider main roads can reduce travel costs, but this benefit may decrease with increasing capacity.

\subsubsection{Optimal congestion pricing}
\begin{figure}[!ht]
    \centering
    \includegraphics[width=0.95\linewidth]{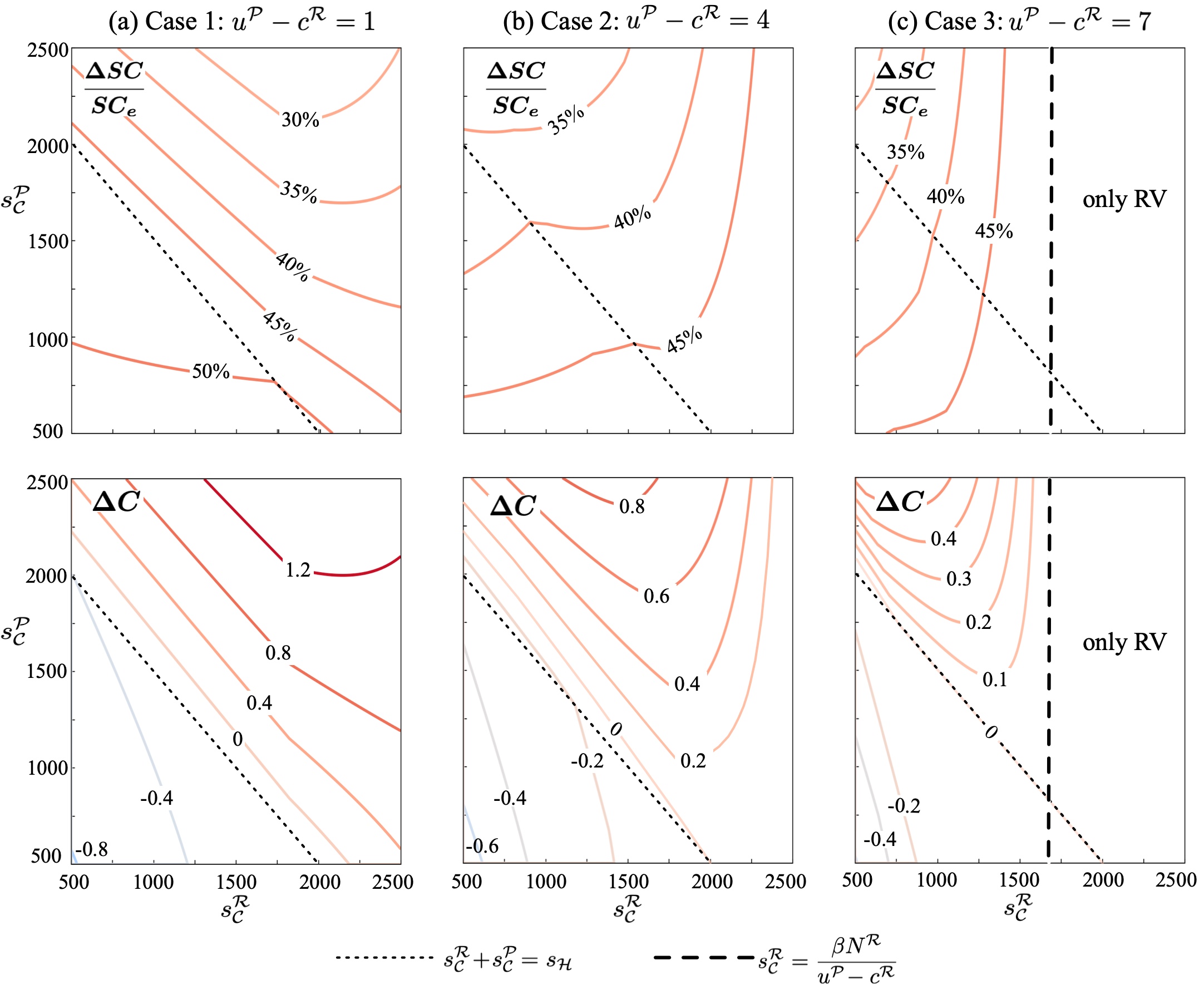}
    \caption{Contour maps for $\frac{\Delta SC}{SC^e}$ and $\Delta C$ with respect to $s_\mathcal{C}^\mathcal{R}$ and $s_\mathcal{C}^\mathcal{P}$}
    \label{fig:DeltaSC_DeltaC}
\end{figure}

Fig. \ref{fig:DeltaSC_DeltaC} illustrates the efficacy of the proposed optimal congestion pricing scheme. It is measured by the social cost reduction rate, denoted by $\frac{\Delta SC}{SC_e}=\frac{SC_e-SC_o}{SC_e}$, as well as the changes in the individual travel cost due to the optimal congestion pricing, denoted by $\Delta C=C_o-C_e$, where the subscripts ``o'' and ``e'' represent the social optimum and the no-toll equilibrium, respectively. 

It can be observed that, except the scenarios with only RV utilization, the social cost reduction rate $\frac{\Delta SC}{SC_e}$ is generally below 50\%. This implies that the performance of the time-varying congestion pricing scheme in our model is not as effective as in the traditional single bottleneck model \citep{vickrey1969,arnott1990TransportationResearchPartB:Methodological, arnott1990JournalofUrbanEconomicsa}. In other words, previous studies may overestimate the efficiency of time-varying congestion pricing when applied to a bi-modal two-tandem bottleneck model. Two possible factors contribute to this result: 1) the mismatch in the capacities of the curbside bottlenecks and the highway bottleneck could reduce the overall effectiveness of congestion pricing; 2) the mode switch of commuters between RV and PV due to congestion pricing may influence overall effectiveness. Additionally, Fig. \ref{fig:DeltaSC_DeltaC} shows that $\frac{\Delta SC}{SC_e}$ decreases as $s_\mathcal{C}^\mathcal{P}$ increases, suggesting that an increase in the capacity of the main road in the urban area may weaken the effectiveness of congestion pricing. This may be attributed to the fact that an increase in $s_\mathcal{C}^\mathcal{P}$ leads to an overall reduction in social cost (as illustrated in Fig. \ref{fig:SC_C_Nr}), and ultimately weakens the effectiveness of congestion pricing.

Furthermore, Fig. \ref{fig:DeltaSC_DeltaC} shows that when $s_\mathcal{C}^\mathcal{R}+s_\mathcal{C}^\mathcal{P}<s_\mathcal{H}$, implementing optimal congestion pricing leads to a reduction in individual travel cost, i.e., $\Delta C<0$. Conversely, when  $s_\mathcal{C}^\mathcal{R}$ and  $s_\mathcal{C}^\mathcal{P}$ are large enough, implementing the optimal congestion pricing can lead to an increase in the individual travel cost, i.e., $\Delta C\geq 0$. These findings underscore the significant influence of road capacity attributes on the performance of optimal congestion pricing, especially its effects on the individual travel cost for commuters.

\subsubsection{The effects of $u^\mathcal{P}-c^\mathcal{R}$}
\begin{figure}[!h]
    \centering
    \includegraphics[width=1\linewidth]{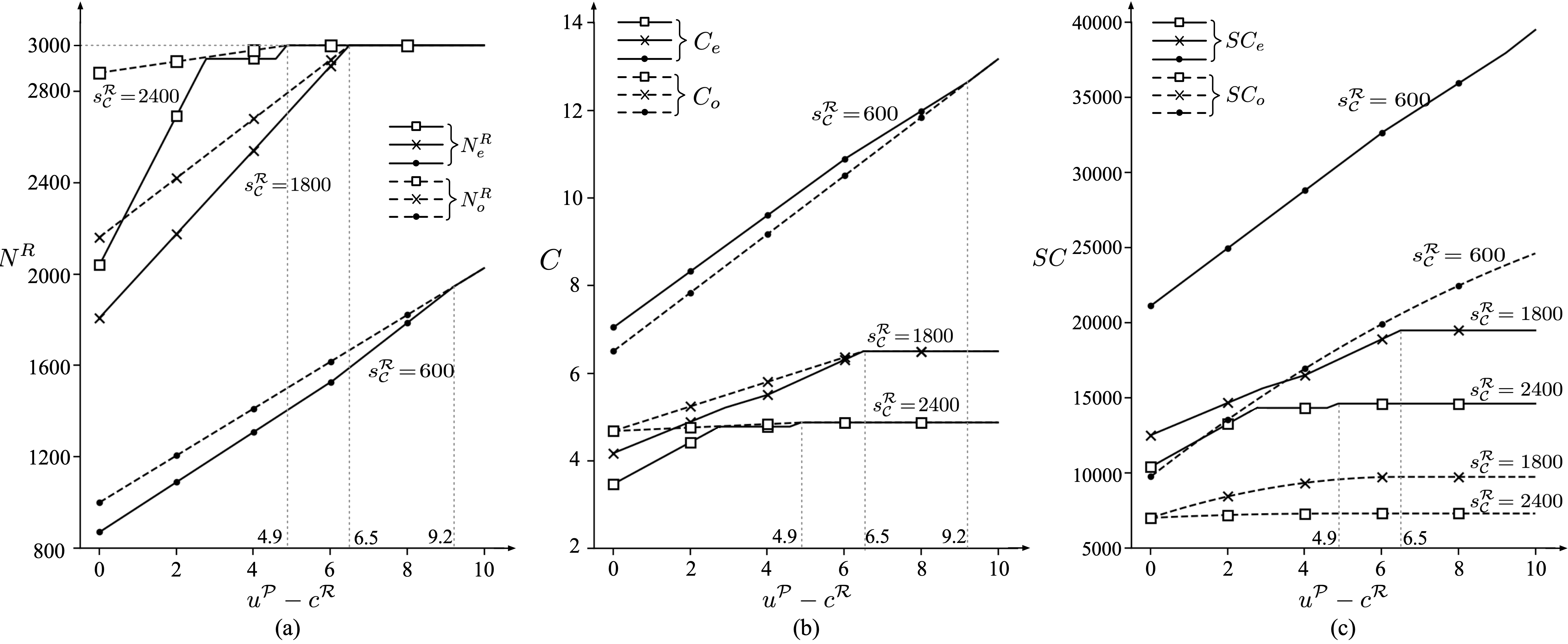}
    \caption{Changes of $N^\mathcal{R}$, $C$, and $SC$ with respect to $(u^\mathcal{P}-c^\mathcal{R})$}
    \label{fig:(uc-cr)_SC&C&Nr}
\end{figure}

Fig. \ref{fig:(uc-cr)_SC&C&Nr} illustrates changes in the modal split $N^\mathcal{R}$, individual travel cost $C$, and social cost $SC$ with respect to the value of $u^\mathcal{P}-c^\mathcal{R}$ before and after implementing the optimal congestion pricing. The solid and dashed lines in Fig. \ref{fig:(uc-cr)_SC&C&Nr} represent the no-toll equilibrium and the social optimum, respectively. The subscripts “o” and “e” of parameters represent the social optimum and the no-toll equilibrium, respectively. The parameter $s_\mathcal{C}^\mathcal{P}$ is set to 1,200 veh/h. 

In general, regardless of whether the congestion pricing is implemented or not, as the value of $u^\mathcal{P}-c^\mathcal{R}$ increases, $N^\mathcal{R}$ increases until it reaches the maximum $N$, both $C$ and $SC$ initially increase and then reach the fixed values when $N^\mathcal{R}=N$. This implies that some current policies to encourage RV use by reducing fixed cost $c^\mathcal{R}$, such as ride-hailing coupons, could indeed increase the modal split of RVs. However, it could also result in an increased individual travel cost and an increased total social cost. Two possible factors contribute to this: 1) the increasing number of RV commuters intensifies the competition for limited curb space, resulting in substantial queuing delays for RV commuters; 2) the more serious congestion spillover effects lead to additional delays for PV commuters. 

Fig. \ref{fig:(uc-cr)_SC&C&Nr}(a) and (b) indicate that for the cases with $s_\mathcal{C}^\mathcal{R} =1800$ veh/h and $s_\mathcal{C}^\mathcal{R} = 2400$ veh/h (as indicated by lines with cross markers and square markers, respectively), the implementation of optimal congestion pricing leads to an increase in $N^\mathcal{R}$ and $C$, while this influence diminishes as the value of $u^\mathcal{P}-c^\mathcal{R}$ increases. However, for the case with $s_\mathcal{C}^\mathcal{R} = 600$ veh/h, the implementation of optimal congestion pricing leads to an increase in $N^\mathcal{R}$ but a decrease in $C$, and this influence diminishes as $u^\mathcal{P}-c^\mathcal{R}$ increases. These findings suggest that when $s_\mathcal{C}^\mathcal{R}$ is small, congestion pricing could reduce the individual travel cost of commuters. Conversely, when $s_\mathcal{C}^\mathcal{R}$ is large enough, the implementation of congestion pricing could lead to an increase in the individual travel cost of commuters. However, this effect gradually diminishes as the value of $u^\mathcal{P} - c^\mathcal{R}$ increases.

Fig. \ref{fig:(uc-cr)_SC&C&Nr}(c) indicates that as the value of $u^\mathcal{P}-c^\mathcal{R}$ increases, the difference between $SC_e$ and $SC_o$ increases and then reaches a fixed value when $N^\mathcal{R}=N$. It suggests that the optimal congestion pricing scheme can reduce more social costs with increasing $u^\mathcal{P}-c^\mathcal{R}$. This may be attributed to the fact that as the value of $u^\mathcal{P}-c^\mathcal{R}$ increases, more commuters switch to using the RV, leading to an overall increase in social cost. Therefore, the optimal congestion pricing scheme, which can completely eliminate congestion spillover effects, would perform better.

\begin{figure}[!h]
    \centering
    \includegraphics[width=0.5\linewidth]{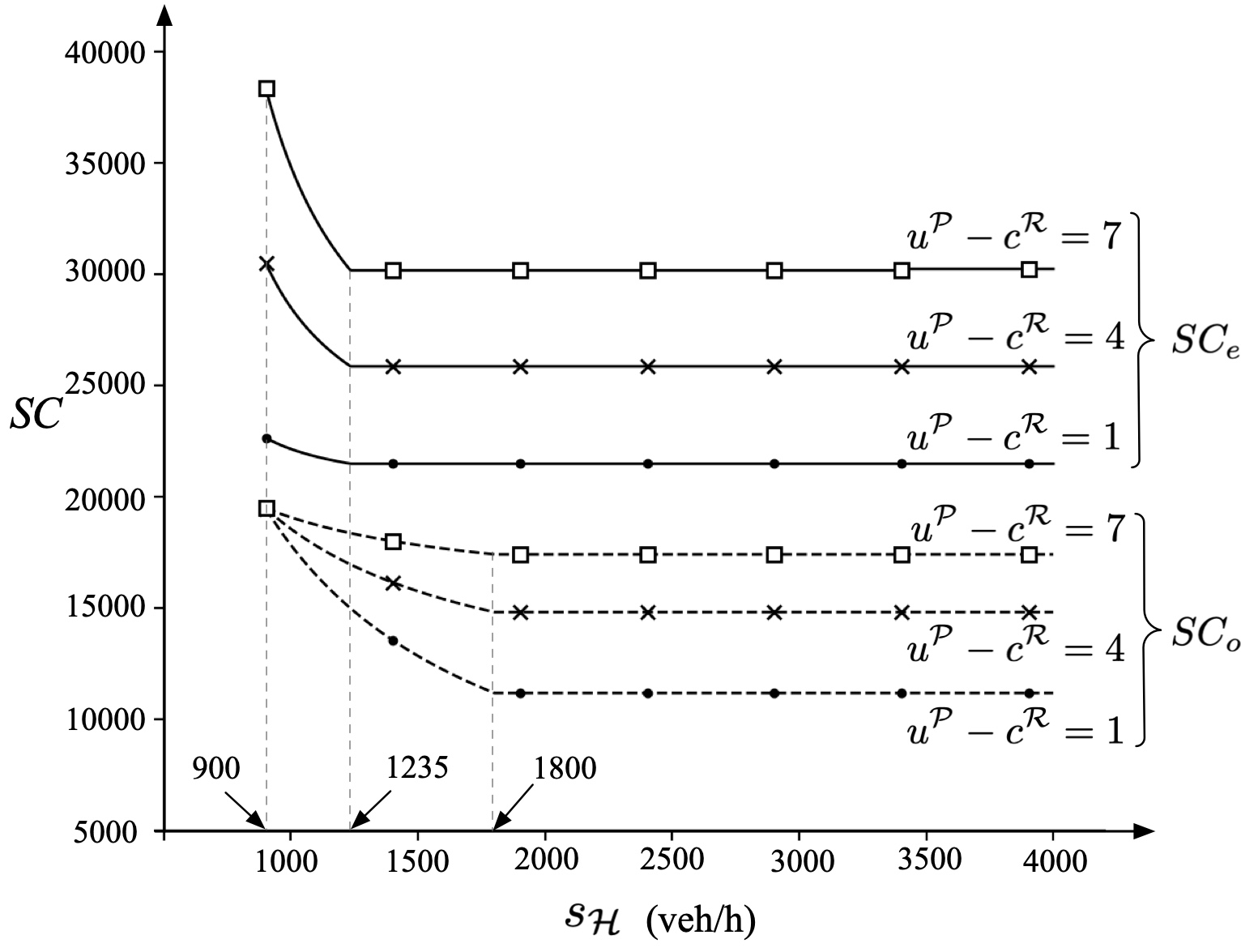}
    \caption{Changes of $SC$ with respect to $s_\mathcal{H}$}
    \label{fig:u2_SC}
\end{figure}

\subsubsection{The effects of the highway bottleneck capacity $s_\mathcal{H}$}
Fig. \ref{fig:u2_SC} shows changes of the social cost SC with respect to the highway bottleneck capacity $s_\mathcal{H}$ under two different cases: the no-toll equilibrium $SC_e$, and the social optimum $SC_o$ where the optimal congestion pricing is implemented. Parameters $s_\mathcal{C}^\mathcal{R}$ and $s_\mathcal{C}^\mathcal{P}$ are set to 900 veh/h and 900 veh/h, respectively.

It can be seen that, as $s_\mathcal{H}$ increases, both $SC_e$ and $SC_o$ initially decrease and then reach the minimum. It is intuitive that the expansion of the highway bottleneck could help reduce the overall social cost. However, as $s_\mathcal{H}>1,235$ veh/h, $SC_e$ no longer decreases with increasing $s_\mathcal{H}$. This suggests that a larger capacity for the highway bottleneck is not always beneficial, and the optimal capacity for the highway bottleneck in the no-toll case is $s_\mathcal{H}=1,235$ veh/h. By contrast, in the social optimum case, the optimal capacity for the highway bottleneck could be $s_\mathcal{H}=1,800$ veh/h, which is equal to the total capacity of the curbside bottlenecks, i.e., $s_\mathcal{H} = s_\mathcal{C}^\mathcal{R}+s_\mathcal{C}^\mathcal{P}$. It suggests that the implementation of the optimal congestion pricing scheme results in a higher optimal capacity for the highway bottleneck. The above findings indicate that implementing optimal congestion pricing combined with optimizing the highway bottleneck capacity could reduce the total social cost effectively.

Furthermore, in combination with the findings in Section \ref{numerical_system_performance} that the marginal benefit of curbside bottleneck expansion decreases as $s_\mathcal{C}^\mathcal{R}$ and $s_\mathcal{C}^\mathcal{P}$ increase (see Fig. \ref{fig:SC_C_Nr}). We have reason to believe that the marginal benefits from expanding the bottleneck, either the highway or the curbside bottlenecks, in the two-tandem bottlenecks corridor may decrease as the capacities of the bottlenecks increase. Furthermore, if the traffic demand is elastic, expanding the bottleneck capacity could bring more induced demand, further reducing the potential benefits \citep{duranton2011Am.Econ.Rev.}.

\subsection{A case study in Hong Kong}
We now apply the proposed model to a real-world case from Hong Kong. First, we offer an overview of the study area and the input data sources used in our analysis. Subsequently, we examine commuting patterns during the morning peak within the selected study area. Finally, we propose strategies to address traffic congestion, including congestion pricing and bottleneck expansion.

\subsubsection{Study area and input data}
\begin{figure}[!ht]
    \centering
    \includegraphics[width=0.5\linewidth]{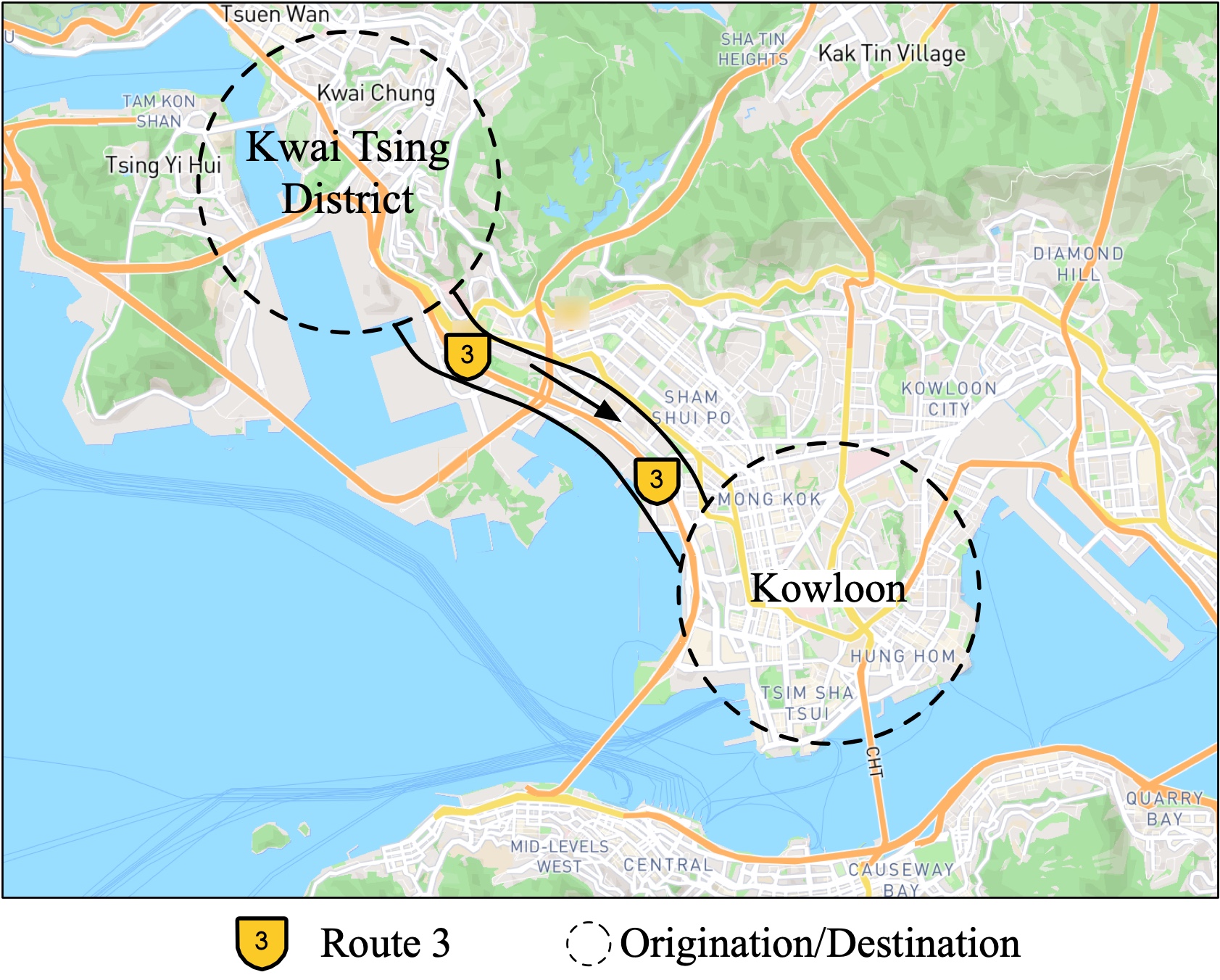}
    \caption{The study area in Hong Kong}
    \label{fig:HK_map}
\end{figure}

We focus on the morning peak commuting flow from Kwai Tsing District to Kowloon in a typical workday in Hong Kong, and we consider a futuristic scenario in which travelers drive or use ride-hailing services to commute when the Mobility-as-a-service (MaaS) becomes available. Specifically, this is a mixed scenario with both solo driving and ride-hailing. The highway named Route 3 links the two districts, as shown in Fig. \ref{fig:HK_map}. Route 3 is a series of highways that have a dual carriageway with three lanes in each direction. The distance between the two districts via Route 3 is approximately 9 km, i.e., $L=9$ km. According to the government's Third Comprehensive Transport Study \citep{thirdTransportStudy}, the practical capacity of Route 3 for one direction is 5,700 veh/h. Consequently, the highway bottleneck capacity can be set to 5,700 veh/h, i.e.,  $s_\mathcal{H}=5,700$ veh/h. A rough estimation of the practical capacity of the main road in Kowloon is 2,100 veh/h. Consequently, the capacity of the curbside bottleneck serving PVs, $s_\mathcal{C}^\mathcal{P}$, can be set to 2,100 veh/h. 

Commuters using ride-hailing services need to drop off in curb spaces near their workplace. We set the area within 500 meters of the central point in Kowloon as the drop-off area for RV commuters. The on-street parking spaces within this area are available for RVs. According to data sourced from Parkopedia \citep{Parkopedia}, the number of available on-street parking spaces is about 90 within this area during the morning peak of a typical workday. A single drop-off process takes approximately 3 minutes, implying that the capacity of these on-street parking spaces is around 1,800 veh/h. Therefore, the capacity of the curbside bottleneck serving RVs, $s_\mathcal{C}^\mathcal{R}$, can be set to 1,800 veh/h. The parameters $\delta^\mathcal{R}$ and $\delta^\mathcal{P}$, which measure the spillover effects of the RV queue on PVs and the PV queue on RVs, are both set to 0.1 in this case.

According to data sourced from \cite{TomTom}, the traffic demand from Kwai Tsing District to Kowloon during the morning peak in a typical workday (from 8:00 am to 9:00 am) is 7,158, i.e,  $N=7,158$ veh. The office hours begin at 9:00 am, i.e., $t^*=9:00$. Following \cite{zhou2022TRC}, the value of travel time $\alpha$ is set to $120$ HKD/h, and the value of schedule delay time $\beta$ is set to $100$ HKD/h based on the assumption $\beta<\alpha$. Other parameters related to the payment for ride-hailing services and the fixed daily cost of car ownership are established following the realistic situations in Hong Kong \citep{taxifeeinHK}. Table \ref{table:case_study_data} provides a summary of these input data.

\begin{table}[!h]
    \caption{Parameters used in the real-world case}\label{table:case_study_data}
    \begin{center}
        \renewcommand{\arraystretch}{1.1}
        \begin{tabular}{l l l}\hline
        \textbf{Parameter} & \textbf{Value}  & \textbf{Description}\\\hline
        $\alpha$ (HKD/h) & 120 & Value of travel time\\
        $\beta$ (HKD/h) & 100 & Value of schedule delay time\\
        $\lambda$ (HKD/km) & 9.5 & Unit price for the trip distance by using RV\\
        $\pi$ (HKD/min) & 1.9 & Unit price for the trip time by using RV\\
        $u^\mathcal{R}$ (HKD) & 27 & Flag-down fee for RV\\
        $u^\mathcal{P}$ (HKD/d) & 200 & Daily fixed cost for car ownership\\
        $L$ (km) & 9 & Total travel distance\\
        $N$ (veh)& 7,158 & Total demand\\
        $s_\mathcal{H}$ (veh/h)& 5,700 & The highway bottleneck capacity\\
        $s_\mathcal{C}^\mathcal{P}$ (veh/h)& 2,100 & The curbside bottleneck capacity for PVs\\
        $s_\mathcal{C}^\mathcal{R}$ (veh/h)& 1,800 & The curbside bottleneck capacity for RVs\\\hline
      \end{tabular}
    \end{center}
\end{table}

\subsubsection{Equilibrium commuting patterns}
\begin{figure}[!h]
    \centering
    \includegraphics[width=0.6\linewidth]{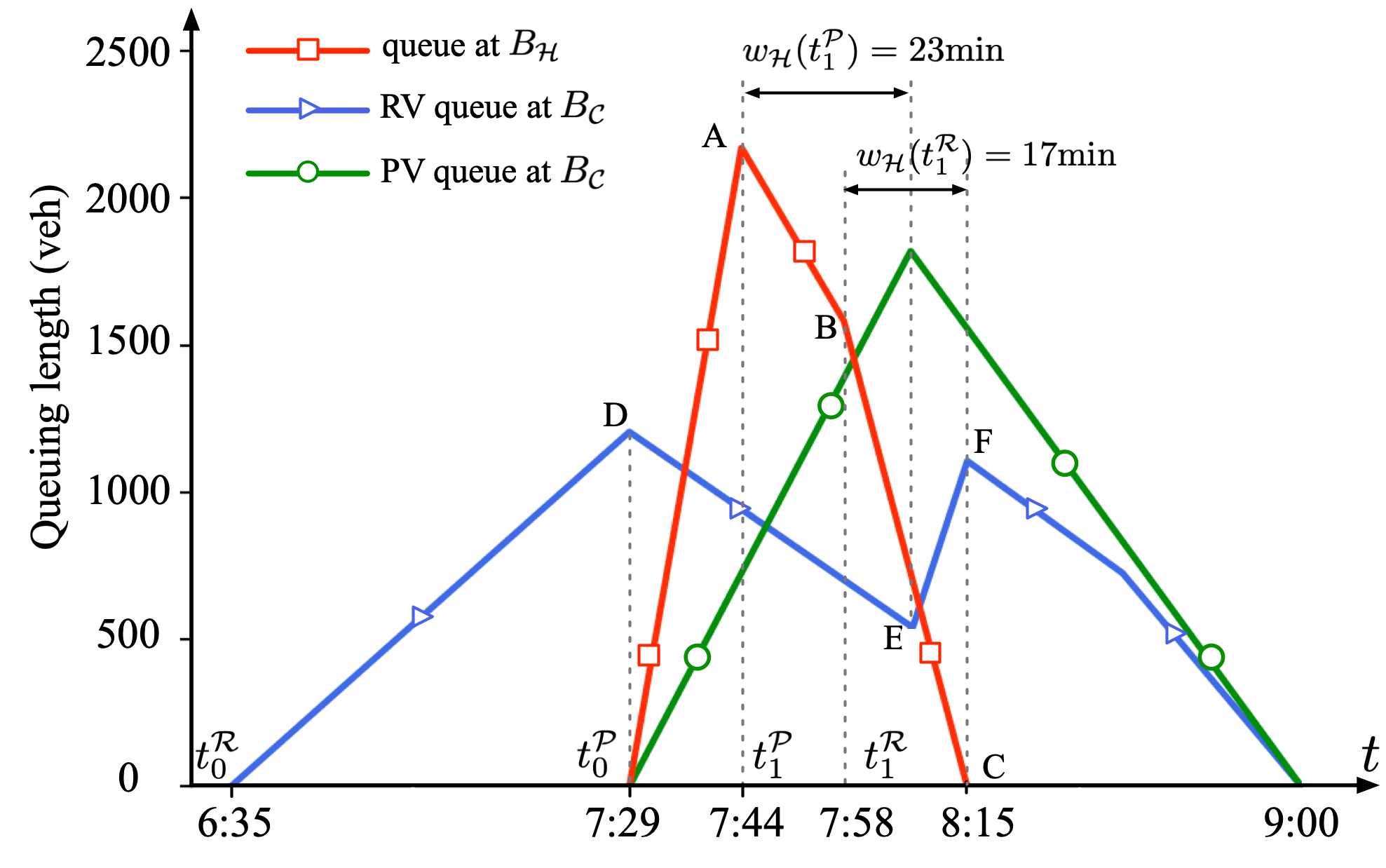}
    \caption{Queue length at bottlenecks}
    \label{fig:HK_queuing}
\end{figure}

Based on the conditions for different equilibrium scenarios presented in Tables \ref{table: conditions_for_8_scenarios}, our analysis indicates that the morning peak commute from Kwai Tsing District to Kowloon in Hong Kong falls into Scenario 5. For illustration purposes, Fig. \ref{fig:HK_queuing} plots the changes in queue length at bottlenecks. Specifically, both the RV mode and the PV mode are utilized, and there is a specific time interval when both modes co-depart from the residential area. Queuing congestion occurs initially at the curbside bottleneck in the urban area and then spreads to the highway bottleneck. 

Fig. \ref{fig:HK_queuing} indicates that the curbside bottleneck $B_\mathcal{C}$ is congested throughout the peak hours [6:35 am, 9:00 am] (as indicated by the blue line with triangle markers and the green line with circle markers), whereas queuing congestion at the highway bottleneck $B_\mathcal{H}$ (as indicated by the red line with square markers) commences when the first PV arrives there at 7:29 am and dissipates once the last RV leaves the highway bottleneck at 8:15 am.

It is worth noting that the decreasing phase of the queue length at the highway bottleneck shows two distinct intervals, as indicated by lines AB and BC. More specifically, in the first interval (i.e., line AB), RV commuters are still departing from home and joining the highway queue, while their departure rate is lower than the highway bottleneck capacity, thus the queue starts to decrease. In the subsequent interval (i.e., line BC), there are no more commuters departing from home, and the queue decreases at a rate of the highway bottleneck capacity. Consequently, the queue length decreases more rapidly than in the first interval. Additionally, within the interval [7:29 am, 8:15 am], RV's queue length at the curbside bottleneck (as indicated by lines DE and EF) experiences a decrease followed by an increase. This is because during the co-departure stage (i.e., line DE), the departure rate of RVs is initially less than the associated service rate of the curbside bottleneck, i.e, $\dot{A}_\mathcal{C}^\mathcal{R}(t)<s_\mathcal{C}^\mathcal{R}$, leading to a decrease in RV's queue length. In the subsequent stage (i.e., line EF), only RVs depart, and their departure rate exceeds the service rate, which leads to an increase in the queue length until all RV commuters leave the residential area.

\subsubsection{Optimal congestion pricing}
The optimal congestion pricing scheme, combining parking pricing and curbside pricing, can be implemented through on-street parking meters and mobile Apps. In fact, parking meters have been widely implemented in Hong Kong. With the Transport Department updating the on-street meters in 2022, there are over 10,700 parking meters installed throughout Hong Kong \citep{HKeMeterinHK}. The deployment of parking meters supporting multiple payment methods and the introduction of mobile Apps like HKeMeter are crucial for implementing dynamic parking pricing and curbside pricing. In this section, we conduct a comparative analysis of the equilibrium results before and after the implementation of the optimal congestion pricing scheme. These insights can help the government evaluate the performance of the congestion pricing scheme and further optimize strategies to alleviate traffic congestion problems.

In this case, the optimal congestion pricing scheme and the equilibrium results are defined by Eqs. \eqref{eq:ueresults_toll_ur1+uc1<u2} and \eqref{eq:toll_level_ur1+uc1<u2} since the condition $s_\mathcal{C}^\mathcal{P}+s_\mathcal{C}^\mathcal{R}<s_\mathcal{H}$ holds. We consider an initial fee of zero, i.e., $f_0=0$. Table \ref{table:case_study_toll} summarizes the equilibrium results before and after the implementation of the optimal congestion pricing scheme. It shows that implementing congestion pricing leads to a delay in the peak hours. The earliest departure times for both RV and PV users are delayed by 6 minutes. The latest departure times for both RV and PV users are delayed by around 1 hour (note that the free-flow travel time is not considered here). Additionally, the number of RV users increases, and the individual travel cost decreases by approximately 9.15 HKD. In addition, compared to the no-toll situation, the implementation of optimal congestion pricing reduces the total social cost by nearly 30\%. These findings highlight that the introduction of optimal congestion pricing in the study area can effectively reduce social cost without bringing any burden to commuters.

\begin{table}[h]
    \caption{Equilibrium results before and after implementing the optimal congestion pricing scheme}\label{table:case_study_toll}
    \begin{center}
        \renewcommand{\arraystretch}{1.1}
        \begin{tabular}{l l l}\hline
        \textbf{Equilibrium results} & \textbf{Before}  & \textbf{After}\\\hline
        RV's congestion fee $f^\mathcal{R}(t)$ (HKD) & & [0,\ 231.67]  \\
        PV's congestion fee $f^\mathcal{P}(t)$ (HKD) & & [0,\ 141.67]  \\
        RV's departure interval & [6:35, 7:58]& [6:41, 9:00] \\
        PV's departure interval & [7:29, 7:44] & [7:35, 9:00] \\
        Number of RV users $N^\mathcal{R}$          & 4027    & 4173 \\
        Individual tralvel cost $C$ (HKD)           & 351.27  & 342.12 \\
        Total social cost $SC$ (HKD)                & 2,515,133.79& 1,752,941.70\\\hline
      \end{tabular}
    \end{center}
\end{table}

\subsubsection{Capacity expansion}\label{subsection_num_Capacity}
\begin{figure}[!h]
    \centering
    \includegraphics[width=0.9\linewidth]{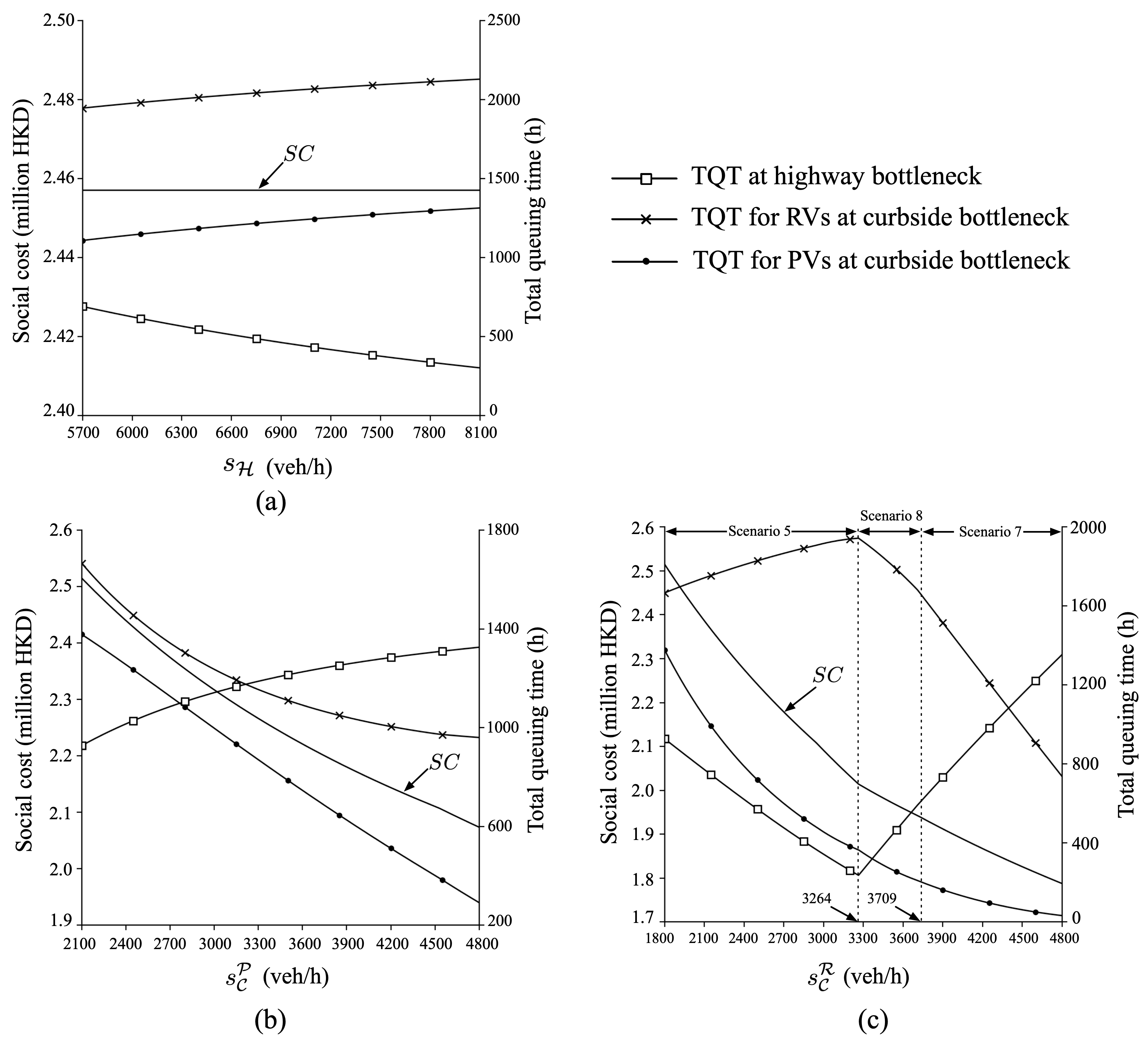}
    \caption{Changes of $SC$ and total queuing time at different bottlenecks with respect to $s_\mathcal{H}$, $s_\mathcal{C}^\mathcal{P}$, or $s_\mathcal{C}^\mathcal{R}$}
    \label{fig:HK_CapacityExpansion}
\end{figure}

Intuitively, a mismatch in capacities between the highway and the curbside bottlenecks could exacerbate urban traffic congestion and weaken the overall effectiveness of congestion pricing. Therefore, it is essential to analyze the influence of each segment's capacity within this two-tandem bottlenecks corridor, as it helps identify the real ``bottleneck'' constraining traffic efficiency in the study area. In this section, we explore the effects of expanding bottleneck capacity on alleviating traffic congestion measured by the social cost $SC$ and the total queuing time, denoted by $TQT$, at different bottlenecks. The $TQT$ is the sum of all commuters' queuing times at the highway bottleneck or at the curbside bottlenecks. Fig. \ref{fig:HK_CapacityExpansion} shows how $SC$ and $TQT$ at different bottlenecks change with increasing $s_\mathcal{H}$, $s_\mathcal{C}^\mathcal{P}$, and $s_\mathcal{C}^\mathcal{R}$, respectively. In particular, with an increase in $s_\mathcal{H}$ or $s_\mathcal{C}^\mathcal{P}$, the equilibrium scenario remains Scenario 5. However, as $s_\mathcal{C}^\mathcal{R}$ increases, the equilibrium scenario transitions from Scenario 5 to Scenario 8 and then to Scenario 7, as shown in Fig. \ref{fig:HK_CapacityExpansion}(c).

Fig. \ref{fig:HK_CapacityExpansion}(a) indicates that, as $s_\mathcal{H}$ increases, $TQT$ at the highway bottleneck decreases, while $TQT$ for RVs and $TQT$ for PVs at the curbside bottlenecks increases. However, the value of $SC$ remains unchanged with increasing $s_\mathcal{H}$. This suggests that expanding the capacity of the highway bottleneck has no significant impact on social cost in this case. Instead, it shifts queuing delay from the highway bottleneck to the curbside bottlenecks. Conversely, as indicated in Figs. \ref{fig:HK_CapacityExpansion}(b) and (c), an increase in $s_\mathcal{C}^\mathcal{P}$ or $s_\mathcal{C}^\mathcal{R}$ results in a significant decrease in social cost $SC$. These findings suggest that the critical “bottleneck” may lie in the curb space or the main road in urban areas. Expanding the curb space or the main road, rather than the highway bottleneck, could effectively reduce the overall social cost. These findings align with Fig. \ref{fig:u2_SC}, which indicates that a larger capacity of the highway bottleneck is not always beneficial.

Fig. \ref{fig:HK_CapacityExpansion}(b) indicates that, as $s_\mathcal{C}^\mathcal{P}$ increases, $TQT$ for PVs and $TQT$ for RVs at the curbside bottlenecks decrease, while $TQT$ at the highway bottleneck increases. This suggests that expanding the capacity of the main road could eliminate queuing congestion at both the curb space and the main road in urban areas, albeit at the cost of exacerbating queuing congestion at the highway bottleneck.

Fig. \ref{fig:HK_CapacityExpansion}(c) illustrates that, as $s_\mathcal{C}^\mathcal{R}$ increases, $TQT$ for PVs at the curbside bottleneck always decreases. However, $TQT$ for RVs at the curbside bottleneck initially increases and then decreases, while $TQT$ at the highway bottleneck first decreases and then increases. Furthermore, the equilibrium scenario moves from Scenario 5 to Scenario 8 and then to Scenario 7 with increasing $s_\mathcal{C}^\mathcal{R}$. These findings suggest that although expanding the capacity of the curb space (i.e., $s_\mathcal{C}^\mathcal{R}$) can reduce the total social cost, it could also alter the equilibrium scenario and result in a longer total queuing time for RV commuters at the curb space. Only when $s_\mathcal{C}^\mathcal{R}> 3,264\ \text{veh/h}$, expanding the capacity of the curb space could reduce the total queuing time for RV commuters at the curb space. However, it also results in an increase in the total queuing time at the highway bottleneck. 

The following analysis can provide insight into the reasons behind the findings in Fig. \ref{fig:HK_CapacityExpansion}(c): when $s_\mathcal{C}^\mathcal{R}\leq 3,264\ \text{veh/h}$, the equilibrium scenario is Scenario 5. As $s_\mathcal{C}^\mathcal{R}$ increases, many PV users switch to using RV, resulting in an increase in $TQT$ for RVs and a decrease in $TQT$ for PVs at the curbside bottlenecks. In the initial stage of morning peak in Scenario 5, only RV commuters leave home and no queue forms at the highway bottleneck. Consequently, as the number of RV users who depart from home in the initial stage increases, $TQT$ at the highway bottleneck decreases. However, when $s_\mathcal{C}^\mathcal{R}>3,264\ \text{veh/h}$, the equilibrium scenario transitions to Scenarios 8 and 7. Queues form at the highway bottleneck from the initial stage when only RV departing from home. A larger $s_\mathcal{C}^\mathcal{R}$ leads to an increased RV users, which results in an increase in $TQT$ at the highway bottleneck during the initial stage and, consequently, throughout the entire morning peak. However, the arrival rate of vehicles at the curbside bottlenecks remains constant and equal to the highway bottleneck capacity $s_\mathcal{H}$. This constant inflow into the curbside bottlenecks, combined with a higher outflow from the curbside bottlenecks due to a larger $s_\mathcal{C}^\mathcal{R}$, results in a reduction in $TQT$ for both PVs and RVs at the curbside bottlenecks.

\newpage
\section{Extension: relaxing the on-late-arrival assumption}\label{section_LateArrival}

In this section, we will discuss the user equilibrium when relaxing the assumption that late arrival is not allowed. Let $\gamma$ denote the value of late arrival time. The generalized travel cost for RV commuters who leave home at time $t$, can be expressed as,
\begin{equation}\label{eq:LateArrival_Cost_RV}
\begin{aligned}
    C_L^\mathcal{R}(t) &= (\alpha+\pi)\big(w_\mathcal{H}(t)+w_\mathcal{C}^\mathcal{R}(t+w_\mathcal{H}(t))\big)+c^\mathcal{R}\\
    &+\max\big\{\beta\big(t^* - t - w_\mathcal{H}(t) - w_\mathcal{C}^\mathcal{R}(t+w_\mathcal{H}(t))\big),\gamma\big(t+ w_\mathcal{H}(t)+w_\mathcal{C}^\mathcal{R}(t+w_\mathcal{H}(t))-t^*\big)\big\},
\end{aligned}
\end{equation} where the subscript ‘L’ of $C_L^\mathcal{R}(t)$ denotes the case considering late arrival. The third term on the RHS of Eq. \eqref{eq:LateArrival_Cost_RV} is the schedule delay cost which depends on whether the commuter arrive early or late at the workplace.

The generalized travel cost for PV commuters who leave home at time $t$ can be expressed as,
\begin{equation}\label{eq:LateArrival_Cost_PV}
\begin{aligned}
    C_L^\mathcal{P}(t) &= \alpha\big(w_\mathcal{H}(t)+w_\mathcal{C}^\mathcal{P}(t+w_\mathcal{H}(t))\big)+u^\mathcal{P}\\
    &+\max\big\{\beta\big(t^* - t - w_\mathcal{H}(t) - w_\mathcal{C}^\mathcal{P}(t+w_\mathcal{H}(t))\big),\gamma\big(t+ w_\mathcal{H}(t)+w_\mathcal{C}^\mathcal{P}(t+w_\mathcal{H}(t))-t^*\big)\big\},
\end{aligned}
\end{equation} where the queuing delay at bottlenecks, i.e., $w_\mathcal{H}(t)$, $w_\mathcal{C}^\mathcal{P}(t+w_\mathcal{H}(t))$, and $w_\mathcal{C}^\mathcal{R}(t+w_\mathcal{H}(t))$ in Eqs. \eqref{eq:LateArrival_Cost_RV} and \eqref{eq:LateArrival_Cost_PV}, are defined in Eqs. \eqref{eq:w2}, \eqref{eq:wcp} and \eqref{eq:wcr_bidirectional}, respectively\footnote{Here we consider the case with bidirectional congestion spillover. Consequently, the queuing delay for PVs at the curbside bottleneck is determined by Eq. \eqref{eq:wcr_bidirectional}.}.

\subsection{Departure intervals of RV and PV commuters}\label{subsec_LateArrival_departure_intervals}
The earliest and latest commuters face no queuing delay at the highway and curbside bottlenecks, as they can reduce their generalized travel costs by departing earlier or later to avoid queuing. Therefore, it is intuitive to deduce that the earliest commuter will arrive at the workplace before $t^*$, while the latest commuter will arrive after $t^*$. Otherwise, their generalized travel costs cannot be equal.

As $u^\mathcal{P}>c^\mathcal{R}$, the earliest and latest commuters are RV commuters, as using the RV mode incurs lower travel costs when there is no queuing delay. However, as queues form at the bottlenecks, commuters will be willing to switch to the PV mode. Therefore, if both modes are used during the morning peak, the departure interval of RV commuters is inclusive of the departure interval of PV commuters, meaning $t_0^\mathcal{R}<t_0^\mathcal{P}<t_1^\mathcal{P}<t_1^\mathcal{R}$. This implies that the departure interval of PV commuters is also the co-departure interval of RV and PV commuters.

\begin{figure}[!h]
    \centering
    \includegraphics[width=\linewidth]{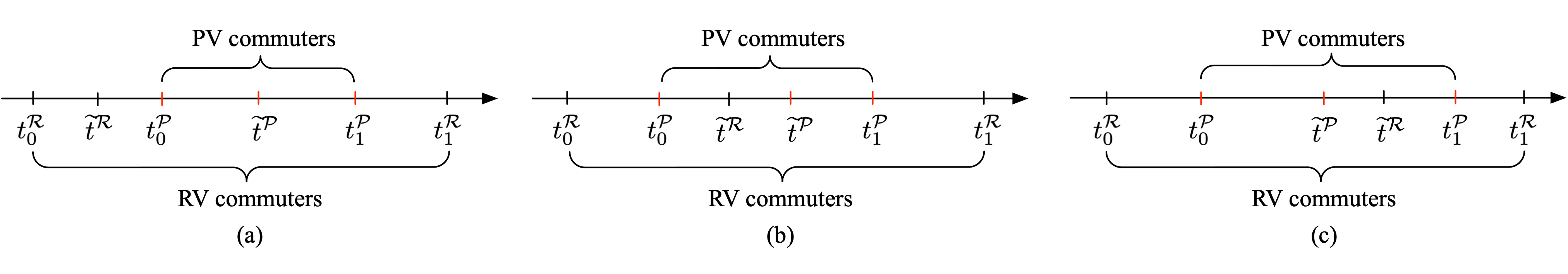}
    \caption{Relationships between departure intervals of RV and PV commuters}
    \label{fig:LateArrival_3cases}
\end{figure}

Let $\tilde{t}^\mathcal{R}$ and $\tilde{t}^\mathcal{P}$ represent the departure times of RV and PV commuters who can arrive at the workplace at the preferred time $t^*$. The relationships among $\tilde{t}^\mathcal{R}$, $\tilde{t}^\mathcal{P}$, $t_0^\mathcal{P}$, and $t_1^\mathcal{P}$ determine whether RV and PV commuters arrive early or late at the workplace during the co-departure stage $(t_0^\mathcal{P},t_1^\mathcal{P})$. For illustration purposes, Fig. \ref{fig:LateArrival_3cases} presents all three cases, illustrating the possible relationships between the departure intervals of RV and PV commuters\footnote{The case of $t_1^\mathcal{P}<\tilde{t}^\mathcal{R}$ does not exist. In this case, during the co-departure stage of $(t_0^\mathcal{P},t_1^\mathcal{P})$, RV commuters always arrive early, while PV commuters may arrive early or late at the workplace. Given that the latest PV commuter arrives late, we have $t_1^\mathcal{P}+w_\mathcal{H}(t_1^\mathcal{P})>t^*$. For the RV commuter who departs at the same time $t_1^\mathcal{P}$, the condition $t_1^\mathcal{P}+w_\mathcal{H}(t_1^\mathcal{P})+w_\mathcal{C}^\mathcal{R}(t_1^\mathcal{P}+w_\mathcal{H}(t_1^\mathcal{P}))>t^*$ always holds. This implies the RV commuter who departs at $t_1^\mathcal{P}$ must arrive late, which conflicts with the condition $t_1^\mathcal{P} < \tilde{t}^\mathcal{R}$. Therefore, this case does not exist.}.

Here we briefly introduce the three cases:

(1) If $\tilde{t}^\mathcal{R}<t_0^\mathcal{P}$, as shown in Fig. \ref{fig:LateArrival_3cases}(a), then PV commuters may arrive early or late, while RV commuters always arrive late at the workplace during the co-departure stage $(t_0^\mathcal{P},t_1^\mathcal{P})$.

(2) If $t_0^\mathcal{P}<\tilde{t}^\mathcal{R}<\tilde{t}^\mathcal{P}$, as shown in Fig. \ref{fig:LateArrival_3cases}(b), then PV and RV commuters may arrive early or late at the workplace during the co-departure stage $(t_0^\mathcal{P},t_1^\mathcal{P})$. The co-departure stage can be divided into three intervals: both RV and PV commuters arrive early within $(t_0^\mathcal{P},\tilde{t}^\mathcal{R})$, PV commuters arrive early while RV commuters arrive late within $(\tilde{t}^\mathcal{R},\tilde{t}^\mathcal{P})$, and both RV and PV commuters arrive late within $(\tilde{t}^\mathcal{P},t_1^\mathcal{P})$.

(3) If $\tilde{t}^\mathcal{P}<\tilde{t}^\mathcal{R}<t_1^\mathcal{P}$,  as shown in Fig. \ref{fig:LateArrival_3cases}(c), then PV and RV commuters may arrive early or late at the workplace during the co-departure stage $(t_0^\mathcal{P},t_1^\mathcal{P})$. The co-departure stage can be divided into three intervals: both RV and PV commuters arrive early within $(t_0^\mathcal{P},\tilde{t}^\mathcal{P})$, RV commuters arrive early while PV commuters arrive late within $(\tilde{t}^\mathcal{P},\tilde{t}^\mathcal{R})$, and both RV and PV commuters arrive late within $(\tilde{t}^\mathcal{R},t_1^\mathcal{P})$.

\subsection{User equilibrium}
The departure situations of RV and PV commuters can be categorized into three cases as shown in Fig. \ref{fig:LateArrival_3cases}. In addition, the queuing situations at bottlenecks can be divided into three cases: queuing only at the curbside bottlenecks $B_\mathcal{C}$, queuing initially at $B_\mathcal{C}$ and then expanding to the highway bottleneck $B_\mathcal{H}$, and queuing at both $B_\mathcal{C}$ and $B_\mathcal{H}$ from the start of the morning peak. Consequently, based on whether both modes of transportation are used, the departure situations of RV and PV commuters, and the queuing situations at bottlenecks, all possible user equilibrium situations when considering late arrival can be divided into 11 scenarios, as shown in Table \ref{table:LateArrival_11_scenarios}. To distinguish from the scenarios without considering late arrival discussed earlier, we denote the scenarios with late arrival as Scenarios L1 - L11.

\begin{table}[h]
    \caption{Overview of 11 scenarios when considering late arrival}\label{table:LateArrival_11_scenarios}
    \begin{center}
        \begin{tabular}{l l l l l}\hline
        \diagbox{\textbf{Queuing situation}}{\textbf{Utilization}}
                            & \multirow{2}{*}{\textbf{only RV}}  &  \multicolumn{3}{c}{\textbf{RV+PV}}\\
                            &          &	$\tilde{t}^\mathcal{R}<t_0^\mathcal{P}$	&	$t_0^\mathcal{P}<\tilde{t}^\mathcal{R}<\tilde{t}^\mathcal{P}$ & $\tilde{t}^\mathcal{P}<\tilde{t}^\mathcal{R}<t_1^\mathcal{P}$ \\\hline
            Only $B_\mathcal{C}$ queuing	            & Scenario L1      &	Scenario L2	& Scenario L3 & Scenario L4\\
            $B_\mathcal{C}$ queuing, expanding to  $B_\mathcal{C}+B_\mathcal{H}$ queuing & $-$             &Scenario L5   &	Scenario L6 & Scenario L7	\\
            $B_\mathcal{C}+B_\mathcal{H}$	queuing	        & Scenario L8      &	Scenario L9	&	Scenario L10 & Scenario L11\\\hline
      \end{tabular}
    \end{center}
\end{table}

Due to space limitations, we do not discuss the occurrence conditions and equilibrium results for all 11 scenarios in detail. Instead, we focus on Scenarios L7 as an example to illustrate the equilibrium patterns for morning peak commuting when the no-late-arrival assumption is relaxed. In Scenario L7, both travel modes are used. The departure intervals of RV and PV commuters is illustrated by Fig. \ref{fig:LateArrival_3cases}(c), implying that RV and PV commuters who depart during the co-departure stage could arrive early or late at the workplace. The queuing congestion initially occurs at the curbside bottleneck, and then expands to the highway bottleneck once the first PV commuter departs from home, implying that $w_\mathcal{C}^\mathcal{R}(t)>0$ and $w_\mathcal{H}(t)=0$ holds for $t\in(t_0^\mathcal{R},t_0^\mathcal{P})$.  

In the initial stage of $(t_0^\mathcal{R},t_0^\mathcal{P})$, only RV commuters depart from home and arrive early at the workplace. In addition, there is no queue at the highway bottleneck, i.e., $w_\mathcal{H}(t)=0$. Substitute Eq. \eqref{eq:wcr_bidirectional} into Eq. \eqref{eq:LateArrival_Cost_RV}, and let $\frac{\partial C_L^\mathcal{R}(t)}{\partial t} = 0$, we can derive the equilibrium departure rate of RV commuters during this period as,
\begin{equation}\label{eq:LateArrival_SL7_initial_rate}
    \dot{A}_\mathcal{C}^\mathcal{R}(t) = \dot{A}_\mathcal{H}(t)=\frac{(\alpha+\pi)s_\mathcal{C}^\mathcal{R}}{\alpha+\pi-\beta},\quad\  t\in(t_0^\mathcal{R},t_0^\mathcal{P}).
\end{equation}

In the final stage of $(t_1^\mathcal{P},t_1^\mathcal{R})$, only RV commuters depart from home and arrive late at the workplace. Substitute Eqs. \eqref{eq:w2} and \eqref{eq:wcr_bidirectional} into Eq. \eqref{eq:LateArrival_Cost_RV}, and let $\frac{\partial C_L^\mathcal{R}(t)}{\partial t} = 0$, we thus have\footnote{If $w_\mathcal{H}(t)=0$ holds for $(t_1^\mathcal{P},t_1^\mathcal{R})$, then $\dot{A}_\mathcal{H} = \dot{A}_\mathcal{C}^\mathcal{R} = \frac{(\alpha+\pi)s_\mathcal{C}^\mathcal{R}}{\alpha+\pi+\gamma}$.},
\begin{equation}
    \begin{cases}
        \dot{A}_\mathcal{H}(t) = \frac{(\alpha+\pi)s_\mathcal{C}^\mathcal{R}}{\alpha+\pi+\gamma}, &t\in(t_1^\mathcal{P},t_1^\mathcal{R}),\\[1.5ex]
        \dot{A}_\mathcal{C}^\mathcal{R}(t) = s_\mathcal{H}, &t\in(t_1^\mathcal{P}+w_\mathcal{H}(t_1^\mathcal{P}),t_1^\mathcal{R}+w_\mathcal{H}(t_1^\mathcal{R})),
    \end{cases}
\end{equation}

In the co-departure stage, both RV and PV commuters arrive early within $(t_0^\mathcal{P},\tilde{t}^\mathcal{P})$; RV commuters arrive early but PV commuters arrive late within $(\tilde{t}^\mathcal{P},\tilde{t}^\mathcal{R})$; both RV and PV commuters arrive late within $(\tilde{t}^\mathcal{R},t_1^\mathcal{P})$. In addition, queues form at the highway bottleneck once the first PV commuter depart from home at time $t_0^\mathcal{P}$. Substitute Eqs. \eqref{eq:w2}, \eqref{eq:wcp} and \eqref{eq:wcr_bidirectional} into Eqs. \eqref{eq:LateArrival_Cost_RV} and \eqref{eq:LateArrival_Cost_PV}, let $\frac{\partial C_L^\mathcal{R}(t)}{\partial t}=\frac{\partial C_L^\mathcal{P}(t)}{\partial t} = 0$, we can derive commuters' total departure rate from home (i.e., total arrival rate at the highway bottleneck) as,
\begin{equation}\label{eq:LateArrival_co_DepartureRate_S7}
    \dot{A}_\mathcal{H}(t) = \begin{cases}
        \frac{\alpha(1-\delta^\mathcal{P})s_\mathcal{C}^\mathcal{P}}{(\alpha-\beta)(1-\delta^\mathcal{R}\delta^\mathcal{P})}+\frac{(\alpha+\pi)(1-\delta^\mathcal{R})s_\mathcal{C}^\mathcal{R}}{(\alpha+\pi-\beta)(1-\delta^\mathcal{R}\delta^\mathcal{P})}, 
        &t\in(t_0^\mathcal{P},\tilde{t}^\mathcal{P}),\\[1.5ex]
        
        \frac{\alpha(1-\delta^\mathcal{P})s_\mathcal{C}^\mathcal{P}}{(\alpha+\gamma)(1-\delta^\mathcal{R}\delta^\mathcal{P})}+\frac{(\alpha+\pi)(1-\delta^\mathcal{R})s_\mathcal{C}^\mathcal{R}}{(\alpha+\pi-\beta)(1-\delta^\mathcal{R}\delta^\mathcal{P})}, &t\in(\tilde{t}^\mathcal{P},\tilde{t}^\mathcal{R}),\\[1.5ex]

        \frac{\alpha(1-\delta^\mathcal{P})s_\mathcal{C}^\mathcal{P}}{(\alpha+\gamma)(1-\delta^\mathcal{R}\delta^\mathcal{P})}+\frac{(\alpha+\pi)(1-\delta^\mathcal{R})s_\mathcal{C}^\mathcal{R}}{(\alpha+\pi+\gamma)(1-\delta^\mathcal{R}\delta^\mathcal{P})}, 
        &t\in(\tilde{t}^\mathcal{R},t_1^\mathcal{P}).
    \end{cases}\\[1.5ex]
\end{equation} and the arrival rates at curbside bottlenecks as,
\begin{equation}\label{eq:LateArrival_co_DepartureRate_S7_2}
\begin{aligned}
&\dot{A}_\mathcal{C}^\mathcal{R}(t) = \begin{cases}
        \frac{(\alpha+\pi)(\alpha-\beta)s_\mathcal{C}^\mathcal{R}-\alpha(\alpha+\pi-\beta)\delta^\mathcal{P}s_\mathcal{C}^\mathcal{P}}{(\alpha+\pi)(\alpha-\beta)(1-\delta^\mathcal{R})s_\mathcal{C}^\mathcal{R}+\alpha(\alpha+\pi-\beta)(1-\delta^\mathcal{P})s_\mathcal{C}^\mathcal{P}}s_\mathcal{H}, &t\in(t_0^\mathcal{P},\tilde{t}^\mathcal{P}+w_\mathcal{H}(\tilde{t}^\mathcal{P})),\\[1.5ex]
        
        \frac{(\alpha+\pi)(\alpha+\gamma)s_\mathcal{C}^\mathcal{R}-\alpha(\alpha+\pi-\beta)\delta^\mathcal{P}s_\mathcal{C}^\mathcal{P}}{(\alpha+\pi)(\alpha+\gamma)(1-\delta^\mathcal{R})s_\mathcal{C}^\mathcal{R}+\alpha(\alpha+\pi-\beta)(1-\delta^\mathcal{P})s_\mathcal{C}^\mathcal{P}}s_\mathcal{H}, &t\in(\tilde{t}^\mathcal{P}+w_\mathcal{H}(\tilde{t}^\mathcal{P}),\tilde{t}^\mathcal{R}+w_\mathcal{H}(\tilde{t}^\mathcal{R})),\\[1.5ex]

        \frac{(\alpha+\pi)(\alpha+\gamma)s_\mathcal{C}^\mathcal{R}-\alpha(\alpha+\pi+\gamma)\delta^\mathcal{P}s_\mathcal{C}^\mathcal{P}}{(\alpha+\pi)(\alpha+\gamma)(1-\delta^\mathcal{R})s_\mathcal{C}^\mathcal{R}+\alpha(\alpha+\pi+\gamma)(1-\delta^\mathcal{P})s_\mathcal{C}^\mathcal{P}}s_\mathcal{H}, &t\in(\tilde{t}^\mathcal{R}+w_\mathcal{H}(\tilde{t}^\mathcal{R}),t_1^\mathcal{P}+w_\mathcal{H}(t_1^\mathcal{P})).
    \end{cases}\\[1.5ex]
    &\dot{A}_\mathcal{C}^\mathcal{P}(t) = s_\mathcal{H}-\dot{A}_\mathcal{C}^\mathcal{R}(t),\ t\in(t_0^\mathcal{P},t_1^\mathcal{P}+w_\mathcal{H}(t_1^\mathcal{P})).
\end{aligned}
\end{equation}

It can be observed that the total departure rate of RV and PV commuters within $(t_0^\mathcal{P},\tilde{t}^\mathcal{P})$ in Scenario L7 (see Eq. \eqref{eq:LateArrival_co_DepartureRate_S7}) is identical to that within the co-departure stage in Scenario 5 (see $\dot{A}_\mathcal{H}(t)$ in Eq. \eqref{eq:bieffect_DepartRate_S5}). In addition, those commuters' arrival rates at curbside bottlenecks within $t\in(t_0^\mathcal{P},\tilde{t}^\mathcal{P}+w_\mathcal{H}(\tilde{t}^\mathcal{P}))$ are also same as those in Scenario 5 (see $\dot{A}_\mathcal{C}^\mathcal{R}(t)$ and $\dot{A}_\mathcal{C}^\mathcal{P}(t)$ in Eq. \eqref{eq:bieffect_DepartRate_S5}). It implies that the period when both RV and PV commuters arrive early in Scenario L7 corresponds to the co-departure stage in Scenario 5. Moreover, the equilibrium departure and arrival rates of commuters remain unchanged when considering late arrival. 

In addition, for periods when RV commuters arrive early while PV commuters arrive late, or when both RV and PV commuters arrive late in Scenario L7, the equilibrium departure and arrival rates are symmetrical to those of the period when both RV and PV commuters arrive early. Only the value of early arrival time $\beta$ is replaced by the negative value of late arrival time $-\gamma$, while the functional forms of the departure and arrival rates remain unchanged. This indicates that, except differences in the value of times, commuter behavior when RV or PV arrive late are essentially similar to that when both arrive early. Commuter behavior in the late arrival intervals can be regarded as analogous and symmetrical to that in the early arrival intervals.

At user equilibrium, all RV commuters have an equal generalized travel cost, i.e., $C_L^\mathcal{R}(t_0^\mathcal{R})=C_L^\mathcal{R}(\tilde{t}^\mathcal{R})=C_L^\mathcal{R}(t_1^\mathcal{R})$. Therefore, the follow equations ensure the intra-equilibrium for RV commuters:
\begin{equation}\label{eq:linear_equations_SL7_RV}
    \begin{cases}
        \beta(t^*-t_0^\mathcal{R})+c^\mathcal{R}=\gamma(t_1^\mathcal{R}-t^*)+c^\mathcal{R},\\
        \beta(t^*-t_0^\mathcal{R})+c^\mathcal{R}=(\alpha+\pi)(t^*-\tilde{t}^\mathcal{R})+c^\mathcal{R}.
    \end{cases}
\end{equation}

Also, all PV commuters have an equal generalized travel cost, i.e., $C_L^\mathcal{P}(t_0^\mathcal{P})=C_L^\mathcal{P}(\tilde{t}^\mathcal{P})=C_L^\mathcal{P}(t_1^\mathcal{P})$. Therefore, the follow equations ensure the intra-equilibrium for PV commuters:
\begin{equation}\label{eq:linear_equations_SL7_PV}
    \begin{cases}
        \beta(t^*-t_0^\mathcal{P})+u^\mathcal{P}=\alpha w_\mathcal{H}(t_1^\mathcal{P})+\gamma(t_1^\mathcal{P}+w_\mathcal{H}(t_1^\mathcal{P})-t^*)+u^\mathcal{P},\\
        \beta(t^*-t_0^\mathcal{P})+u^\mathcal{P}=\alpha(t^*-\tilde{t}^\mathcal{P})+u^\mathcal{P}.
    \end{cases}
\end{equation}

At user equilibrium, all RV commuters arrive at the workplace within the interval of $(t_0^\mathcal{R},t_1^\mathcal{R})$. All PV commuters arrive within the interval of $(t_0^\mathcal{P},t_1^\mathcal{P})$. In addition, RV and PV commuters have an equal generalized travel cost. Therefore, the follow equations ensure the inter-equilibrium for RV and PV commuters:
\begin{equation}\label{eq:LateArrival_SL7_interEquilibrium}
    \begin{cases}
        \int_{{t}_0^\mathcal{R}}^{t_1^\mathcal{R}} \dot{A}_\mathcal{H}(t) dt=N,\\[1.5ex]

        \int_{{t}_0^\mathcal{R}}^{t_1^\mathcal{R}} \dot{A}_\mathcal{C}^\mathcal{R}(t) dt=N^\mathcal{R},\\[1.5ex]

        \beta(t^*-t_0^\mathcal{R})+c^\mathcal{R}=\beta(t^*-t_0^\mathcal{P})+u^\mathcal{P}.
    \end{cases}
\end{equation}

By combining the intra-equilibrium conditions of Eqs. \eqref{eq:linear_equations_SL7_RV}-\eqref{eq:linear_equations_SL7_PV}, and the inter-equilibrium conditions of Eqs. \eqref{eq:LateArrival_SL7_interEquilibrium}, one can derive the equilibrium results of Scenario L7, including the critical departure times (i.e, $t_0^\mathcal{R},\tilde{t}^\mathcal{R},t_1^\mathcal{R},t_0^\mathcal{P},\tilde{t}^\mathcal{P},t_1^\mathcal{P}$ ) and the mode split $N^\mathcal{R}$. 

\begin{figure}[!h]
    \centering
    \includegraphics[width=0.9\linewidth]{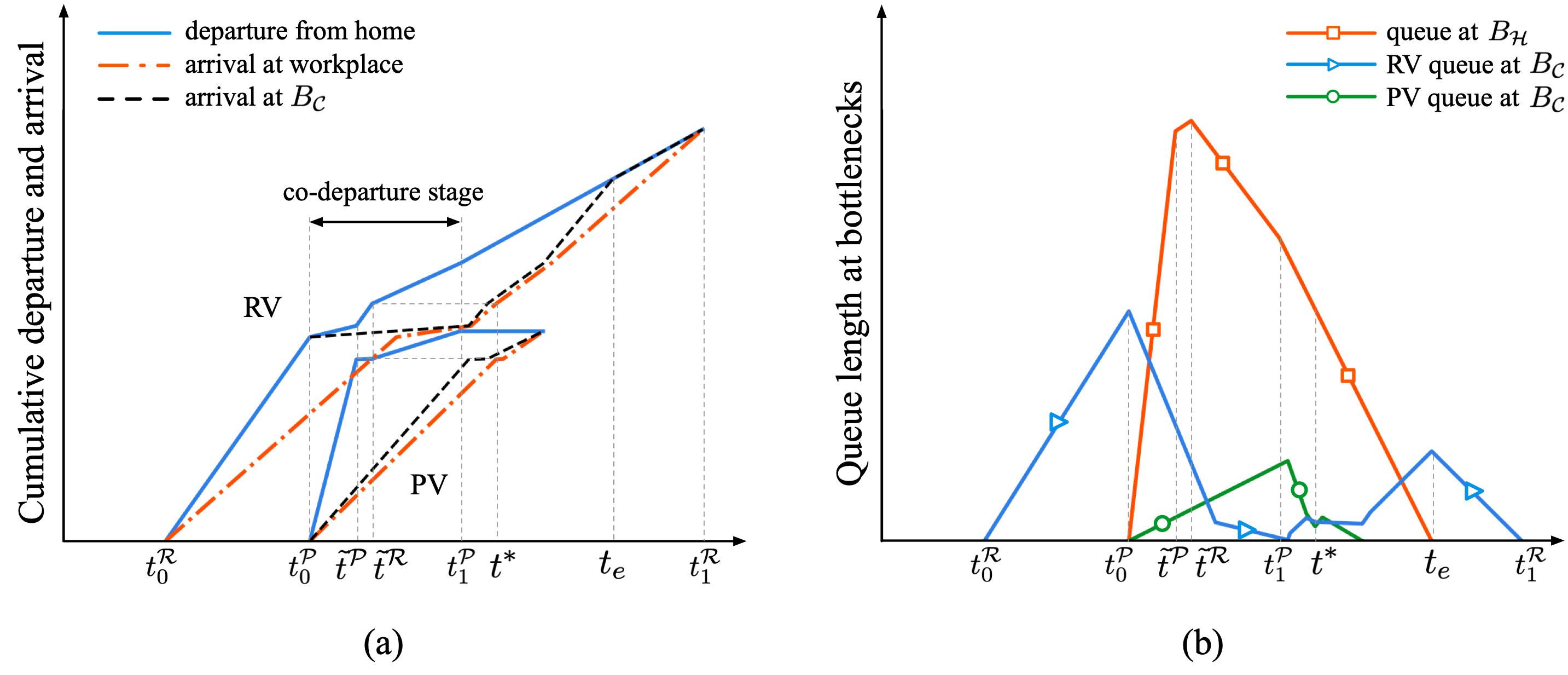}
    \caption{Equilibrium patterns for Scenario L7: (a) cumulative departure and arrival of commuters; (b) queue length at bottlenecks over time}
    \label{fig:LateArrival_SL7}
\end{figure}

Based on an numerical example with $N=6500\ \text{veh}, s_\mathcal{H}=5500\ \text{veh/h}, s_\mathcal{C}^\mathcal{R}=4000\ \text{veh/h}, s_\mathcal{C}^\mathcal{P}=4500\ \text{veh/h}$\footnote{Other parameters for this numerical example are as follows: $u^\mathcal{P}-c^\mathcal{R} = \$ 1, 
\alpha=\$4/h, \beta=\$3/h, \pi=\$4/h, \gamma=\$4.8/h$, $\delta^\mathcal{R}=\delta^\mathcal{P}=0.3$.}, Fig. \ref{fig:LateArrival_SL7} illustrates the cumulative departure and arrival of commuters during the morning peak, and changes in queue length at bottlenecks over time for Scenario L7. Fig. \ref{fig:LateArrival_SL7}(a) shows that, in the initial stage of $(t_0^\mathcal{R},t_0^\mathcal{P})$, only RV commuters depart from home, and RV queue forms exclusive at the curbside bottleneck as shown in Fig. \ref{fig:LateArrival_SL7}(b). In the final stage of $(t_1^\mathcal{P},t_1^\mathcal{R})$, only RV commuters depart, and those who depart after $t_e$ face no queuing delay at the highway bottleneck because the highway bottleneck queue dissipates at time $t_e$ as shown in Fig. \ref{fig:LateArrival_SL7}(b). Additionally, in Fig. \ref{fig:LateArrival_SL7}(a), the co-departure stage of $(t_0^\mathcal{P}, t_1^\mathcal{P})$ shows three different departure rates for RV and PV commuters (i.e., slopes of the cumulative departure curves). These correspond to different departure intervals, which ensure that both RV and PV commuters arrive early, RV commuters arrive early while PV commuters arrive late, or both arrive late at the workplace. Aa a result, the arrival patterns of those commuters at the bottlenecks (black dashed lines) and at the workplace (red dash-dotted lines), as well as the queue length at bottlenecks over time shown in Fig. \ref{fig:LateArrival_SL7}(b) are more intricate than the case without considering late arrival.

\subsection{Optimal congestion pricing}
It should be noted that the design of optimal congestion pricing when considering late arrival follows the same logic as without considering late arrival. It requires discussions for the two cases of $s_\mathcal{C}^\mathcal{R} + s_\mathcal{C}^\mathcal{P} \leq s_\mathcal{H}$  and $s_\mathcal{C}^\mathcal{R} + s_\mathcal{C}^\mathcal{P} > s_\mathcal{H}$, respectively. For Scenario L7 with $s_\mathcal{H}=5500\ \text{veh/h}, s_\mathcal{C}^\mathcal{R}=4000\ \text{veh/h}$, and $s_\mathcal{C}^\mathcal{P}=4500\ \text{veh/h}$ shown in Fig. \ref{fig:LateArrival_SL7}, $s_\mathcal{C}^\mathcal{R} + s_\mathcal{C}^\mathcal{P} > s_\mathcal{H}$ holds. Therefore, both the individual departure rates of RV and PV commuters, as well as their total departure rate during the co-departure stage, should not exceed the service rates of the bottlenecks. Therefore, in the initial and last stages when only RV commuters depart from home, the departure rate of RV commuters should be equal to $s_\mathcal{C}^\mathcal{R}$ under the optimal congestion pricing scheme. In the co-departure stage, the total departure rate of RV and PV commuters should be equal to $s_\mathcal{H}$. Let $\theta$ be the ratio of PV commuters' departure rate to the curbside bottleneck capacity during the co-departure stage, i.e., $\theta = \dot {A}_\mathcal{C}^\mathcal{P}(t)/s_\mathcal{C}^\mathcal{P}$. The departure rates of PV and RV commuters during the co-departure stage can thus be expressed as $\theta s_\mathcal{C}^\mathcal{P}$ and $s_\mathcal{H}-\theta s_\mathcal{C}^\mathcal{P}$, respectively. The parameter $\theta$ falls within the range of $\big[(s_\mathcal{H}-s_\mathcal{C}^\mathcal{R})/s_\mathcal{C}^\mathcal{P},1\big]$, ensuring full utilization of both bottlenecks without queuing delay. 

Under the optimal congestion pricing scheme, the generalized travel costs for RV and PV commuters who arrive at the workplace at time $t$ can be expressed as,
\begin{equation}\label{eq:LateArrival_toll_Cost_define}
\begin{aligned}
    \hat{C}_L^\mathcal{R}(t) &= \begin{cases}
        f_L^\mathcal{R}(t)+c^\mathcal{R}+\beta(t^* - t),\ t\in[\hat{t}_0^\mathcal{R},t^*],\\
        f_L^\mathcal{R}(t)+c^\mathcal{R}+\gamma(t - t^*),\ t\in(t^*,\hat{t}_1^\mathcal{R}].
    \end{cases}\\
    \hat{C}_L^\mathcal{P}(t) &= \begin{cases}
        f_L^\mathcal{P}(t)+u^\mathcal{P}+\beta(t^* - t),\ t\in[\hat{t}_0^\mathcal{P},t^*],\\
        f_L^\mathcal{P}(t)+u^\mathcal{P}+\gamma(t - t^*),\ t\in(t^*,\hat{t}_1^\mathcal{P}].
    \end{cases}\\
\end{aligned}
\end{equation} where $\hat{t}_0^\mathcal{R}$ and $\hat{t}_0^\mathcal{P}$ denote the arrival times of the earliest RV and PV commuters after charging the congestion fee. $\hat{t}_1^\mathcal{R}$ and $\hat{t}_1^\mathcal{P}$ denote the arrival times of the latest RV and PV commuters after charging the congestion fee.

Set $\frac{\partial \hat{C}_L^\mathcal{R}(t)}{\partial t}=\frac{\partial \hat{C}_L^\mathcal{P}(t)}{\partial t}=0$, one can obtain, 
\begin{equation}\label{eq:LateArrival_pricing_rate}
    \frac{\partial f_L^\mathcal{R}(t)}{\partial t}=\frac{\partial f_L^\mathcal{P}(t)}{\partial t}=\begin{cases}
        \beta, &t\leq t^*,\\
        -\gamma, &t> t^*.
    \end{cases}
\end{equation} Eq. \eqref{eq:LateArrival_pricing_rate} shows that the congestion fees $f_L^\mathcal{R}(t)$ and $f_L^\mathcal{P}(t)$ are piecewise linear functions with respective to the arrival time. The congestion fees increase at a rate of $\beta$ during the early arrival period, while decrease at a rate of $\gamma$ during the late arrival period.

Let $f_L^\mathcal{R}(\hat{t}_0^\mathcal{R})$ and $f_L^\mathcal{P}(\hat{t}_0^\mathcal{P})$ be the initial fees for RV and PV commuters, $f_L^\mathcal{R}(\hat{t}_1^\mathcal{R})$ and $f_L^\mathcal{P}(\hat{t}_1^\mathcal{P})$ be the congestion fees paid by the latest RV and PV commuters. Therefore, the optimization for congestion pricing when considering late arrival can be formulated as the following social cost minimization problem: 
\begin{alignat}{2}
    \begin{split}
    \mathop{\min}_{f_L^\mathcal{R}(\hat{t}_0^\mathcal{R}),\ f_L^\mathcal{P}(\hat{t}_0^\mathcal{P}),\ f_L^\mathcal{R}(\hat{t}_1^\mathcal{R}),\ f_L^\mathcal{P}(\hat{t}_1^\mathcal{P})} SC &=\int_{\hat{t}_0^\mathcal{R}}^{\hat{t}_0^\mathcal{P}} \big(\hat{C}^\mathcal{R}(t)-f_L^\mathcal{R}(t)\big) s_\mathcal{C}^\mathcal{R} dt + \int_{\hat{t}_1^\mathcal{P}}^{\hat{t}_1^\mathcal{R}} \big(\hat{C}^\mathcal{R}(t)-f_L^\mathcal{R}(t)\big) s_\mathcal{C}^\mathcal{R} dt\\
    &+\int_{\hat{t}_0^\mathcal{P}}^{\hat{t}_1^\mathcal{P}}\big(\hat{C}^\mathcal{R}(t)-f_L^\mathcal{R}(t)\big) (s_\mathcal{H}-\theta s_{C}^\mathcal{P}) dt+\int_{\hat{t}_0^\mathcal{P}}^{\hat{t}_1^\mathcal{P}} \big(\hat{C}^\mathcal{P}(t)-f_L^{P}(t)\big) \theta s_\mathcal{C}^\mathcal{P} dt,
    \end{split}\label{eq:LateArrival_minimizationSC_ur1+uc1>u2}\\
\mbox{s.t.}\quad
    &\quad\beta(t^*-\hat{t}_0^\mathcal{R})+f_L^\mathcal{R}(\hat{t}_0^\mathcal{R})+c^\mathcal{R} = \gamma(\hat{t}_1^\mathcal{R}-t^*)+f_L^\mathcal{R}(\hat{t}_1^\mathcal{R})+c^\mathcal{R}\label{eq:LateArrival_constrain2a}, \\
    &\quad\beta(t^*-\hat{t}_0^\mathcal{P})+f_L^\mathcal{P}(\hat{t}_0^\mathcal{P})+u^\mathcal{P} = \gamma(\hat{t}_1^\mathcal{P}-t^*)+f_L^\mathcal{P}(\hat{t}_1^\mathcal{P})+u^\mathcal{P}\label{eq:LateArrival_constrain2b},\\
    &\quad\beta(t^*-\hat{t}_0^\mathcal{R})+f_L^\mathcal{R}(\hat{t}_0^\mathcal{R})+c^\mathcal{R} = 
    \beta(t^*-\hat{t}_0^\mathcal{P})+f_L^\mathcal{P}(\hat{t}_0^\mathcal{P})+u^\mathcal{P}\label{eq:LateArrival_constrain2c},\\
    &\quad(\hat{t}_0^\mathcal{P}-\hat{t}_0^\mathcal{R}+\hat{t}_1^\mathcal{R}-\hat{t}_1^\mathcal{P})s_\mathcal{C}^\mathcal{R} + (\hat{t}_1^\mathcal{P}-\hat{t}_0^\mathcal{P})(s_\mathcal{H}-\theta s_\mathcal{C}^\mathcal{P})=N^\mathcal{R}\label{eq:LateArrival_constrain2d}, \\
    &\quad(\hat{t}_1^\mathcal{P}-\hat{t}_0^\mathcal{P})\theta s_\mathcal{C}^\mathcal{P} = N^\mathcal{P}\label{eq:LateArrival_constrain2e}, \\
    &\quad N^\mathcal{R}+N^\mathcal{P} = N \label{eq:LateArrival_constrain2f}.
\end{alignat} where Eqs. \eqref{eq:LateArrival_constrain2a}-\eqref{eq:LateArrival_constrain2c} indicate that the earliest and latest RV commuters, as well as the earliest and latest PV commuters, have equal generalized travel costs. Eqs. \eqref{eq:LateArrival_constrain2d}-\eqref{eq:LateArrival_constrain2f} are population constraints.

By solving Eqs. \eqref{eq:LateArrival_constrain2a}-\eqref{eq:LateArrival_constrain2f} and substituting the results of $\hat{t}_0^\mathcal{R}$, $\hat{t}_0^\mathcal{P}$, $\hat{t}_1^\mathcal{R}$, and $\hat{t}_1^\mathcal{P}$ into Eq. \eqref{eq:LateArrival_minimizationSC_ur1+uc1>u2}, then setting $\frac{\partial SC}{\partial f_L^\mathcal{R}(\hat{t}_0^\mathcal{R})} = \frac{\partial SC}{\partial f_L^\mathcal{P}(\hat{t}_0^\mathcal{P})} = \frac{\partial SC}{\partial f_L^\mathcal{R}(\hat{t}_1^\mathcal{R})} = \frac{\partial SC}{\partial f_L^\mathcal{P}(\hat{t}_1^\mathcal{P})} = 0$, one can obtain the first-order optimality conditions for this minimization problem\footnote{It can be proven, based on the numerical example in Fig. \ref{fig:LateArrival_SL7}, that the Hessian matrix of the objective function in Eq. \eqref{eq:LateArrival_minimizationSC_ur1+uc1>u2} is semi-positive definite. Thus, the objective function is convex, and the optimal solution to the minimization problem can be determined by solving the first-order optimality conditions.}:
\begin{equation}\label{eq:LateArrival_firstoptimalitycondition}
    \begin{cases}
        f_L^\mathcal{R}(\hat{t}_0^\mathcal{R})=f_L^\mathcal{R}(\hat{t}_1^\mathcal{R}),\\[0.8ex]
        f_L^\mathcal{P}(\hat{t}_0^\mathcal{P}) = f_L^\mathcal{P}(\hat{t}_1^\mathcal{P}),\\[0.8ex]
        f_L^\mathcal{P}(\hat{t}_0^\mathcal{P})-f_L^\mathcal{R}(\hat{t}_0^\mathcal{R}) = \frac{(u^\mathcal{P}-c^\mathcal{R})(s_\mathcal{C}^\mathcal{R}+\theta s_\mathcal{C}^\mathcal{P}-s_\mathcal{H})}{s_\mathcal{H}-s_\mathcal{C}^\mathcal{R}}.
    \end{cases}
\end{equation}

Eq. \eqref{eq:LateArrival_firstoptimalitycondition} indicates that the earliest and latest RV commuters pay an equal congestion fee. Also, the earliest and latest PV commuters pay an equal congestion fee. Define $\Delta f_L = f_L^\mathcal{P}(\hat{t}_0^\mathcal{P})-f_L^\mathcal{R}(\hat{t}_0^\mathcal{R})$ as the difference in the initial fee between PV and RV commuters. Substituting Eq. \eqref{eq:LateArrival_firstoptimalitycondition} into Eq. \eqref{eq:LateArrival_minimizationSC_ur1+uc1>u2}, and taking the first-order derivative of $SC$ with respective to $\Delta f_L$, we have,
\begin{equation}\label{eq:LateArrival_dSC/df}
    \frac{\partial SC}{\partial \Delta f_L} = \frac{N(s_\mathcal{H}-s_\mathcal{C}^\mathcal{R})}{s_\mathcal{H}}-\frac{(\Delta f_L+u^\mathcal{P}-c^\mathcal{R})s_\mathcal{C}^\mathcal{R}(s_\mathcal{H}-s_\mathcal{C}^\mathcal{R})(\beta+\gamma)}{\beta\gamma s_\mathcal{H}}.
\end{equation} 

For the bi-modal utilization scenarios, we have $\Delta f_L \leq \frac{N\beta\gamma}{s_\mathcal{C}^\mathcal{R}(\beta+\gamma)} - (u^\mathcal{P} - c^\mathcal{R})$ when considering late arrival\footnote{If $\Delta f_L +(u^\mathcal{P} - c^\mathcal{R}) > \frac{N\beta\gamma}{s_\mathcal{C}^\mathcal{R}(\beta+\gamma)}$, then the cost difference between the first PV and RV commuters under the optimal congestion pricing, $\Delta f_L +(u^\mathcal{P} - c^\mathcal{R})$, exceeds the generalized travel cost for RV commuters when all commuters choose the RV mode, i.e., $\frac{N\beta\gamma}{s_\mathcal{C}^\mathcal{R}(\beta+\gamma)}$. Therefore, for the bi-modal utilization scenarios, the condition $\Delta f_L \leq \frac{N\beta\gamma}{s_\mathcal{C}^\mathcal{R}(\beta+\gamma)} - (u^\mathcal{P} - c^\mathcal{R})$ holds.}. Based on Eq. \eqref{eq:LateArrival_dSC/df}, it is not hard to prove that $\frac{\partial SC}{\partial \Delta f}\geq 0$ holds for the bi-modal utilization scenarios with late arrival. Additionally, the value of $\Delta f_L$ should fall within the range of $\big[0,\frac{(u^\mathcal{P}-c^\mathcal{R})(s_\mathcal{C}^\mathcal{R}+ s_\mathcal{C}^\mathcal{P}-s_\mathcal{H})}{(s_\mathcal{H}-s_\mathcal{C}^\mathcal{R})}\big]$ to ensure full utilization of both bottlenecks without queuing delay\footnote{As $\theta\in\big[\frac{(s_\mathcal{H}-s_\mathcal{C}^\mathcal{R})}{s_\mathcal{C}^\mathcal{P}},1\big]$, the range of $\Delta f_L$ can be derived based on Eq. \eqref{eq:LateArrival_firstoptimalitycondition}.}. Consequently, the optimal $\Delta f_L^*$ to minimize the social cost takes the value of its left boundary, i.e., $\Delta f_L^* = 0$, implying an identical initial fee for RV and PV commuters. Let $f_0$ represent an arbitrary initial fee, thus we have $f_L^\mathcal{R}(\hat{t}_0^\mathcal{R})=f_L^\mathcal{P}(\hat{t}_0^\mathcal{P})=f_0$. Therefore, the optimal congestion pricing for this case can be expressed as,
\begin{equation}\label{eq:LateArrival_toll_level_ur1+uc1>u2}
    \begin{aligned}
        &f^\mathcal{R}(t)=
            \begin{cases}
            f_0,\ &t\leq\hat{t}_0^\mathcal{R}\ \text{or}\ t>\hat{t}_1^\mathcal{R},\\
            f_0+\beta(t-\hat{t}_0^\mathcal{R}),\ &t\in(\hat{t}_0^\mathcal{R}, t^*],\\
            f_0+\gamma(\hat{t}_1^\mathcal{R}-t),\ &t\in(t^*,\hat{t}_1^\mathcal{R}].\\
            \end{cases}\\
        &f^\mathcal{P}(t)=
            \begin{cases}
            f_0,\ &t\leq\hat{t}_0^\mathcal{P}\ \text{or}\ t>\hat{t}_1^\mathcal{P},\\
            f_0+\beta(t-\hat{t}_0^\mathcal{P}),\ &t\in(\hat{t}_0^\mathcal{P}, t^*],\\
            f_0+\gamma(\hat{t}_1^\mathcal{P}-t),\ &t\in(t^*,\hat{t}_1^\mathcal{P}].\\
            \end{cases}\\
    \end{aligned}
\end{equation}

One can further derive the equilibrium results with optimal congestion pricing by solving Eqs. \eqref{eq:LateArrival_constrain2a}-\eqref{eq:LateArrival_constrain2f} in combination with Eqs. \eqref{eq:LateArrival_toll_level_ur1+uc1>u2}. Based on Eqs. \eqref{eq:LateArrival_toll_level_ur1+uc1>u2}, for the co-departure stage of RV and PV commuters, we have,
\begin{equation}
    \label{eq:LateArrival_toll_level_diff_co_departure}
    f^\mathcal{R}(t)-f^\mathcal{P}(t) = \beta(\hat{t}_0^\mathcal{P}-\hat{t}_0^\mathcal{R})= \gamma(\hat{t}_1^\mathcal{R}-\hat{t}_1^\mathcal{P})=u^\mathcal{P}-c^\mathcal{R},\ \forall t\in[\hat{t}_0^\mathcal{P},\hat{t}_1^\mathcal{P}].
\end{equation} 

\begin{figure}[!h]
    \centering
    \includegraphics[width=0.5\linewidth]{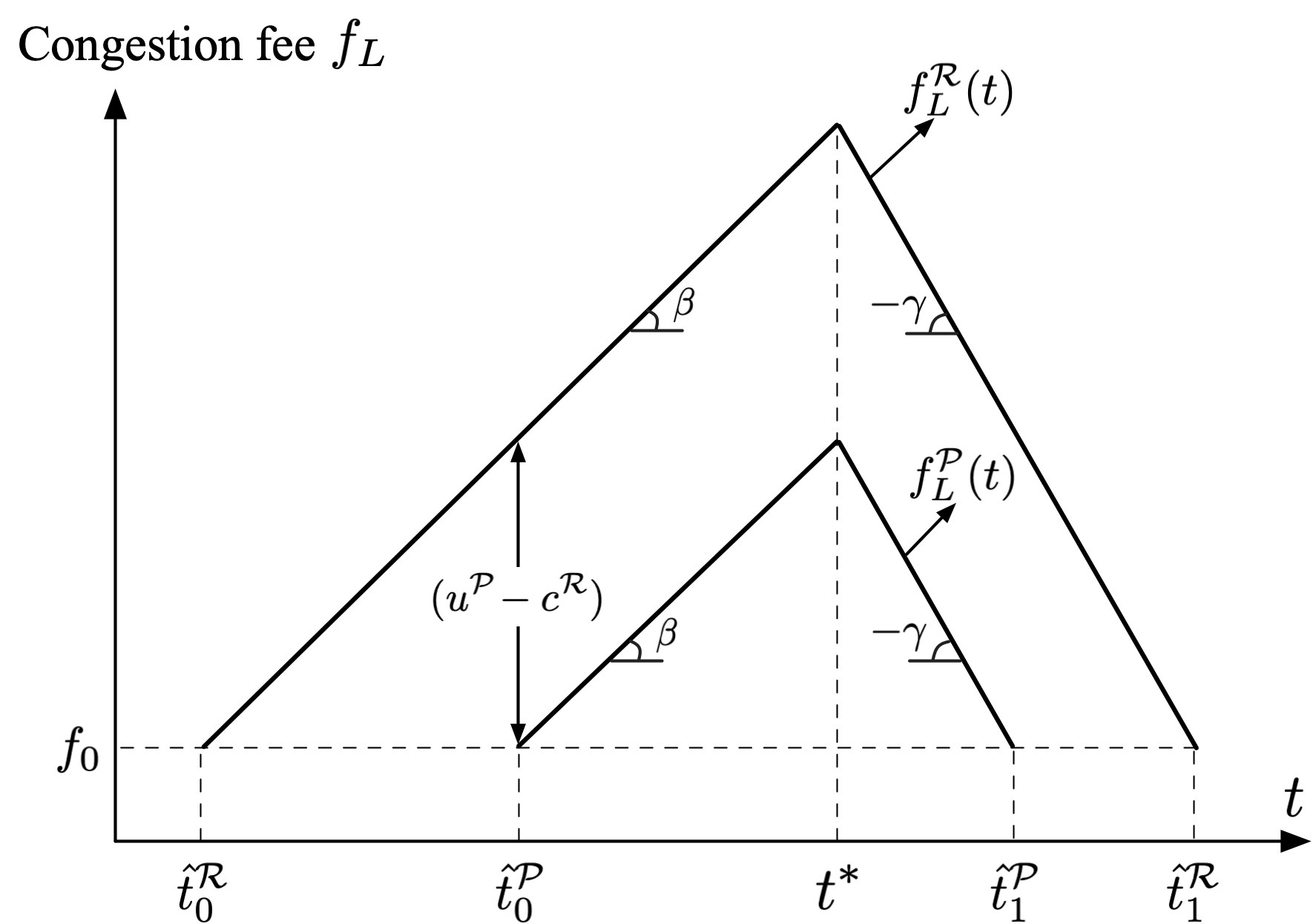}
    \caption{The optimal congestion pricing for RV and PV commuters with late arrival}
    \label{fig:LateArrival_congestion_fee}
\end{figure}

Eq. \eqref{eq:LateArrival_toll_level_diff_co_departure} indicates that even when considering late arrival, the optimal congestion fee difference on RV and PV commuters during the co-departure stage remains $u^\mathcal{P} - c^\mathcal{R}$. Additionally, RV and PV commuters have a uniform initial fee and same marginal charge rates as indicated by Eqs. \eqref{eq:LateArrival_toll_level_ur1+uc1>u2}. This suggests that, even when late arrival is considered, Proposition \ref{Proposition 4} still holds. For illustration purpose, Fig. \ref{fig:LateArrival_congestion_fee} plots the optimal congestion fee with respect to time $t$. It shows that, during the co-departure stage, the RV commuters should be charged with just a markup of $u^\mathcal{P} - c^\mathcal{R}$ over the PV commuters.

\section{Conclusion and further studies}\label{section_conclusion}
This paper investigates the congestion effects of curbside drop-offs by ride-hailing services. 
We propose a novel bi-modal two-tandem bottleneck model that examines congestion patterns on the highway and the curb space during the morning peak. Furthermore, we propose an optimal congestion pricing scheme, which combines curbside pricing at curb spaces and parking pricing at CBD parking lots, to effectively manage both curb space and highway congestion. We also apply the model to a real-world case in Hong Kong and explore the efficacy of bottleneck capacity expansion to alleviate traffic congestion. Some important and insightful findings are obtained. 
\begin{itemize}
    \item First, if the curbside capacity is below the threshold, the bottleneck congestion occurs initially at the curb space, and then extends to the highway bottleneck over time; otherwise, both the curbside and highway bottlenecks are congested from the start of the morning peak (Proposition \ref{Proposition 1}). In addition, commuters initially opt for a lower-cost transportation mode, but as congestion delay increases, they may switch to another mode of transportation. We have identified the specific conditions under which both modes (i.e., RV and PV) are utilized during the morning peak (Proposition \ref{Proposition 2}). Eight equilibrium scenarios are identified and solved with closed-form solutions.
    
    \item Second, we propose an optimal congestion pricing scheme combining curbside pricing with dynamic parking fees to eliminate bottleneck congestion in the curb space and highway bottlenecks simultaneously (Proposition \ref{Proposition 4}).  However, the congestion pricing may not work similarly well to the conventional single bottleneck model \citep{vickrey1969,arnott1990TransportationResearchPartB:Methodological, arnott1990JournalofUrbanEconomicsa}, potentially wasting road capacity when there is a mismatch between the highway and curbside bottlenecks. 
    
    \item Third, the numerical results indicate that expanding the curb space or the main road in urban area, rather than the highway bottleneck, can effectively reduce the total social cost, whereas resulting in a shift of queuing congestion to the highway bottleneck. Furthermore, a real-world case from Hong Kong shows that the morning peak commute from Kwai Tsing District to Kowloon could fall into Scenario 5, where both RV and PV are utilized. Queuing congestion initially forms at the curbside bottleneck before spreading to the highway bottleneck. The optimal congestion pricing can reduce the commuter's generalized travel cost, and lower the total social cost by nearly 30\%. 

    \item Last, relaxing the no-late-arrival assumption introduces more complex and varied equilibrium scenarios. However, commuter behavior in the late arrival intervals, especially the equilibrium departure/arrival rates, remains analogous and symmetrical to that in the early arrival intervals. The optimal congestion pricing scheme remains straightforward, and Proposition \ref{Proposition 4} still holds when considering late arrival.
\end{itemize}

In terms of future studies, we plan to extend this study in the following directions: (1) Implementing dynamic or time-varying tolls faces high implementation costs and may confuse travelers due to their bounded rationality. Therefore, we will focus on exploring more economical and practical step pricing schemes as an alternative \citep{laih1994,li&lam2017TransportationResearchPartB:Methodological}; (2) Relaxing assumptions regarding bottleneck capacities ($s_\mathcal{C}^\mathcal{R}\leq s_\mathcal{H}$ and $s_\mathcal{C}^\mathcal{P}\leq s_\mathcal{H}$ ) allows us to explore a wide range of possibilities and scenarios; (3) Examining the pick-up and drop-off processes of ride-hailing services together when modeling the morning-evening commuting behavior \citep{zhang&yang2005TransportationResearchPartA:PolicyandPractice,zhang&huang2008TransportationResearchPartB:Methodological}; (4) Curb space is not only used by ride-hailing services, but also used as on-street parking spaces, truck loading zones, bus stops, and scooter lanes \citep{mitman2018,jaller2021,lim&masoud2024,liu&qian2024TransportationScience}. Therefore, we plan to incorporate the competition for curb space from various traffic modes into our model and study the optimization of curb space allocation; (5) The application of traffic congestion management schemes based on the bottleneck model to real large-scale networks has long been a key focus for researchers. A DTA-based bottleneck modeling approach \citep{depalma&marchal2002Netw.Spat.Econ.,depalma&kilani2005TransportationResearchPartA:PolicyandPractice,depalma&lindsey2006TransportPolicy,depalma&lindsey2008J.Transp.Econ.Policy} offers a novel perspective for addressing this issue, making it a promising direction for future research.

\section*{Acknowledgments}
Y. Deng and Z. Li were supported by the National Natural Science Foundation of China (72131008), the Fundamental Research Funds for the Central Universities (2021GCRC014), and the Interdisciplinary Research Program of HUST (2023JCYJ021). S. Qian was supported by the U.S. National Science Foundation Grant No. CMMI-1931827. Y. Deng and W. Ma were also supported by grants from the Research Grants Council of the Hong Kong Special Administrative Region, China (Project No. PolyU/25209221 and PolyU/15206322).

\appendix
\section{Derivation of different user equilibrium scenarios}
\setcounter{equation}{0}
\renewcommand\theequation{A.\arabic{equation}}
\label{Appendix A}

Here we will provide a detailed derivation of the user equilibrium results, including the equilibrium departure and arrival rates, as well as the critical time points of the morning peak, for the eight scenarios without imposing any congestion pricing.

\subsubsection*{(1) Scenario 1}
In Scenario 1, where only RV is utilized and the queue forms at the curbside bottleneck exclusively, we can treat it as a traditional bottleneck model with a single traffic mode. The arrival rate of RVs at the CBD is equal to the service rate of the curbside bottleneck, i.e., $s_\mathcal{C}^\mathcal{R}$. Therefore, the user equilibrium results can be obtained by solving the following system of equations:
\begin{equation}
    \label{eq:linearsystem_S1}
    \begin{cases}
        \frac{\partial c^\mathcal{R}(t)}{\partial t}=(\alpha+\pi-\beta)\dot{w}_\mathcal{C}^\mathcal{R}(t)-\beta=0,\\
        \beta(t^*-t_0^\mathcal{R})+c^\mathcal{R}=(\alpha+\pi)(t^*-t_1^\mathcal{R})+c^\mathcal{R},\\
        (t^*-t_0^\mathcal{R})s_\mathcal{C}^\mathcal{R}=N.\\
    \end{cases}
\end{equation}

Therefore, the user equilibrium results for Scenario 1 are as follows:
\begin{equation}
    \begin{cases}
        t_0^\mathcal{R} = t^*-\frac{N}{s_\mathcal{C}^\mathcal{R}},\\[1.2ex]
        
        t_1^\mathcal{R} = t^*-\frac{\beta}{\alpha+\pi}\frac{N}{s_\mathcal{C}^\mathcal{R}},\\[1.2ex]

        \dot{A}_\mathcal{C}^\mathcal{R}(t)=\frac{(\alpha+\pi)s_\mathcal{C}^\mathcal{R}}{\alpha+\pi-\beta}, & t\in[t_0^\mathcal{R},t_1^\mathcal{R}].\\[1.2ex]
    \end{cases}
\end{equation}

In Scenario 1, congestion occurs only at the curbside bottleneck. It implies that $\dot{A}_\mathcal{C}^\mathcal{R}(t)<s_\mathcal{H}$ holds for $t\in[t_0^\mathcal{R},t_1^\mathcal{R}]$, we thus have,
\begin{equation}
\label{eq:AppendixB_conditions_S1_1}
        \frac{s_\mathcal{C}^\mathcal{R}}{s_\mathcal{H} }\leq \frac{\alpha+\pi-\beta}{\alpha+\pi}.
\end{equation}

Since only RV is used during the morning peak, the cost of using PV is always higher than that of using RV, we thus have,
\begin{equation}
\label{eq:AppendixB_conditions_S1_2}
        (u^\mathcal{P}-c^\mathcal{R})\geq\frac{\beta N}{s_\mathcal{C}^\mathcal{R}}.
\end{equation}

Eqs. \eqref{eq:AppendixB_conditions_S1_1} and \eqref{eq:AppendixB_conditions_S1_2} present the occurrence conditions for Scenario 1.

\subsubsection*{(2) Scenario 2 and 4}
In Scenarios 2 and 4, the departure intervals of RV and PV commuters do not overlap and only RV commuters will depart from home within $(t_0^\mathcal{R},t_1^\mathcal{R})$. Consequently, the following system of equations ensures RV commuters' intra-equilibrium within this period:  
\begin{equation}
    \label{eq:linearsystem_S2and4_1}
    \begin{cases}
        \frac{\partial c^\mathcal{R}(t)}{\partial t}=(\alpha+\pi-\beta)\dot{w}_\mathcal{C}^\mathcal{R}(t)-\beta=0,\\
        
        \beta(t^*-t_0^\mathcal{R})+c^\mathcal{R}=(\alpha+\pi)(t^*-t_1^\mathcal{R})+c^\mathcal{R},\\
        
        (t^*-t_0^\mathcal{R})s_\mathcal{C}^\mathcal{R}=N^\mathcal{R}.\\
    \end{cases}
\end{equation}

During the period $(t_0^\mathcal{P},t_1^\mathcal{P})$ when only PV commuters depart from home, we further:
\begin{itemize}
    \item suppose queues form exclusively at the curbside bottleneck, then Scenario 2 occurs. The following equation system ensures the intra-equilibrium for PV commuters:
\begin{equation}
    \label{eq:linearsystem_S2and4_2_onlyB1}
    \begin{cases}
        \frac{\partial C^\mathcal{P}(t)}{\partial t}=(\alpha-\beta)\dot{w}_\mathcal{C}^\mathcal{P}(t)-\beta=0,\\
        \beta(t^*-t_0^\mathcal{P})+u^\mathcal{P}=\alpha(t^*-t_1^\mathcal{P})+u^\mathcal{P},\\
        \int_{{t}_0^\mathcal{P}}^{t_1^\mathcal{P}} \dot{A}_\mathcal{C}^\mathcal{P}(t) dt=N-N^\mathcal{R}.
    \end{cases}
\end{equation}
    \item suppose queues form at both the highway and the curbside bottlenecks (the case where queue forms exclusively at the highway bottleneck is not feasible as $s_\mathcal{H}>s_\mathcal{C}^\mathcal{P}$ holds in this paper), then Scenario 4 occurs. The following equation system ensures the intra-equilibrium for PV commuters:
\begin{equation}
    \label{eq:linearsystem_S2and4_2_B1+B2}
    \begin{cases}
        \frac{\partial C^\mathcal{P}(t)}{\partial t}=(\alpha-\beta)\big(\dot{w}_\mathcal{H}(t)+\dot{w}_\mathcal{C}^\mathcal{P}(t+w_\mathcal{H}(t))(1+\dot{w}_\mathcal{H}(t))\big)-\beta=0,\\
        \beta(t^*-t_0^\mathcal{P})+u^\mathcal{P}=\alpha(t^*-t_1^\mathcal{P})+u^\mathcal{P},\\
        \int_{{t}_0^\mathcal{P}}^{t_1^\mathcal{P}} \dot{A}_\mathcal{H}(t) dt=N-N^\mathcal{R},\\
        \dot{A}_\mathcal{C}^\mathcal{P}(t)=s_\mathcal{H}.
    \end{cases}
\end{equation}
\end{itemize}

Combining Eqs. \eqref{eq:linearsystem_S2and4_1}-\eqref{eq:linearsystem_S2and4_2_B1+B2} with the inter-equilibrium conditions of Eq. \eqref{eq:jointUE}, we can derive the user equilibrium results for Scenarios 2 and 4 as follows,
\begin{itemize}
    \item Scenario 2: \begin{equation}
    \label{eq:UEresults_S2}
    \begin{cases}
        t_0^\mathcal{R} = t^*-\frac{N^\mathcal{R}}{s_\mathcal{C}^\mathcal{R}},\\
        t_1^\mathcal{R} = t^* -\frac{\beta}{\alpha+\pi}\frac{N^\mathcal{R}}{s_\mathcal{C}^\mathcal{R}},\\
        t_0^\mathcal{P} = t^* -\frac{N-N^\mathcal{R}}{s_\mathcal{C}^\mathcal{P}},\\
        t_1^\mathcal{P} = t^* -\frac{\beta}{\alpha}\frac{N-N^\mathcal{R}}{s_\mathcal{C}^\mathcal{P}},\\
        
        N^\mathcal{R} = \frac{N\beta+(u^\mathcal{P}-c^\mathcal{R})s_\mathcal{C}^\mathcal{P}}{\beta(s_\mathcal{C}^\mathcal{R}+s_\mathcal{C}^\mathcal{P})}s_\mathcal{C}^\mathcal{R}.\\
    \end{cases}
    \begin{cases}
        \dot{A}_\mathcal{C}^\mathcal{R}(t)=\frac{(\alpha+\pi)s_\mathcal{C}^\mathcal{R}}{\alpha+\pi-\beta},& t\in[t_0^\mathcal{R},t_1^\mathcal{R}],\\
        
        \dot{A}_\mathcal{C}^\mathcal{P}(t)= \frac{\alpha s_\mathcal{C}^\mathcal{P}}{\alpha-\beta},&  t\in[t_0^\mathcal{P},t_1^\mathcal{P}].\\
    \end{cases}
    \end{equation}
    \item Scenario 4:\begin{equation}
    \label{eq:UEresults_S4}
    \begin{cases}
        t_0^\mathcal{R} = t^*-\frac{N^\mathcal{R}}{s_\mathcal{C}^\mathcal{R}},\\
        t_1^\mathcal{R} = t^* -\frac{\beta}{\alpha+\pi}\frac{N^\mathcal{R}}{s_\mathcal{C}^\mathcal{R}},\\
        
        t_0^\mathcal{P} = t^* -\frac{N-N^\mathcal{R}}{s_\mathcal{C}^\mathcal{P}},\\
        t_1^\mathcal{P} = t^* -\frac{\beta}{\alpha}\frac{N-N^\mathcal{R}}{s_\mathcal{C}^\mathcal{P}},\\

        N^\mathcal{R} = \frac{N\beta+(u^\mathcal{P}-c^\mathcal{R})s_\mathcal{C}^\mathcal{P}}{\beta(s_\mathcal{C}^\mathcal{P}+s_\mathcal{C}^\mathcal{R})}s_\mathcal{C}^\mathcal{R}.\\
    \end{cases}
    \begin{cases}
        \dot{A}_\mathcal{H}(t)= \frac{\alpha s_\mathcal{C}^\mathcal{P}}{\alpha-\beta}, t\in[t_0^\mathcal{P},t_1^\mathcal{P}],\\
        
        \dot{A}_\mathcal{C}^\mathcal{R}(t)=\frac{(\alpha+\pi)s_\mathcal{C}^\mathcal{R}}{\alpha+\pi-\beta}, t\in[t_0^\mathcal{R},t_1^\mathcal{R}],\\
        
        \dot{A}_\mathcal{C}^\mathcal{P}(t)=s_\mathcal{H}, t\in[t_0^\mathcal{P},t_1^\mathcal{P}+w_\mathcal{H}(t_1^\mathcal{P})].\\
    \end{cases}
    \end{equation}
\end{itemize}

In Scenarios 2 and 4, congestion occurs only at the curbside bottleneck at the start of the morning peak. According to Proposition \ref{Proposition 1}, we have,
\begin{equation}
\label{eq:AppendixB_conditions_S2and4_1}
        \frac{s_\mathcal{C}^\mathcal{R}}{s_\mathcal{H} }\leq \frac{\alpha+\pi-\beta}{\alpha+\pi}.
\end{equation}

Since RV and PV are used with non-overlapping departure intervals, implying that the condition $t_0^\mathcal{P}\geq t_1^\mathcal{R}$ holds. Combining it with Eqs. \eqref{eq:UEresults_S2} and \eqref{eq:UEresults_S4}, one can obtain,
\begin{equation}
\label{eq:AppendixB_conditions_S2and4}
    (u^\mathcal{P}-c^\mathcal{R})\in\big[\frac{N\beta(\alpha+\pi-\beta)}{\beta s_\mathcal{C}^\mathcal{P}+(\alpha+\pi)s_\mathcal{C}^\mathcal{R}},\frac{\beta N}{s_\mathcal{C}^\mathcal{R}}\big).
\end{equation}

Different from Scenario 2, queues form at the highway bottleneck after the first PV commuter departs from home in Scenario 4, indicated by the inequality of $\dot{A}_\mathcal{H}(t)>s_\mathcal{H}$. Consequently, we can obtain the following condition to differentiate Scenario 4 from Scenario 2,
\begin{equation}
    \label{eq:distinct_condition_2and4}
    \frac{s_\mathcal{C}^\mathcal{P}}{s_\mathcal{H}}>\frac{\alpha-\beta}{\alpha}.
\end{equation}

Therefore, Scenario 4 occurs if Eqs. \eqref{eq:AppendixB_conditions_S2and4_1}-\eqref{eq:distinct_condition_2and4} are satisfied; otherwise, if Eqs. \eqref{eq:AppendixB_conditions_S2and4_1} and \eqref{eq:AppendixB_conditions_S2and4} hold but Eq. \eqref{eq:distinct_condition_2and4} is not met, Scenario 2 occurs.

\subsubsection*{(3) Scenario 3}
In Scenario 3, RV and PV are used with overlapping departure intervals. In addition, the queue forms solely at the curbside bottlenecks during the morning peak. Consequently, the following equation system ensures the intra-equilibrium for RV commuters:
\begin{equation}
    \label{eq:linearsystem_S3_RV}
    \begin{cases}
        \frac{\partial c^\mathcal{R}(t)}{\partial t}=(\alpha+\pi-\beta)\dot{w}_\mathcal{C}^\mathcal{R}(t)-\beta=0,\\
        \beta(t^*-t_0^\mathcal{R})+c^\mathcal{R}=(\alpha+\pi)(t^*-t_1^\mathcal{R})+c^\mathcal{R},\\
        (t^*-t_0^\mathcal{R})s_\mathcal{C}^\mathcal{R}=N^\mathcal{R}.\\
    \end{cases}
\end{equation}

The queuing delay for PVs is affected by the RV queue. As a result, the following equation system ensures PV commuters' intra-equilibrium:
\begin{equation}
    \label{eq:linearsystem_S3_PV}
    \begin{cases}
        \frac{\partial C^\mathcal{P}(t)}{\partial t}=(\alpha-\beta)\dot{w}_\mathcal{C}^\mathcal{P}(t)-\beta=0,\\
        \beta(t^*-t_0^\mathcal{P})+u^\mathcal{P}=\alpha(t^*-t_1^\mathcal{P})+u^\mathcal{P},\\
        \int_{{t}_0^\mathcal{P}}^{t_1^\mathcal{P}} \dot{A}_\mathcal{C}^\mathcal{P}(t) dt=N-N^\mathcal{R}.
    \end{cases}
\end{equation}

Whether RV or PV commuters depart during the final stage depends on the departure times of the latest RV and the latest PV commuters. If $t_1^\mathcal{P}\leq t_1^\mathcal{R}$, then only RV commuters depart after $t_1^\mathcal{P}$; if $t_1^\mathcal{P}>t_1^\mathcal{R}$, then only PV commuters depart after $t_1^\mathcal{R}$. For the two cases, by combining Eqs. \eqref{eq:linearsystem_S3_RV} and \eqref{eq:linearsystem_S3_PV} with the inter-equilibrium conditions of Eq. \eqref{eq:jointUE}, we can derive the user equilibrium results for Scenario 3 as follows,
\begin{itemize}
    \item If $t_1^\mathcal{P}\leq t_1^\mathcal{R}$, then only RV commuters depart after $t_1^\mathcal{P}$, we have
        \begin{equation}
            \label{eq:UEresults_S3_1}
            \begin{cases}
                t_0^\mathcal{R} = t^*-\frac{N^\mathcal{R}}{s_\mathcal{C}^\mathcal{R}},\\
                t_1^\mathcal{R} = t^* -\frac{\beta}{\alpha+\pi}\frac{N^\mathcal{R}}{s_\mathcal{C}^\mathcal{R}},\\
                t_0^\mathcal{P} = t^* -\frac{N-N^\mathcal{R}}{s_\mathcal{C}^\mathcal{P}-\delta^\mathcal{R}\frac{(\alpha+\pi)(\alpha-\beta)}{\alpha(\alpha+\pi-\beta)}s_\mathcal{C}^\mathcal{R}},\\
                t_1^\mathcal{P} = t^* -\frac{\beta}{\alpha}\frac{N-N^\mathcal{R}}{s_\mathcal{C}^\mathcal{P}-\delta^\mathcal{R}\frac{(\alpha+\pi)(\alpha-\beta)}{\alpha(\alpha+\pi-\beta)}s_\mathcal{C}^\mathcal{R}},\\

                N^\mathcal{R} = \frac{N\beta+(u^\mathcal{P}-c^\mathcal{R})\big(s_\mathcal{C}^\mathcal{P}-\frac{\delta^\mathcal{R}(\alpha+\pi)(\alpha-\beta)}{\alpha(\alpha+\pi-\beta)}s_\mathcal{C}^\mathcal{R}\big)}{\beta\big(s_\mathcal{C}^\mathcal{R}+s_\mathcal{C}^\mathcal{P}-\frac{\delta^\mathcal{R}(\alpha+\pi)(\alpha-\beta)}{\alpha(\alpha+\pi-\beta)}s_\mathcal{C}^\mathcal{R}\big)}s_\mathcal{C}^\mathcal{R}.
            \end{cases}
            \begin{cases}
                \dot{A}_\mathcal{C}^\mathcal{P}(t)= 
                \frac{\alpha s_\mathcal{C}^\mathcal{P}}{\alpha-\beta}-\frac{\delta^\mathcal{R}(\alpha+\pi)s_\mathcal{C}^\mathcal{R}}{\alpha+\pi-\beta},& t\in[t_0^\mathcal{P},t_1^\mathcal{P}],\\
                \dot{A}_\mathcal{C}^\mathcal{R}(t)=\frac{(\alpha+\pi)s_\mathcal{C}^\mathcal{R}}{\alpha+\pi-\beta}, &t\in[t_0^\mathcal{R},t_1^\mathcal{R}].
            \end{cases}
        \end{equation}
    \item If $t_1^\mathcal{P}> t_1^\mathcal{R}$, then only PV commuters depart after $t_1^\mathcal{R}$, we have
        \begin{equation}
            \label{eq:UEresults_S3_2}
            \begin{cases}
                t_0^\mathcal{R} = t^*-\frac{N^\mathcal{R}}{s_\mathcal{C}^\mathcal{R}},\\
                t_1^\mathcal{R} = t^* -\frac{\beta}{\alpha+\pi}\frac{N^\mathcal{R}}{s_\mathcal{C}^\mathcal{R}},\\
                t_0^\mathcal{P} = t^* -\frac{N-\frac{\alpha+\pi-\beta(1-\delta^\mathcal{R})}{\alpha+\pi-\beta}N^\mathcal{R}}{s_\mathcal{C}^\mathcal{P}-\delta^\mathcal{R}\frac{\alpha+\pi}{\alpha+\pi-\beta}s_\mathcal{C}^\mathcal{R}},\\
                t_1^\mathcal{P} = t^* -\frac{\beta}{\alpha}\frac{N-\frac{\alpha+\pi-\beta(1-\delta^\mathcal{R})}{\alpha+\pi-\beta}N^\mathcal{R}}{s_\mathcal{C}^\mathcal{P}-\delta^\mathcal{R}\frac{\alpha+\pi}{\alpha+\pi-\beta}s_\mathcal{C}^\mathcal{R}},\\
        
                 N^\mathcal{R} = \frac{N\beta+(u^\mathcal{P}-c^\mathcal{R})\big(s_\mathcal{C}^\mathcal{P}-\frac{(\alpha+\pi)\delta^\mathcal{R}}{\alpha+\pi-\beta}s_\mathcal{C}^\mathcal{R}\big)}
                {\beta\big(s_\mathcal{C}^\mathcal{P}+s_\mathcal{C}^\mathcal{R}(1-\delta^\mathcal{R})\big)}s_\mathcal{C}^\mathcal{R}.\\
            \end{cases}
            \begin{cases}
                \dot{A}_\mathcal{C}^\mathcal{P}(t)= \begin{cases}
                \frac{\alpha s_\mathcal{C}^\mathcal{P}}{\alpha-\beta}-\frac{\delta^\mathcal{R}(\alpha+\pi)s_\mathcal{C}^\mathcal{R}}{\alpha+\pi-\beta},& t\in[t_0^\mathcal{P},t_1^\mathcal{R}),\\
                \frac{\alpha s_\mathcal{C}^\mathcal{P}}{\alpha-\beta}& t\in[t_1^\mathcal{R},t_1^\mathcal{P}],
                \end{cases}\\
                \dot{A}_\mathcal{C}^\mathcal{R}(t)=\frac{(\alpha+\pi)s_\mathcal{C}^\mathcal{R}}{\alpha+\pi-\beta}, t\in[t_0^\mathcal{R},t_1^\mathcal{R}].
            \end{cases}
        \end{equation}
\end{itemize}

During the co-departure period when both RV and PV commuters depart from home, the total departure rate of them can not exceed the service rate of the highway bottleneck in Scenario 3. It implies that $\dot{A}_\mathcal{C}^\mathcal{R}(t) + \dot{A}_\mathcal{C}^\mathcal{P}(t) \leq s_\mathcal{H}$ should hold for the co-departure period. Therefore, we have the following condition:
\begin{equation}
    \label{eq:distinct_condition_3}
    \frac{s_\mathcal{C}^\mathcal{P}}{s_\mathcal{H}}\leq\frac{\alpha-\beta}{\alpha}-\frac{(1-\delta^\mathcal{R})(\alpha+\pi)(\alpha-\beta)s_\mathcal{C}^\mathcal{R}}{\alpha(\alpha+\pi-\beta)s_\mathcal{H}}.
\end{equation}

In Scenario 3, the queue forms solely at the curbside bottleneck at the start of the morning peak. According to Proposition \ref{Proposition 1}, we have,
\begin{equation}
\label{eq:AppendixB_conditions_S3}
    \frac{s_\mathcal{C}^\mathcal{R}}{s_\mathcal{H} }\leq \frac{\alpha+\pi-\beta}{\alpha+\pi}.
\end{equation}

RV and PV are used with overlapping departure intervals, implying that the condition $t_0^\mathcal{P}< t_1^\mathcal{R}$ holds. Combining it with Eqs. \eqref{eq:UEresults_S3_1} and \eqref{eq:UEresults_S3_2}, one can obtain,
\begin{equation}
\label{eq:AppendixB_conditions_S3_2}
    (u^\mathcal{P}-c^\mathcal{R})\in\big(0,\frac{N\beta(\alpha+\pi-\beta)}{\beta s_\mathcal{C}^\mathcal{P}+(\alpha+\pi)s_\mathcal{C}^\mathcal{R}}\big).
\end{equation}

Therefore, Eqs. \eqref{eq:distinct_condition_3}-\eqref{eq:AppendixB_conditions_S3_2} present the occurrence conditions for Scenario 3. 

\subsubsection*{(4) Scenario 5}
In Scenario 5, RV and PV are used with overlapping departure intervals. In the initial stage $(t_0^\mathcal{R},t_0^\mathcal{P})$, only RV commuters depart from home and the queue of RV forms exclusively at the curbside bottleneck. Therefore, we have the following condition:
\begin{equation}
    \label{eq:linearsystem_S5_beginningphase}
    \begin{aligned}
        \frac{\partial C^\mathcal{R}(t)}{\partial t}=(\alpha+\pi-\beta)\dot{w}_\mathcal{C}^\mathcal{R}(t)-\beta=0.\\
    \end{aligned}
\end{equation}

During the co-departure stage, RV and PV commuters depart from home, and queues form at both the highway and the curbside bottlenecks. Therefore, we have the following conditions:
\begin{equation}
    \label{eq:linearsystem_S5_codepartphase}
\begin{cases}
        \frac{\partial C^\mathcal{R}(t)}{\partial t}=(\alpha+\pi-\beta)\big(\dot{w}_\mathcal{H}(t)+\dot{w}_\mathcal{C}^\mathcal{R}(t+w_\mathcal{H}(t))(1+\dot{w}_\mathcal{H}(t))\big)-\beta=0,\\[1.2ex]
        \frac{\partial C^\mathcal{P}(t)}{\partial t}=(\alpha-\beta)\big(\dot{w}_\mathcal{H}(t)+\dot{w}_\mathcal{C}^\mathcal{P}(t+w_\mathcal{H}(t))(1+\dot{w}_\mathcal{H}(t))\big)-\beta=0,\\[1.2ex]
        \dot{A}_\mathcal{C}^\mathcal{R}(t)+\dot{A}_\mathcal{C}^\mathcal{P}(t)=s_\mathcal{H}.
\end{cases}
\end{equation}

During the final stage, the highway bottleneck and the curbside bottleneck are still congested. There are two possible cases depending on the departure times of the latest RV and the latest PV commuters:
\begin{itemize}
\item if $t_1^\mathcal{P}\leq t_1^\mathcal{R}$, then only RV commuters depart after $t_1^\mathcal{P}$, we thus have,
        \begin{equation}
        \label{eq:finalphase_t1c<t1r}
        \begin{cases}
            \frac{\partial C^\mathcal{R}(t)}{\partial t}=(\alpha+\pi-\beta)\big(\dot{w}_\mathcal{H}(t)+\dot{w}_\mathcal{C}^\mathcal{R}(t+w_\mathcal{H}(t))(1+\dot{w}_\mathcal{H}(t))\big)-\beta=0,\\
            \dot{A}_\mathcal{C}^\mathcal{R}(t)=s_\mathcal{H}.
        \end{cases}
        \end{equation}
\item if $t_1^\mathcal{P}>t_1^\mathcal{R}$, then only PV commuters depart after $t_1^\mathcal{R}$, we thus have,
    \begin{equation}
    \label{eq:finalphase_t1c>t1r}    
    \begin{cases}
        \frac{\partial C^\mathcal{P}(t)}{\partial t}=(\alpha-\beta)\big(\dot{w}_\mathcal{H}(t)+\dot{w}_\mathcal{C}^\mathcal{P}(t+w_\mathcal{H}(t))(1+\dot{w}_\mathcal{H}(t))\big)-\beta=0,\\
        \dot{A}_\mathcal{C}^\mathcal{P}(t)=s_\mathcal{H}.
    \end{cases}
    \end{equation}
\end{itemize}

Based on Eqs. \eqref{eq:linearsystem_S5_beginningphase}-\eqref{eq:finalphase_t1c>t1r}, we can derive the equilibrium arrival rates of RV and PV commuters at bottlenecks as:
\begin{itemize}
    \item if $t_1^\mathcal{P}\leq t_1^\mathcal{R}$,
    \begin{equation}
    \begin{aligned}
        &\dot{A}_\mathcal{H}(t) = \begin{cases}
            \frac{(\alpha+\pi)s_\mathcal{C}^\mathcal{R}}{\alpha+\pi-\beta},  & t\in[t_0^\mathcal{R},t_0^\mathcal{P})\cup(t_1^\mathcal{P},t_1^\mathcal{R}] ,\\[1.2ex]
            \frac{\alpha s_\mathcal{C}^\mathcal{P}}{\alpha-\beta}+\frac{(1-\delta^\mathcal{R})(\alpha+\pi)s_\mathcal{C}^\mathcal{R}}{(\alpha+\pi-\beta)}, & t\in[t_0^\mathcal{P},t_1^\mathcal{P}].
        \end{cases}\\
        &\dot{A}_\mathcal{C}^\mathcal{R}(t) = \begin{cases}
            \frac{(\alpha+\pi)s_\mathcal{C}^\mathcal{R}}{\alpha+\pi-\beta},  & t\in[t_0^\mathcal{R},t_0^\mathcal{P}), \\
            \frac{(\alpha+\pi)(\alpha-\beta)s_\mathcal{C}^\mathcal{R}s_\mathcal{H}}{\alpha(\alpha+\pi-\beta)s_\mathcal{C}^\mathcal{P}+(1-\delta^\mathcal{R})(\alpha+\pi)(\alpha-\beta)s_\mathcal{C}^\mathcal{R}}, & t\in\big[t_0^\mathcal{P},t_1^\mathcal{P}+w_\mathcal{H}(t_1^\mathcal{P})\big],\\
            s_\mathcal{H}, & t\in\big(t_1^\mathcal{P}+w_\mathcal{H}(t_1^\mathcal{P}),t_1^\mathcal{R}+w_\mathcal{H}(t_1^\mathcal{R})\big].
        \end{cases}\\
        &\dot{A}_\mathcal{C}^\mathcal{P}(t)=s_\mathcal{H}-\frac{(\alpha+\pi)(\alpha-\beta)s_\mathcal{C}^\mathcal{R}s_\mathcal{H}}{\alpha(\alpha+\pi-\beta)s_\mathcal{C}^\mathcal{P}+(1-\delta^\mathcal{R})(\alpha+\pi)(\alpha-\beta)s_\mathcal{C}^\mathcal{R}},t\in\big[t_0^\mathcal{P},t_1^\mathcal{P}+w_\mathcal{H}(t_1^\mathcal{P})\big].
    \end{aligned}
    \label{eq:arrivalrate_S5_t1c<t1r}
    \end{equation}

    \item if $t_1^\mathcal{P}>t_1^\mathcal{R}$,
    \begin{equation}
    \begin{aligned}
        &\dot{A}_\mathcal{H}(t) = \begin{cases}
            \frac{(\alpha+\pi)s_\mathcal{C}^\mathcal{R}}{\alpha+\pi-\beta}, & t\in[t_0^\mathcal{R},t_0^\mathcal{P}),\\[1.2ex]
            
            \frac{\alpha s_\mathcal{C}^\mathcal{P}}{\alpha-\beta}+\frac{(1-\delta^\mathcal{R})(\alpha+\pi)s_\mathcal{C}^\mathcal{R}}{(\alpha+\pi-\beta)}, & t\in[t_0^\mathcal{P},t_1^\mathcal{R}],\\[1.2ex]

            \frac{\alpha s_\mathcal{C}^\mathcal{P}}{\alpha-\beta}, & t\in[t_1^\mathcal{R},t_1^\mathcal{P}].
        \end{cases}\\
        &\dot{A}_\mathcal{C}^\mathcal{R}(t) = \begin{cases}
            \frac{(\alpha+\pi)s_\mathcal{C}^\mathcal{R}}{\alpha+\pi-\beta}, & t\in[t_0^\mathcal{R},t_0^\mathcal{P}),\\
            \frac{(\alpha+\pi)(\alpha-\beta)s_\mathcal{C}^\mathcal{R}s_\mathcal{H}}{\alpha(\alpha+\pi-\beta)s_\mathcal{C}^\mathcal{P}+(1-\delta^\mathcal{R})(\alpha+\pi)(\alpha-\beta)s_\mathcal{C}^\mathcal{R}}, & t\in\big[t_0^\mathcal{P},t_1^\mathcal{R}+w_\mathcal{H}(t_1^\mathcal{R})\big].\\
        \end{cases}\\
        &\dot{A}_\mathcal{C}^\mathcal{P}(t)=\begin{cases}
            s_\mathcal{H}-\frac{(\alpha+\pi)(\alpha-\beta)s_\mathcal{C}^\mathcal{R}s_\mathcal{H}}{\alpha(\alpha+\pi-\beta)s_\mathcal{C}^\mathcal{P}+(1-\delta^\mathcal{R})(\alpha+\pi)(\alpha-\beta)s_\mathcal{C}^\mathcal{R}}, & t\in\big[t_0^\mathcal{P},t_1^\mathcal{R}+w_\mathcal{H}(t_1^\mathcal{R})\big),\\
            s_\mathcal{H}, & t\in\big[t_1^\mathcal{R}+w_\mathcal{H}(t_1^\mathcal{R}),t_1^\mathcal{P}+w_\mathcal{H}(t_1^\mathcal{P})\big].\\
        \end{cases}
    \end{aligned}
    \label{eq:arrivalrate_S5_t1c>t1r}
 \end{equation}
\end{itemize}

At user equilibrium, all commuters have the same generalized travel cost, i.e., $C^\mathcal{R}(t_0^\mathcal{R})=C^\mathcal{R}(t_1^\mathcal{R})=C^\mathcal{P}(t_0^\mathcal{P})=C^\mathcal{P}(t_1^\mathcal{P})$. All RV commuters arrive at CBD within the time interval $(t_0^\mathcal{R},t^*)$ at a rate of $s_\mathcal{C}^\mathcal{R}$. All PV commuters pass through the highway bottleneck and arrive at the curbside bottleneck within the time interval $(t_0^\mathcal{P},t_1^\mathcal{P}+w_\mathcal{H}(t_1^\mathcal{P}))$. Therefore, we have the following system of equations:
\begin{equation}
        \begin{cases}
            \beta(t^*-t_0^\mathcal{R})+c^\mathcal{R}=\beta(t^*-t_0^\mathcal{P})+u^\mathcal{P},\\
            \beta(t^*-t_0^\mathcal{R})+c^\mathcal{R}=(\alpha+\pi)(t^*-t_1^\mathcal{R})+c^\mathcal{R},\\
            \beta(t^*-t_0^\mathcal{P})+u^\mathcal{P}=\alpha(t^*-t_1^\mathcal{P})+u^\mathcal{P},\\
            (t^*-t_0^\mathcal{R})s_\mathcal{C}^\mathcal{R}=N^\mathcal{R},\\[0.8ex]
            \int_{{t}_0^\mathcal{P}}^{t_1^\mathcal{P}+w_\mathcal{H}(t_1^\mathcal{P})} \dot{A}_\mathcal{C}^\mathcal{P}(t) dt=N-N^\mathcal{R}.
        \end{cases}
    \label{eq:linearsystem_S5}
\end{equation}

Substituting Eqs. \eqref{eq:arrivalrate_S5_t1c<t1r} and \eqref{eq:arrivalrate_S5_t1c>t1r} into the Eqs. \eqref{eq:linearsystem_S5}, we can obtain the equilibrium results for Scenario 5:
\begin{itemize}
    \item if $t_1^\mathcal{P}\leq t_1^\mathcal{R}$,
    \begin{equation}
        \begin{cases}\label{eq:S5_UEresults_1}
            t_0^\mathcal{R}=t^*-\frac{N^\mathcal{R}}{s_\mathcal{C}^\mathcal{R}},\\
            t_1^\mathcal{R}=t^*-\frac{\beta}{\alpha+\pi}\frac{N^\mathcal{R}}{s_\mathcal{C}^\mathcal{R}},\\

            t_0^\mathcal{P} = t^* -\frac{N-N^\mathcal{R}}{s_\mathcal{C}^\mathcal{P}-\delta^\mathcal{R}\frac{(\alpha+\pi)(\alpha-\beta)}{\alpha(\alpha+\pi-\beta)}s_\mathcal{C}^\mathcal{R}},\\
             t_1^\mathcal{P} = t^* -\frac{\beta}{\alpha}\frac{N-N^\mathcal{R}}{s_\mathcal{C}^\mathcal{P}-\delta^\mathcal{R}\frac{(\alpha+\pi)(\alpha-\beta)}{\alpha(\alpha+\pi-\beta)}s_\mathcal{C}^\mathcal{R}},\\

            N^\mathcal{R} = \frac{N\beta+(u^\mathcal{P}-c^\mathcal{R})\big(s_\mathcal{C}^\mathcal{P}-\frac{\delta^\mathcal{R}(\alpha+\pi)(\alpha-\beta)}{\alpha(\alpha+\pi-\beta)}s_\mathcal{C}^\mathcal{R}\big)}{\beta\big(s_\mathcal{C}^\mathcal{R}+s_\mathcal{C}^\mathcal{P}-\frac{\delta^\mathcal{R}(\alpha+\pi)(\alpha-\beta)}{\alpha(\alpha+\pi-\beta)}s_\mathcal{C}^\mathcal{R}\big)}s_\mathcal{C}^\mathcal{R}.
        \end{cases}
    \end{equation}
    
    \item if $t_1^\mathcal{P}>t_1^\mathcal{R}$,
    \begin{equation}\label{eq:S5_UEresults_2}
        \begin{cases}
            t_0^\mathcal{R}=t^*-\frac{N^\mathcal{R}}{s_\mathcal{C}^\mathcal{R}},\\
            t_1^\mathcal{R}=t^*-\frac{\beta}{\alpha+\pi}\frac{N^\mathcal{R}}{s_\mathcal{C}^\mathcal{R}},\\
            
            t_0^\mathcal{P} = t^* -\frac{N-\frac{\alpha+\pi-\beta(1-\delta^\mathcal{R})}{\alpha+\pi-\beta}N^\mathcal{R}}{s_\mathcal{C}^\mathcal{P}-\delta^\mathcal{R}\frac{\alpha+\pi}{\alpha+\pi-\beta}s_\mathcal{C}^\mathcal{R}},\\
                t_1^\mathcal{P} = t^* -\frac{\beta}{\alpha}\frac{N-\frac{\alpha+\pi-\beta(1-\delta^\mathcal{R})}{\alpha+\pi-\beta}N^\mathcal{R}}{s_\mathcal{C}^\mathcal{P}-\delta^\mathcal{R}\frac{\alpha+\pi}{\alpha+\pi-\beta}s_\mathcal{C}^\mathcal{R}},\\
        
                 N^\mathcal{R} = \frac{N\beta+(u^\mathcal{P}-c^\mathcal{R})\big(s_\mathcal{C}^\mathcal{P}-\frac{(\alpha+\pi)\delta^\mathcal{R}}{\alpha+\pi-\beta}s_\mathcal{C}^\mathcal{R}\big)}
                {\beta\big(s_\mathcal{C}^\mathcal{P}+s_\mathcal{C}^\mathcal{R}(1-\delta^\mathcal{R})\big)}s_\mathcal{C}^\mathcal{R}.\\
        \end{cases}
    \end{equation}
\end{itemize}

Similar to Scenario 3, a queue forms solely at the curbside bottleneck at the start of the morning peak. According to Proposition \ref{Proposition 1}, we have,
\begin{equation}
\label{eq:AppendixB_conditions_S5}
    \frac{s_\mathcal{C}^\mathcal{R}}{s_\mathcal{H} }\leq \frac{\alpha+\pi-\beta}{\alpha+\pi}.
\end{equation}

RV and PV are used with overlapping departure intervals, implying that the condition $t_0^\mathcal{P}< t_1^\mathcal{R}$ holds. Combining it with Eqs. \eqref{eq:S5_UEresults_1} and \eqref{eq:S5_UEresults_2}, one can obtain,
\begin{equation}
\label{eq:AppendixB_conditions_S5_2}
    (u^\mathcal{P}-c^\mathcal{R})\in\big(0,\frac{N\beta(\alpha+\pi-\beta)}{\beta s_\mathcal{C}^\mathcal{P}+(\alpha+\pi)s_\mathcal{C}^\mathcal{R}}\big).
\end{equation}

Different from Scenario 3, during the co-departure stage when both RV and PV commuters depart from home, the condition $\dot{A}_\mathcal{H}(t)>s_\mathcal{H}$ holds to ensure that queues form at the highway bottleneck. We thus have,
\begin{equation}
    \label{eq:distinct_condition_5}
    \frac{s_\mathcal{C}^\mathcal{P}}{s_\mathcal{H}}>\frac{\alpha-\beta}{\alpha}-\frac{(1-\delta^\mathcal{R})(\alpha+\pi)(\alpha-\beta)s_\mathcal{C}^\mathcal{R}}{\alpha(\alpha+\pi-\beta)s_\mathcal{H}}.
\end{equation}

Therefore, Eqs. \eqref{eq:AppendixB_conditions_S5}-\eqref{eq:distinct_condition_5} give the occurrence condition for Scenario 5. 

\subsubsection*{(6) Scenario 6}
In Scenario 6, only the RV mode will be used, and queues will form at both the highway bottleneck and the curbside bottleneck during the morning peak. Therefore, the arrival rate of RVs at the curbside bottleneck is equal to the service rate of the highway bottleneck, i.e., $\dot {A}_\mathcal{C}^\mathcal{R}(t)=s_\mathcal{H}$. The arrival rate of RVs at CBD is equal to the service rate of the curbside bottleneck $s_\mathcal{C}^\mathcal{R}$. Therefore, the equilibrium results for Scenario 6 can be obtained by solving the following equation system:
\begin{equation}
    \begin{cases}
        \frac{\partial C^\mathcal{R}(t)}{\partial t}=(\alpha+\pi-\beta)\big(\dot{w}_\mathcal{H}(t)+\dot{w}_\mathcal{C}^\mathcal{R}(t+w_\mathcal{H}(t))(1+\dot{w}_\mathcal{H}(t))\big)-\beta=0,\\
        \beta(t^*-t_0^\mathcal{R})+c^\mathcal{R}=(\alpha+\pi)(t^*-t_1^\mathcal{R})+c^\mathcal{R},\\
        (t^*-t_0^\mathcal{R})s_\mathcal{C}^\mathcal{R}=N.
    \end{cases}
\end{equation}
The equilibrium results are as follows:
\begin{equation}
    \begin{cases}
        t_0^\mathcal{R}=t^*-\frac{N}{s_\mathcal{C}^\mathcal{R}},\\
        t_1^\mathcal{R}=t^*-\frac{\beta N}{(\alpha+\pi)s_\mathcal{C}^\mathcal{R}}.\\
    \end{cases}
    \begin{cases}
        \dot{A}_\mathcal{H}(t)=\frac{(\alpha+\pi)s_\mathcal{C}^\mathcal{R}}{\alpha+\pi-\beta}, &t\in[t_0^\mathcal{R},t_1^\mathcal{R}],\\
        \dot{A}_\mathcal{C}^\mathcal{R}(t) = s_\mathcal{H}, &t\in[t_0^\mathcal{R},t_1^\mathcal{R}+w_\mathcal{H}(t_1^\mathcal{R})].
    \end{cases}
\end{equation}

Similar to Scenario 1, the following condition ensures that only RV is used during the morning peak, 
\begin{equation}
\label{eq:AppendixB_conditions_S6_1}
        (u^\mathcal{P}-c^\mathcal{R})\geq\frac{\beta N}{s_\mathcal{C}^\mathcal{R}}.
\end{equation}

However, congestion occurs at the highway and the curbside bottlenecks in Scenario 6. According to Proposition \ref{Proposition 1}, we have,
\begin{equation}
\label{eq:AppendixB_conditions_S6_2}
        \frac{s_\mathcal{C}^\mathcal{R}}{s_\mathcal{H} }> \frac{\alpha+\pi-\beta}{\alpha+\pi}.
\end{equation}

Therefore, Eqs. \eqref{eq:AppendixB_conditions_S6_1} and \eqref{eq:AppendixB_conditions_S6_2} give the occurrence conditions for Scenario 6.

\subsubsection*{(7) Scenario 7}
In Scenario 7, RV and PV are utilized with non-overlapping departure intervals. During the initial stage of $(t_0^\mathcal{R}, t_1^\mathcal{R})$, only RV commuters depart from home, and the RV queue forms at both the highway and the curbside bottlenecks. The arrival rate of RVs at the curbside bottleneck equals the service rate of the highway bottleneck, i.e., $\dot{A}_\mathcal{C}^\mathcal{R}(t)=s_\mathcal{H}$. In addition, the arrival rate of RVs at CBD is equal to the service rate of the curbside bottleneck $s_\mathcal{C}^\mathcal{R}$. Therefore, the intra-equilibrium for RV commuters can be expressed as the following system of equations:
\begin{equation}
    \label{eq:linearsystem_S7_RV}
    \begin{cases}
        \frac{\partial c^\mathcal{R}(t)}{\partial t}=(\alpha+\pi-\beta)\big(\dot{w}_\mathcal{H}(t)+\dot{w}_\mathcal{C}^\mathcal{R}(t+w_\mathcal{H}(t))(1+\dot{w}_\mathcal{H}(t))\big)-\beta=0,\\
        \beta(t^*-t_0^\mathcal{R})+c^\mathcal{R}=(\alpha+\pi)(t^*-t_1^\mathcal{R})+c^\mathcal{R},\\
        (t^*-t_0^\mathcal{R})s_\mathcal{C}^\mathcal{R}=N^\mathcal{R},\\
        \dot{A}_\mathcal{C}^\mathcal{R}(t)=s_\mathcal{H}.
    \end{cases}
\end{equation}

During the period $(t_0^\mathcal{P}, t_1^\mathcal{P})$ when only PV commuters leave home, the highway bottleneck might be congested. As $t_0^\mathcal{P}\leq t_1^\mathcal{R}+w_\mathcal{H}(t_1^\mathcal{R})$, the first PV commuter encounters a queuing delay of $w_\mathcal{H}(t_0^\mathcal{P})= \frac{A_\mathcal{H}(t_0^\mathcal{P}) - s_\mathcal{H}(t_0^\mathcal{P}-t_0^\mathcal{R})}{s_\mathcal{H}}$ at the highway bottleneck; otherwise, the first PV commuter faces no queue at the highway bottleneck, i.e., $w_\mathcal{H}(t_0^\mathcal{P})=0$. Therefore, the intra-equilibrium for PV commuters can be expressed as the following equation system:
\begin{itemize}
    \item if $t_0^\mathcal{P}\leq t_1^\mathcal{R}+w_\mathcal{H}(t_1^\mathcal{R})$,
    \begin{equation}
        \label{eq:linearsystem_S7_PV_1}
        \begin{cases}
            \frac{\partial C^\mathcal{P}(t)}{\partial t}=(\alpha-\beta)\big(\dot{w}_\mathcal{H}(t)+\dot{w}_\mathcal{C}^\mathcal{P}(t+w_\mathcal{H}(t))(1+\dot{w}_\mathcal{H}(t))\big)-\beta=0,\\
            (\alpha-\beta) w_\mathcal{H}(t_0^\mathcal{P})+\beta\big(t^*-t_0^\mathcal{P}\big)+u^\mathcal{P}=\alpha(t^*-t_1^\mathcal{P})+u^\mathcal{P},\\
            \int_{{t}_0^\mathcal{P}}^{t_1^\mathcal{P}} \dot{A}_\mathcal{H}(t) dt=N-N^\mathcal{R},\\
            \dot{A}_\mathcal{C}^\mathcal{P}(t)=s_\mathcal{H}.\\
        \end{cases}
    \end{equation}
    \item if $t_0^\mathcal{P}> t_1^\mathcal{R}+w_\mathcal{H}(t_1^\mathcal{R})$
    \begin{equation}\label{eq:linearsystem_S7_PV_2}
        \begin{cases}
            \frac{\partial C^\mathcal{P}(t)}{\partial t}=(\alpha-\beta)\big(\dot{w}_\mathcal{H}(t)+\dot{w}_\mathcal{C}^\mathcal{P}(t+w_\mathcal{H}(t))(1+\dot{w}_\mathcal{H}(t))\big)-\beta=0,\\
            \beta\big(t^*-t_0^\mathcal{P}\big)+u^\mathcal{P}=\alpha(t^*-t_1^\mathcal{P})+u^\mathcal{P},\\
            \int_{{t}_0^\mathcal{P}}^{t_1^\mathcal{P}} \dot{A}_\mathcal{H}(t) dt=N-N^\mathcal{R},\\
            \dot{A}_\mathcal{C}^\mathcal{P}(t)=s_\mathcal{H}.\\
        \end{cases}
    \end{equation}
\end{itemize}

Combining Eqs. \eqref{eq:linearsystem_S7_RV} and \eqref{eq:linearsystem_S7_PV_2} with the inter-equilibrium conditions of Eq. \eqref{eq:jointUE}, we can derive the user equilibrium results for Scenario 7:
\begin{itemize}
    \item if $t_0^\mathcal{P}\leq t_1^\mathcal{R}+w_\mathcal{H}(t_1^\mathcal{R})$,
\begin{equation}
    \label{eq:UEresults_S7_1}
    \begin{cases}
        t_0^\mathcal{R} = t^*-\frac{N^\mathcal{R}}{s_\mathcal{C}^\mathcal{R}},\\
        t_1^\mathcal{R} = t^* -\frac{\beta}{\alpha+\pi}\frac{N^\mathcal{R}}{s_\mathcal{C}^\mathcal{R}},\\
        
        t_0^\mathcal{P} = t^* -\frac{N^\mathcal{R}}{s_\mathcal{C}^\mathcal{R}} + \frac{s_\mathcal{H}(u^\mathcal{P}-c^\mathcal{R})(\alpha+\pi-\beta)}{\alpha(\alpha+\pi-\beta)s_\mathcal{H}-(\alpha+\pi)(\alpha-\beta)s_\mathcal{C}^\mathcal{R}},\\

        t_1^\mathcal{P} = t_0^\mathcal{P} + \frac{(N-N^\mathcal{R})(\alpha-\beta)}{\alpha s_\mathcal{C}^\mathcal{P}},\\
        
        N^\mathcal{R} = \frac{s_\mathcal{H}s_\mathcal{C}^\mathcal{R}}{s_\mathcal{H}(s_\mathcal{C}^\mathcal{R}+s_\mathcal{C}^\mathcal{P})-s_\mathcal{C}^\mathcal{P}s_\mathcal{C}^\mathcal{R}}N.\\
    \end{cases}
    \begin{cases}
        \dot{A}_\mathcal{H}(t)=\begin{cases}
            \frac{(\alpha+\pi)s_\mathcal{C}^\mathcal{R}}{\alpha+\pi-\beta}, & t\in(t_0^\mathcal{R},t_1^\mathcal{R}),\\
            \frac{\alpha s_\mathcal{C}^\mathcal{P}}{\alpha-\beta}, & t\in(t_0^\mathcal{P},t_1^\mathcal{P}),
        \end{cases}\\
        \dot{A}_\mathcal{C}^\mathcal{R}(t)=s_\mathcal{H},\  t\in(t_0^\mathcal{R},t_1^\mathcal{R}+w_\mathcal{H}(t_1^\mathcal{R})),\\
        \dot{A}_\mathcal{C}^\mathcal{P}(t)= s_\mathcal{H},\  t\in(t_0^\mathcal{P}+w_\mathcal{H}(t_0^\mathcal{P}),t_1^\mathcal{P}+w_\mathcal{H}(t_1^\mathcal{P})).
    \end{cases}
\end{equation}

\item if $t_0^\mathcal{P}> t_1^\mathcal{R}+w_\mathcal{H}(t_1^\mathcal{R})$,
\begin{equation}
    \label{eq:UEresults_S7_2}
    \begin{cases}
        t_0^\mathcal{R} = t^*-\frac{N^\mathcal{R}}{s_\mathcal{C}^\mathcal{R}},\\
        t_1^\mathcal{R} = t^* -\frac{\beta}{\alpha+\pi}\frac{N^\mathcal{R}}{s_\mathcal{C}^\mathcal{R}},\\
        
        t_0^\mathcal{P} = t^* -\frac{N-N^\mathcal{R}}{s_\mathcal{C}^\mathcal{P}},\\

        t_1^\mathcal{P} = t^* - \frac{\beta}{\alpha}\frac{N-N^\mathcal{R}}{s_\mathcal{C}^\mathcal{P}},\\
        
        N^\mathcal{R} = \frac{N\beta+(u^\mathcal{P}-c^\mathcal{R})s_\mathcal{C}^\mathcal{P}}{\beta(s_\mathcal{C}^\mathcal{R}+s_\mathcal{C}^\mathcal{P})}s_\mathcal{C}^\mathcal{R}.\\
    \end{cases}
    \begin{cases}
        \dot{A}_\mathcal{H}(t)=\begin{cases}
            \frac{(\alpha+\pi)s_\mathcal{C}^\mathcal{R}}{\alpha+\pi-\beta}, & t\in(t_0^\mathcal{R},t_1^\mathcal{R}),\\
            \frac{\alpha s_\mathcal{C}^\mathcal{P}}{\alpha-\beta}, & t\in(t_0^\mathcal{P},t_1^\mathcal{P}),
        \end{cases}\\
        \dot{A}_\mathcal{C}^\mathcal{R}(t)=s_\mathcal{H},\  t\in(t_0^\mathcal{R},t_1^\mathcal{R}+w_\mathcal{H}(t_1^\mathcal{R})),\\
        \dot{A}_\mathcal{C}^\mathcal{P}(t)= s_\mathcal{H},\  t\in(t_0^\mathcal{P}+w_\mathcal{H}(t_0^\mathcal{P}),t_1^\mathcal{P}+w_\mathcal{H}(t_1^\mathcal{P})).
    \end{cases}
\end{equation}
\end{itemize}

In Scenarios 7, congestion occurs at the highway and the curbside bottlenecks. According to Proposition \ref{Proposition 1}, we have,
\begin{equation}
\label{eq:AppendixB_conditions_S7_1}
        \frac{s_\mathcal{C}^\mathcal{R}}{s_\mathcal{H} }>\frac{\alpha+\pi-\beta}{\alpha+\pi}.
\end{equation}

RV and PV are used with non-overlapping departure intervals, implying that the condition $t_0^\mathcal{P}\geq t_1^\mathcal{R}$ holds. Combining it with Eqs. \eqref{eq:UEresults_S7_1} and \eqref{eq:UEresults_S7_2}, one can obtain,
\begin{equation}
\label{eq:AppendixB_conditions_S7_2}
    (u^\mathcal{P}-c^\mathcal{R})\in\bigg[\frac{N\big( \alpha(\alpha+\pi-\beta)s_\mathcal{H}-(\alpha+\pi)(\alpha-\beta)s_\mathcal{C}^\mathcal{R}\big)}{(\alpha+\pi)\big(s_\mathcal{H}(s_\mathcal{C}^\mathcal{R}+s_\mathcal{C}^\mathcal{P})-s_\mathcal{C}^\mathcal{R}s_\mathcal{C}^\mathcal{P}\big)},\frac{\beta N}{s_\mathcal{C}^\mathcal{R}}\bigg).
\end{equation}

Therefore, Eqs. \eqref{eq:AppendixB_conditions_S7_1} and \eqref{eq:AppendixB_conditions_S7_2} give the occurrence conditions for Scenario 7.

\subsubsection*{(8) Scenario 8}
In Scenario 8, queues form at both the highway bottleneck and the curbside bottlenecks. In addition, the departure intervals of RV and PV commuters from home overlap. We can divide the morning peak into three stages: the initial stage where only RV commuters leave home, the co-departure stage where both RV and PV commuters leave home, and the final stage where only RV commuters depart if $t_1^\mathcal{P}<t_1^\mathcal{R}$ holds or only PV commuters depart if $t_1^\mathcal{P}>t_1^\mathcal{R}$ holds.

During the initial stage of $(t_0^\mathcal{R},t_0^\mathcal{P})$, only RV commuters leave home, and queues form at both the highway and the curbside bottlenecks after the first RV commuter leaves. The equilibrium conditions during this period can be expressed as:
\begin{equation}
    \label{eq:S8_initialstage}
    \begin{cases}
        \frac{\partial C^\mathcal{R}(t)}{\partial t}=(\alpha+\pi-\beta)\big(\dot{w}_\mathcal{H}(t)+\dot{w}_\mathcal{C}^\mathcal{R}(t+w_\mathcal{H}(t))(1+\dot{w}_\mathcal{H}(t))\big)-\beta=0,\\
        \dot{A}_\mathcal{C}^\mathcal{R}(t)=s_\mathcal{H}.
    \end{cases}
\end{equation}

During the co-departure stage, both RV and PV commuters depart from home, and the highway bottleneck remains congested. The PV queue forms at the curbside bottleneck once the first PV arrives there, which is influenced by the RV queue. The equilibrium conditions during this stage can be expressed as:
\begin{equation}
    \label{eq:S8_codeparturestage}
\begin{cases}
        \frac{\partial C^\mathcal{R}(t)}{\partial t}=(\alpha+\pi-\beta)\big(\dot{w}_\mathcal{H}(t)+\dot{w}_\mathcal{C}^\mathcal{R}(t+w_\mathcal{H}(t))(1+\dot{w}_\mathcal{H}(t))\big)-\beta=0,\\[1ex]
        \frac{\partial C^\mathcal{P}(t)}{\partial t}=(\alpha-\beta)\big(\dot{w}_\mathcal{H}(t)+\dot{w}_\mathcal{C}^\mathcal{P}(t+w_\mathcal{H}(t))(1+\dot{w}_\mathcal{H}(t))\big)-\beta=0,\\[1ex]
        \dot{A}_\mathcal{C}^\mathcal{R}(t)+\dot{A}_\mathcal{C}^\mathcal{P}(t)=s_\mathcal{H}.
\end{cases}
\end{equation}

During the final stage, both the curbside bottleneck and the highway bottleneck remain congested. If $t_1^\mathcal{P}\leq t_1^\mathcal{R}$, only RV commuters leave home during this period; otherwise, only PV commuters leave home. Therefore, the equilibrium conditions for this stage can be expressed as:
\begin{itemize}
\item if $t_1^\mathcal{P}\leq t_1^\mathcal{R}$,
    \begin{equation}
    \label{eq:S8_finalphase_t1c<t1r}
    \begin{cases}
        \frac{\partial C^\mathcal{R}(t)}{\partial t}=(\alpha+\pi-\beta)\big(\dot{w}_\mathcal{H}(t)+\dot{w}_\mathcal{C}^\mathcal{R}(t+w_\mathcal{H}(t))(1+\dot{w}_\mathcal{H}(t))\big)-\beta=0,\\            \dot{A}_\mathcal{C}^\mathcal{R}(t)=s_\mathcal{H}.
    \end{cases}
    \end{equation}
\item if $t_1^\mathcal{P}>t_1^\mathcal{R}$,
    \begin{equation}
    \label{eq:S8_finalphase_t1c>t1r}    
    \begin{cases}
        \frac{\partial C^\mathcal{P}(t)}{\partial t}=(\alpha-\beta)\big(\dot{w}_\mathcal{H}(t)+\dot{w}_\mathcal{C}^\mathcal{P}(t+w_\mathcal{H}(t))(1+\dot{w}_\mathcal{H}(t))\big)-\beta=0,\\
        \dot{A}_\mathcal{C}^\mathcal{P}(t)=s_\mathcal{H}.
    \end{cases}
    \end{equation}
\end{itemize}

Based on Eqs. \eqref{eq:S8_initialstage}-\eqref{eq:S8_finalphase_t1c>t1r}, we can derive the equilibrium arrival rates of commuters at bottlenecks as:
\begin{itemize}
    \item if $t_1^\mathcal{P}\leq t_1^\mathcal{R}$,
    \begin{equation}
    \label{eq:arrivalrate_S8_t1c<t1r}
    \begin{aligned}
        &\dot{A}_\mathcal{H}(t) = \begin{cases}
            \frac{(\alpha+\pi)s_\mathcal{C}^\mathcal{R}}{\alpha+\pi-\beta},  & t\in[t_0^\mathcal{R},t_0^\mathcal{P})\cup(t_1^\mathcal{P},t_1^\mathcal{R}], \\
            \frac{\alpha s_\mathcal{C}^\mathcal{P}}{\alpha-\beta}+\frac{(1-\delta^\mathcal{R})(\alpha+\pi)s_\mathcal{C}^\mathcal{R}}{\alpha+\pi-\beta}, & t\in[t_0^\mathcal{P},t_1^\mathcal{P}].
        \end{cases}\\
        &\dot{A}_\mathcal{C}^\mathcal{R}(t) = \begin{cases}
            s_\mathcal{H}, & t\in\big[t_0^\mathcal{R},t_0^\mathcal{P}+w_\mathcal{H}(t_0^\mathcal{P})\big)
                \cup\big(t_1^\mathcal{P}+w_\mathcal{H}(t_1^\mathcal{P}),t_1^\mathcal{R}+w_\mathcal{H}(t_1^\mathcal{R})\big], \\
            \frac{(\alpha+\pi)(\alpha-\beta)s_\mathcal{C}^\mathcal{R}s_\mathcal{H}}
                {\alpha(\alpha+\pi-\beta)s_\mathcal{C}^\mathcal{P}+(1-\delta^\mathcal{R})(\alpha+\pi)(\alpha-\beta)s_\mathcal{C}^\mathcal{R}}, & t\in\big[t_0^\mathcal{P}+w_\mathcal{H}(t_0^\mathcal{P}),t_1^\mathcal{P}+w_\mathcal{H}(t_1^\mathcal{P})\big].\\
        \end{cases}\\
        &\dot{A}_\mathcal{C}^\mathcal{P}(t)=s_\mathcal{H}-\frac{(\alpha+\pi)(\alpha-\beta)s_\mathcal{C}^\mathcal{R}s_\mathcal{H}}{\alpha(\alpha+\pi-\beta)s_\mathcal{C}^\mathcal{P}+(1-\delta^\mathcal{R})(\alpha+\pi)(\alpha-\beta)s_\mathcal{C}^\mathcal{R}},t\in\big[t_0^\mathcal{P}+w_\mathcal{H}(t_0^\mathcal{P}),t_1^\mathcal{P}+w_\mathcal{H}(t_1^\mathcal{P})\big].
    \end{aligned}
    \end{equation}

    \item if $t_1^\mathcal{P}>t_1^\mathcal{R}$,
    \begin{equation}
    \label{eq:arrivalrate_S8_t1c>t1r}
    \begin{aligned}
        &\dot{A}_\mathcal{H}(t) = \begin{cases}
            \frac{(\alpha+\pi)s_\mathcal{C}^\mathcal{R}}{\alpha+\pi-\beta}, & t\in[t_0^\mathcal{R},t_0^\mathcal{P}),\\
            
            \frac{\alpha s_\mathcal{C}^\mathcal{P}}{\alpha-\beta}+\frac{(1-\delta^\mathcal{R})(\alpha+\pi)s_\mathcal{C}^\mathcal{R}}{(\alpha+\pi-\beta)}, & t\in[t_0^\mathcal{P},t_1^\mathcal{R}),\\
            
            
            \frac{\alpha s_\mathcal{C}^\mathcal{P}}{\alpha-\beta}, & t\in[t_1^\mathcal{R},t_1^\mathcal{P}].            
        \end{cases}\\
        &\dot{A}_\mathcal{C}^\mathcal{R}(t) = \begin{cases}
            s_\mathcal{H}, & t\in[t_0^\mathcal{R},t_0^\mathcal{P}+w_\mathcal{H}(t_0^\mathcal{P})),\\
            
            \frac{(\alpha+\pi)(\alpha-\beta)s_\mathcal{C}^\mathcal{R}s_\mathcal{H}}{\alpha(\alpha+\pi-\beta)s_\mathcal{C}^\mathcal{P}+(1-\delta^\mathcal{R})(\alpha+\pi)(\alpha-\beta)s_\mathcal{C}^\mathcal{R}}, & t\in\big[t_0^\mathcal{P}+w_\mathcal{H}(t_0^\mathcal{P}),t_1^\mathcal{R}+w_\mathcal{H}(t_1^\mathcal{R})\big].\\
        \end{cases}\\
        &\dot{A}_\mathcal{C}^\mathcal{P}(t)=\begin{cases}

            s_\mathcal{H}-\frac{(\alpha+\pi)(\alpha-\beta)s_\mathcal{C}^\mathcal{R}s_\mathcal{H}}{\alpha(\alpha+\pi-\beta)s_\mathcal{C}^\mathcal{P}+(1-\delta^\mathcal{R})(\alpha+\pi)(\alpha-\beta)s_\mathcal{C}^\mathcal{R}}, & t\in\big[t_0^\mathcal{P}+w_\mathcal{H}(t_0^\mathcal{P}),t_1^\mathcal{R}+w_\mathcal{H}(t_1^\mathcal{R})\big),\\
            
            s_\mathcal{H}, & t\in\big[t_1^\mathcal{R}+w_\mathcal{H}(t_1^\mathcal{R}),t_1^\mathcal{P}+w_\mathcal{H}(t_1^\mathcal{P})\big].
        \end{cases}
    \end{aligned}
 \end{equation}
\end{itemize}

At user equilibrium, all commuters have the same generalized travel cost, i.e., $C^\mathcal{R}(t_0^\mathcal{R})=C^\mathcal{R}(t_1^\mathcal{R})=C^\mathcal{P}(t_0^\mathcal{P})=C^\mathcal{P}(t_1^\mathcal{P})$. All RV commuters arrive at CBD within $[t_0^\mathcal{R},t^*]$ at a rate of $s_\mathcal{C}^\mathcal{R}$, all PV commuters arrive at the curbside bottleneck within $[t_0^\mathcal{P}+w_\mathcal{H}(t_0^\mathcal{P}),t_1^\mathcal{P}+w_\mathcal{H}(t_1^\mathcal{P})]$. Therefore, we have the following system of equations:
\begin{equation}
\label{eq:linearsystem_S8_UE}
        \begin{cases}
            \beta(t^*-t_0^\mathcal{R})+c^\mathcal{R}=(\alpha+\pi)(t^*-t_1^\mathcal{R})+c^\mathcal{R},\\
            (\alpha-\beta)w_\mathcal{H}(t_0^\mathcal{P})+\beta(t^*-t_0^\mathcal{P})+u^\mathcal{P}=\alpha(t^*-t_1^\mathcal{P})+u^\mathcal{P},\\
            \beta(t^*-t_0^\mathcal{R})+c^\mathcal{R}=(\alpha-\beta)w_\mathcal{H}(t_0^\mathcal{P})+\beta(t^*-t_0^\mathcal{P})+u^\mathcal{P},\\
            (t^*-t_0^\mathcal{R})s_\mathcal{C}^\mathcal{R}=N^\mathcal{R},\\[1.2ex]
            \int_{{t}_0^\mathcal{P}+w_\mathcal{H}(t_0^\mathcal{P})}^{t_1^\mathcal{P}+w_\mathcal{H}(t_1^\mathcal{P})} \dot{A}_\mathcal{C}^\mathcal{P}(t) dt=N-N^\mathcal{R}.
        \end{cases}
\end{equation}

By substituting Eqs. \eqref{eq:arrivalrate_S8_t1c<t1r} and \eqref{eq:arrivalrate_S8_t1c>t1r} into Eq. \eqref{eq:linearsystem_S8_UE}, we can derive the equilibrium results as follows:
\begin{itemize}
    \item if $t_1^\mathcal{P}\leq t_1^\mathcal{R}$,
    \begin{equation}\label{eq:ue_results_S8_1}
        \begin{cases}
            t_0^\mathcal{R}=t^*-\frac{N^\mathcal{R}}{s_\mathcal{C}^\mathcal{R}},\\
            t_1^\mathcal{R}=t^*-\frac{\beta}{\alpha+\pi}\frac{N^\mathcal{R}}{s_\mathcal{C}^\mathcal{R}},\\
            
            t_0^\mathcal{P} = t^*-\frac{Ns_\mathcal{H}(\alpha+\pi-\beta)}{s_\mathcal{C}^\mathcal{R}\widetilde{s}_\mathcal{C}^\mathcal{P}(\alpha+\pi)}+\frac{N^\mathcal{R}\big(s_\mathcal{H}(s_\mathcal{C}^\mathcal{R}+\widetilde{s}_\mathcal{C}^\mathcal{P})\frac{\alpha+\pi-\beta}{\alpha+\pi} - s_\mathcal{C}^\mathcal{R}\widetilde{s}_\mathcal{C}^\mathcal{P}  \big)}{\widetilde{s}_\mathcal{C}^\mathcal{P}(s_\mathcal{C}^\mathcal{R})^2},\\

            t_1^\mathcal{P} = t_0^\mathcal{P}+\frac{\alpha-\beta}{\alpha}\frac{N-N^\mathcal{R}}{\widetilde{s}_\mathcal{C}^\mathcal{P}},\\
            
            N^\mathcal{R}=\bigg(\frac{N}{s_\mathcal{C}^\mathcal{R}+\widetilde{s}_\mathcal{C}^\mathcal{P}} + \frac{\frac{(u^\mathcal{P}-c^\mathcal{R})(\alpha+\pi)}{\alpha(\alpha+\pi-\beta)}s_\mathcal{C}^\mathcal{R}\widetilde{s}_\mathcal{C}^\mathcal{P}}{(s_\mathcal{C}^\mathcal{R}+\widetilde{s}_\mathcal{C}^\mathcal{P}) \big(s_\mathcal{H}-s_\mathcal{C}^\mathcal{R}\frac{(\alpha+\pi)(\alpha-\beta)}{\alpha(\alpha+\pi-\beta)} \big)}\bigg)s_\mathcal{C}^\mathcal{R}.\\
        \end{cases}
    \end{equation} where the parameter $\widetilde{s}_\mathcal{C}^\mathcal{P}$ represents the actual service rate of the curbside bottleneck serving PVs during the co-departure stage, and $\widetilde{s}_\mathcal{C}^\mathcal{P}=\frac{\dot{A}_\mathcal{C}^\mathcal{P}(t)}{\dot{A}_\mathcal{C}^\mathcal{P}(t)+\delta^\mathcal{R}\dot{A}_\mathcal{C}^{\mathcal{R}}(t)}s_\mathcal{C}^\mathcal{P} =s_\mathcal{C}^\mathcal{P}-s_\mathcal{C}^\mathcal{R}\frac{\delta^\mathcal{R}(\alpha+\pi)(\alpha-\beta)}{\alpha(\alpha+\pi-\beta)}$ holds.
    
    \item if $t_1^\mathcal{P}>t_1^\mathcal{R}$,
    \begin{equation}\label{eq:ue_results_S8_2}
        \begin{cases}
            t_0^\mathcal{R}=t^*-\frac{N^\mathcal{R}}{s_\mathcal{C}^\mathcal{R}},\\
            t_1^\mathcal{R}=t^*-\frac{\beta}{\alpha+\pi}\frac{N^\mathcal{R}}{s_\mathcal{C}^\mathcal{R}},\\
            
            t_0^\mathcal{P}=t^*-\frac{Ns_\mathcal{H}(\alpha+\pi-\beta)}{s_\mathcal{C}^\mathcal{R}(s_\mathcal{C}^\mathcal{P}-\delta^\mathcal{R}s_\mathcal{H})(\alpha+\pi)}-\frac{N^\mathcal{R}\big( s_\mathcal{C}^\mathcal{R}\big(s_\mathcal{C}^\mathcal{P}-s_\mathcal{H}\frac{\alpha+\pi-\beta(1-\delta^\mathcal{R})}{\alpha+\pi}\big) - s_\mathcal{H}s_\mathcal{C}^\mathcal{P}\frac{\alpha+\pi-\beta}{\alpha+\pi}  \big)}{(s_\mathcal{C}^\mathcal{R})^2(s_\mathcal{C}^\mathcal{P}-\delta^\mathcal{R}s_\mathcal{H})},\\
            
            t_1^\mathcal{P}=t_0^\mathcal{P}+\frac{N(\alpha-\beta)}{\alpha(s_\mathcal{C}^\mathcal{P}-\delta^\mathcal{R}s_\mathcal{H})}-\frac{N^\mathcal{R}\big(s_\mathcal{H}\delta^\mathcal{R}+s_\mathcal{C}^\mathcal{R}(1-\delta^\mathcal{R})   \big)(\alpha-\beta)}{\alpha s_\mathcal{C}^\mathcal{R}(s_\mathcal{C}^\mathcal{P}-\delta^\mathcal{R}s_\mathcal{H})},\\
            
            N^\mathcal{R}=\bigg(\frac{N}{s_\mathcal{C}^\mathcal{P}+s_\mathcal{C}^\mathcal{R}(1-\delta^\mathcal{R})}+\frac{\frac{(u^\mathcal{P}-c^\mathcal{R})(\alpha+\pi)}{\alpha(\alpha+\pi-\beta)}s_\mathcal{C}^\mathcal{R}(s_\mathcal{C}^\mathcal{P}-\delta^\mathcal{R}s_\mathcal{H})}{(s_\mathcal{C}^\mathcal{P}+s_\mathcal{C}^\mathcal{R}(1-\delta^\mathcal{R}))\big(s_\mathcal{H}-s_\mathcal{C}^\mathcal{R}\frac{(\alpha+\pi)(\alpha-\beta)}{\alpha(\alpha+\pi-\beta)}\big)}\bigg)s_\mathcal{C}^\mathcal{R}.\\
        \end{cases}
    \end{equation}
\end{itemize}

In Scenarios 8, congestion occurs at the highway and the curbside bottlenecks. According to Proposition \ref{Proposition 1}, we have,
\begin{equation}
\label{eq:AppendixB_conditions_S8_1}
        \frac{s_\mathcal{C}^\mathcal{R}}{s_\mathcal{H} }>\frac{\alpha+\pi-\beta}{\alpha+\pi}.
\end{equation}

Since the departure intervals of RV and PV commuters overlap, implying that the condition $t_0^\mathcal{P}< t_1^\mathcal{R}$ holds. Combining it with Eqs. \eqref{eq:ue_results_S8_1} and \eqref{eq:ue_results_S8_2}, one can obtain,
\begin{equation}
\label{eq:AppendixB_conditions_S8_2}
    (u^\mathcal{P}-c^\mathcal{R})\in\bigg(0,\frac{N\big( \alpha(\alpha+\pi-\beta)s_\mathcal{H}-(\alpha+\pi)(\alpha-\beta)s_\mathcal{C}^\mathcal{R}\big)}{(\alpha+\pi)\big(s_\mathcal{H}(s_\mathcal{C}^\mathcal{R}+s_\mathcal{C}^\mathcal{P})-s_\mathcal{C}^\mathcal{R}s_\mathcal{C}^\mathcal{P}\big)}\bigg).
\end{equation}

Therefore, Eqs. \eqref{eq:AppendixB_conditions_S8_1} and \eqref{eq:AppendixB_conditions_S8_2} give the occurrence conditions for Scenario 8.

\section{Proof of Proposition 2}
\setcounter{equation}{0}

\renewcommand\theequation{B.\arabic{equation}}
\label{Appendix B}
Since the fixed cost of RV is lower than that of PV (i.e., $u^\mathcal{P}>c^\mathcal{R}$), if only one mode of transportation is used, it must necessarily be RV, i.e., $N^\mathcal{R}=N$, we thus have
\begin{equation}\label{eq:ts-t0R=N/scr}
    t^*-t_0^\mathcal{R}=\frac{N}{s_C^\mathcal{R}}.
\end{equation}

In addition, the condition $C^\mathcal{P}(t)\geq C^\mathcal{R}(t)$ holds for $\forall t\in(t_0^\mathcal{R},t^*]$, meaning that the travel cost of PV is always higher than that of RV for any departure time $t\in(t_0^\mathcal{R},t^*]$, we then have,
\begin{equation}
    \label{eq:A.4}
    (\alpha-\beta)w_\mathcal{H}(t)+\beta(t^*-t)+u^\mathcal{P}\geq\beta(t^*-t_0^\mathcal{R})+c^\mathcal{R},
\end{equation} where $w_\mathcal{H}(t)$ denotes the highway bottleneck queuing time for PV commuters who depart at time $t$. According to Proposition \ref{Proposition 1}, $w_\mathcal{H}(t)$ can be expressed as 
\begin{equation}
    w_\mathcal{H}(t)=\begin{cases}
        \frac{\frac{\alpha+\pi}{\alpha+\pi-\beta}s_\mathcal{C}^\mathcal{R}-s_\mathcal{H}}{s_\mathcal{H}}(t-t_0^\mathcal{R}),\ &\text{if}\ \frac{s_\mathcal{C}^\mathcal{R}}{s_\mathcal{H}}>\frac{\alpha+\pi-\beta}{\alpha+\pi}  \ \text{and}\ t_0^\mathcal{P}\leq t_1^\mathcal{R},\\
        
        \frac{\frac{\alpha+\pi}{\alpha+\pi-\beta}s_\mathcal{C}^\mathcal{R}-s_\mathcal{H}}{s_\mathcal{H}}(t_1^\mathcal{R}-t_0^\mathcal{R})-(t-t_1^\mathcal{R}),\ &\text{if}\ \frac{s_\mathcal{C}^\mathcal{R}}{s_\mathcal{H}}>\frac{\alpha+\pi-\beta}{\alpha+\pi}  \ \text{and}\ t_0^\mathcal{P}\in(t_1^\mathcal{R},t_1^\mathcal{R}+w_\mathcal{H}(t_1^\mathcal{R})),\\

        0, &\text{otherwise}.
    \end{cases}
\end{equation}

If only RV is used during the morning peak, Eq. \eqref{eq:A.4} should hold for any $t\in(t_0^\mathcal{R},t^*]$, thus we have,
\begin{equation}
\label{eq:A.5}
    u^\mathcal{P}-c^\mathcal{R}\geq\beta(t^*-t_0^\mathcal{R}).
\end{equation}
By substituting Eq. \eqref{eq:ts-t0R=N/scr} into Eq. \eqref{eq:A.5}, one can obtain,
\begin{equation}\label{eq:A.6}
    u^\mathcal{P}-c^\mathcal{R}\geq\frac{\beta N}{s_\mathcal{C}^\mathcal{R}}.
\end{equation}

Eq. \eqref{eq:A.6} presents the critical condition under which only RV will be used. If Eq. \eqref{eq:A.6} is not satisfied, both the RV and PV modes will be utilized during the morning peak.

(ii) When both modes are in use during the morning peak, the inter-equilibrium conditions of Eq. \eqref{eq:jointUE} should hold. Since the earliest departure time of RV commuters is always earlier than that of PV commuters, i.e., $t_0^\mathcal{R}<t_0^\mathcal{P}$, whether their departure intervals overlap depends on the departure times of the last RV commuter and the first PV commuter. When the first PV commuter departs earlier than the last RV commuter, i.e., $t_0^\mathcal{P} < t_1^\mathcal{R}$, their departure intervals overlap; otherwise, their departure intervals do not overlap. 

To identify the critical condition between the above two cases, we consider the situation of $t_0^\mathcal{P} = t_1^\mathcal{R}$, implying that the first PV commuter and the last RV commuter depart simultaneously. Combining $t_0^\mathcal{P} = t_1^\mathcal{R}$ with the equilibrium results of Scenarios 3, 5, and 8, one can obtain
\begin{equation}
\label{eq:A.7}
    \begin{cases}
        (u^\mathcal{P}-c^\mathcal{R})=\frac{N\beta(\alpha+\pi-\beta)}{\beta s_\mathcal{C}^\mathcal{P}+(\alpha+\pi)s_\mathcal{C}^\mathcal{R}}, & \text{if}\  \frac{s_\mathcal{C}^\mathcal{R}}{s_\mathcal{H} }\leq \frac{\alpha+\pi-\beta}{\alpha+\pi}, \\[1.5ex]

        (u^\mathcal{P}-c^\mathcal{R})=\frac{N\big( \alpha(\alpha+\pi-\beta)s_\mathcal{H}-(\alpha+\pi)(\alpha-\beta)s_\mathcal{C}^\mathcal{R}\big)}{(\alpha+\pi)\big(s_\mathcal{H}(s_\mathcal{C}^\mathcal{R}+s_\mathcal{C}^\mathcal{P})-s_\mathcal{C}^\mathcal{R}s_\mathcal{C}^\mathcal{P}\big)}, & \text{if}\ \frac{s_\mathcal{C}^\mathcal{R}}{s_\mathcal{H} }> \frac{\alpha+\pi-\beta}{\alpha+\pi}.
    \end{cases}
\end{equation}

Eq. \eqref{eq:A.7} provides the condition under which the last RV commuter and the first PV commuter depart simultaneously. Based on this, when  $u^\mathcal{P} - c^\mathcal{R}$ exceeds the threshold, PV commuters will delay their departure, resulting in non-overlapping departure intervals. Therefore, the condition for non-overlapping departure intervals of RV and PV commuters is as follows,
\begin{equation}\label{eq:A.8}
    \begin{cases}
        (u^\mathcal{P}-c^\mathcal{R})\geq\frac{N\beta(\alpha+\pi-\beta)}{\beta s_\mathcal{C}^\mathcal{P}+(\alpha+\pi)s_\mathcal{C}^\mathcal{R}}, & \text{if}\  \frac{s_\mathcal{C}^\mathcal{R}}{s_\mathcal{H} }\leq \frac{\alpha+\pi-\beta}{\alpha+\pi}, \\[1.5ex]

        (u^\mathcal{P}-c^\mathcal{R})\geq\frac{N\big( \alpha(\alpha+\pi-\beta)s_\mathcal{H}-(\alpha+\pi)(\alpha-\beta)s_\mathcal{C}^\mathcal{R}\big)}{(\alpha+\pi)\big(s_\mathcal{H}(s_\mathcal{C}^\mathcal{R}+s_\mathcal{C}^\mathcal{P})-s_\mathcal{C}^\mathcal{R}s_\mathcal{C}^\mathcal{P}\big)}, & \text{if}\ \frac{s_\mathcal{C}^\mathcal{R}}{s_\mathcal{H} }> \frac{\alpha+\pi-\beta}{\alpha+\pi}.
    \end{cases}
\end{equation}

This completes the proof of Proposition \ref{Proposition 2}.

\section{Derivation of typical scenarios with bidirectional congestion spillover}
\setcounter{equation}{0}
\renewcommand\theequation{C.\arabic{equation}}
\label{Appendix C}

For Scenario 3, the following system of equations ensures the user equilibrium:
\begin{equation}\label{eq:bieffect_S3_UE_EquationSystem}
    \begin{cases}
        t_0^\mathcal{P}-t_0^\mathcal{R}=\frac{u^\mathcal{P}-c^\mathcal{R}}{\beta},\\
        (\alpha+\pi)(t^*-t_1^\mathcal{R})+c^\mathcal{R} = \beta(t^*-t_0^\mathcal{R})+c^\mathcal{R},\\
        \alpha(t^*-t_1^\mathcal{P})+u^\mathcal{P} = \beta(t^*-t_0^\mathcal{P})+u^\mathcal{P},\\
        \int_{t_0^\mathcal{R}}^{t_1^\mathcal{R}} \dot{A}_\mathcal{C}^\mathcal{R}(t)dt=N^\mathcal{R},\\
        \int_{t_0^\mathcal{P}}^{t_1^\mathcal{P}} \dot{A}_\mathcal{C}^\mathcal{P}(t)dt=N-N^\mathcal{R}.
    \end{cases}
\end{equation} where the first line in Eq. \eqref{eq:bieffect_S3_UE_EquationSystem} represents the inter-equilibrium condition. The second line indicates that the first and the last RV commuters have an equal generalized travel cost, and the third line applies the same condition to PV commuters. The last two lines are population constraints.

Substituting the equilibrium departure rates of commuters, i.e., Eqs. \eqref{eq:bieffect_DepartRate_S3} and \eqref{eq:bieffect_CoDepartRate_S3}, into Eq. \eqref{eq:bieffect_S3_UE_EquationSystem}, one can derive the equilibrium results for Scenario 3 as follows,
\begin{equation}\label{eq:bieffect_S3_UE_results}
    \begin{cases}
        t_0^\mathcal{R} = t^*-\frac{(\alpha+\pi)\big(N^\mathcal{R}\beta-(u^\mathcal{P}-c^\mathcal{R})(r_1-r_2^\mathcal{R})\big)}{\beta(\alpha+\pi-\beta)r_2^\mathcal{R}},\\[1.2ex]
        
        t_1^\mathcal{R} = \frac{\alpha+\pi-\beta}{\alpha+\pi}t^*+\frac{\beta}{\alpha+\pi}t_0^\mathcal{R},\\[1.2ex]
        
        t_0^\mathcal{P} = t_0^\mathcal{R}+\frac{u^\mathcal{P}-c^\mathcal{R}}{\beta},\\[1.2ex]
        
        t_1^\mathcal{P} = \frac{\alpha-\beta}{\alpha}t^*+\frac{\beta}{\alpha}t_0^\mathcal{P},\\[1.2ex]
        
        N^\mathcal{R} = \frac{\alpha\beta(\alpha+\pi-\beta)Nr_2^\mathcal{R}-(u^\mathcal{P}-c^\mathcal{R})\big(r_1\big(\beta\pi r_3-\alpha(\alpha+\pi-\beta)r_2^\mathcal{P}\big) +r_3r_2^\mathcal{R}\beta(\alpha-\beta)  \big)}{\beta\big(\alpha(\alpha+\pi-\beta)(r_2^\mathcal{P}+r_2^\mathcal{R})-\beta\pi r_3\big)}.
    \end{cases}
\end{equation} where the parameters $r_1$ and $r_3$ denote the departure rates of commuters from home in the initial and final stages, respectively, as defined by Eq. \eqref{eq:bieffect_DepartRate_S3}.
$r_2^\mathcal{R}$ and $r_2^\mathcal{P}$ denote the departure rates of RV and PV commuters during the co-departure stage, respectively, as defined by Eq. \eqref{eq:bieffect_CoDepartRate_S3}.

For Scenario 5, the following system of equations ensures the user equilibrium:
\begin{equation}\label{eq:bieffect_S5_UE_EquationSystem}
    \begin{cases}
        t_0^\mathcal{P}-t_0^\mathcal{R}=\frac{u^\mathcal{P}-c^\mathcal{R}}{\beta},\\[1.2ex]
        
        (\alpha+\pi)(t^*-t_1^\mathcal{R})+c^\mathcal{R} = \beta(t^*-t_0^\mathcal{R})+c^\mathcal{R},\\[1.2ex]
        
        \alpha(t^*-t_1^\mathcal{P})+u^\mathcal{P} = \beta(t^*-t_0^\mathcal{P})+u^\mathcal{P},\\[1.2ex]
        
        \int_{t_0^\mathcal{R}}^{t_1^\mathcal{R}+w_\mathcal{H}(t_1^\mathcal{R})} \dot{A}_\mathcal{C}^\mathcal{R}(t)dt=N^\mathcal{R},\\[1.2ex]
        
        \int_{t_0^\mathcal{P}}^{t_1^\mathcal{P}+w_\mathcal{H}(t_1^\mathcal{P})} \dot{A}_\mathcal{C}^\mathcal{P}(t)dt=N-N^\mathcal{R}.
    \end{cases}
\end{equation}where the first line in Eq. \eqref{eq:bieffect_S5_UE_EquationSystem} represents the inter-equilibrium condition. The second line indicates that the first and the last RV commuters have an equal generalized travel cost, and the third line applies the same condition to PV commuters. The last two lines are population constraints.

Substituting the equilibrium departure rates of commuters, i.e., Eqs. \eqref{eq:bieffect_DepartRate_S5}, into Eq. \eqref{eq:bieffect_S5_UE_EquationSystem}, one can derive the equilibrium results for Scenario 5 as follows,
\begin{equation}\label{eq:bieffect_S5_UE_results}
    \begin{cases}
        t_0^\mathcal{R} = t^*-\frac{N^\mathcal{R}\beta s_\mathcal{H}+(u^\mathcal{P}-c^\mathcal{R})(r_2r_2^\mathcal{R}-r_1s_\mathcal{H})}{\frac{\beta(\alpha+\pi-\beta)}{\alpha+\pi}r_2r_2^\mathcal{R}},\\[1.2ex]
        
        t_1^\mathcal{R} = \frac{\alpha+\pi-\beta}{\alpha+\pi}t^*+\frac{\beta}{\alpha+\pi}t_0^\mathcal{R},\\[1.2ex]
        
        t_0^\mathcal{P} = t_0^\mathcal{R}+\frac{u^\mathcal{P}-c^\mathcal{R}}{\beta},\\[1.2ex]
        
        t_1^\mathcal{P} = \frac{\alpha-\beta}{\alpha}t^*+\frac{\beta}{\alpha}t_0^\mathcal{P},\\[1.2ex]
        
        N^\mathcal{R} = \frac{\alpha\beta(\alpha+\pi-\beta)Nr_2r_2^\mathcal{R}-(u^\mathcal{P}-c^\mathcal{R})\big( r_1s_\mathcal{H}\big(\beta\pi r_3-\alpha(\alpha+\pi-\beta)r_2\big)+r_2r_2^\mathcal{R}\big( \alpha(\alpha+\pi-\beta) r_1+\beta(\alpha-\beta)r_3 \big) \big)}{\beta s_\mathcal{H}\big( \alpha(\alpha+\pi-\beta)r_2-\beta\pi r_3  \big)}.
    \end{cases}
\end{equation} where the parameters $r_1$, $r_2$, and $r_3$ denote the departure rates of commuters from home in the initial, co-departure, and final stages, respectively, i.e, $\dot{A}_\mathcal{H}(t)$ defined by Eq. \eqref{eq:bieffect_DepartRate_S5}. $r_2^\mathcal{R}$ and $r_2^\mathcal{P}$ denote the arrival rates of RV and PV commuters at the curbside bottlenecks during the co-departure stage, i.e., $\dot{A}_\mathcal{C}^\mathcal{R}(t)$ and $\dot{A}_\mathcal{C}^\mathcal{P}(t)$ defined by Eq. \eqref{eq:bieffect_DepartRate_S5}.

For Scenario 8, the following system of equations ensures the user equilibrium:
\begin{equation}\label{eq:bieffect_S8_UE_EquationSystem}
    \begin{cases}
        t_0^\mathcal{P}-t_0^\mathcal{R}=
        \frac{u^\mathcal{P}-c^\mathcal{R}}{\alpha-\frac{(\alpha+\pi)(\alpha-\beta)s_\mathcal{C}^\mathcal{R}}{(\alpha+\pi-\beta)s_\mathcal{H}}},\\[1.2ex]
        
        \beta(t^*-t_0^\mathcal{R})+c^\mathcal{R}=(\alpha+\pi)(t^*-t_1^\mathcal{R})+c^\mathcal{R},\\[1.2ex]
        
        (\alpha-\beta)w_\mathcal{H}(t_0^\mathcal{P})+\beta(t^*-t_0^\mathcal{P})+u^\mathcal{P}=\alpha(t^*-t_1^\mathcal{P})+u^\mathcal{P},\\[1.2ex]
        
        \int_{{t}_0^\mathcal{R}}^{t_1^\mathcal{R}+w_\mathcal{H}(t_1^\mathcal{R})} \dot{A}_\mathcal{C}^\mathcal{R}(t) dt=N^\mathcal{R},\\[1.2ex]
        
        \int_{{t}_0^\mathcal{P}+w_\mathcal{H}(t_0^\mathcal{P})}^{t_1^\mathcal{P}+w_\mathcal{H}(t_1^\mathcal{P})} \dot{A}_\mathcal{C}^\mathcal{P}(t) dt=N-N^\mathcal{R}.
    \end{cases}
\end{equation} where the first three lines in Eq. \eqref{eq:bieffect_S8_UE_EquationSystem} represent the inter-equilibrium condition and the intra-equilibrium conditions for RV and PV commuters. The last two lines are population constraints.

Substituting the equilibrium departure rates of commuters, i.e., Eqs. \eqref{eq:bieffect_DepartRate_S8}, into Eq. \eqref{eq:bieffect_S8_UE_EquationSystem}, one can derive the equilibrium results for Scenario 8 as follows,
\begin{equation}\label{eq:bieffect_S8_UE_results}
    \begin{cases}
        t_0^\mathcal{R} = t^* - \frac{N^\mathcal{R}s_\mathcal{H}\big((\alpha - \beta)r_1-\alpha s_\mathcal{H}\big) + (u^\mathcal{P}-c^\mathcal{R})(r_1 s_\mathcal{H} - r_2 r_2^\mathcal{R})s_\mathcal{H}}{\frac{\alpha  + \pi- \beta}{\alpha + \pi}\big((\alpha - \beta)r_1-\alpha s_\mathcal{H}\big)r_2 r_2^\mathcal{R}},\\[1.2ex]
        
        t_1^\mathcal{R} = \frac{\alpha+\pi-\beta}{\alpha+\pi}t^*+\frac{\beta}{\alpha+\pi}t_0^\mathcal{R},\\[1.2ex]
        
        t_0^\mathcal{P} = t_0^\mathcal{R}+\frac{(u^\mathcal{P}-c^\mathcal{R})s_\mathcal{H}}{\alpha s_\mathcal{H}-(\alpha-\beta)r_1},\\[1.2ex]
        
        t_1^\mathcal{P} = \frac{\alpha-\beta}{\alpha}t^*+\frac{\beta}{\alpha}t_0^\mathcal{P}-\frac{(\alpha-\beta)(u^\mathcal{P}-c^\mathcal{R})(r_1-s_\mathcal{H})}{\alpha\big(\alpha s_\mathcal{H}-(\alpha-\beta)r_1\big)},\\[1.2ex]
        
        N^\mathcal{R} = \frac{\alpha(\alpha+\pi-\beta)Nr_2r_2^\mathcal{R}\big(\alpha s_\mathcal{H}-(\alpha-\beta)r_1\big) - (u^\mathcal{P}-c^\mathcal{R})\big( r_1(s_\mathcal{H})^2\big(\beta\pi r_3-\alpha(\alpha+\pi-\beta)r_2\big) +\alpha(\alpha+\pi)r_2r_2^\mathcal{R}\big(s_\mathcal{H}(s_\mathcal{C}^\mathcal{R}+s_\mathcal{C}^\mathcal{P})-s_\mathcal{C}^\mathcal{R}s_\mathcal{C}^\mathcal{P}\big)    \big)}{s_\mathcal{H}\big( \alpha s_\mathcal{H}-(\alpha-\beta)r_1 \big)\big(\alpha r_2-\frac{\beta\pi}{\alpha+\pi-\beta}r_3\big)}.
    \end{cases}
\end{equation} where the parameters $r_1$, $r_2$, and $r_3$ denote the departure rates of commuters from home in the initial, co-departure, and final stages, respectively, i.e, $\dot{A}_\mathcal{H}(t)$ defined by Eq. \eqref{eq:bieffect_DepartRate_S8}. $r_2^\mathcal{R}$ and $r_2^\mathcal{P}$ denote the arrival rates of RV and PV commuters at the curbside bottlenecks during the co-departure stage, i.e., $\dot{A}_\mathcal{C}^\mathcal{R}(t)$ and $\dot{A}_\mathcal{C}^\mathcal{P}(t)$ defined by Eq. \eqref{eq:bieffect_DepartRate_S8}.

\printcredits

\bibliographystyle{trbad}

\bibliography{cas-refs}

\end{document}